\documentclass[a4paper,twoside,12pt]{article} 
\usepackage{amsmath}
\usepackage{amssymb}
\usepackage{amsthm}
\usepackage{amsxtra}

\input xy
\xyoption{all}
\newdir{ >}{{}*!/-3mm/@{>}}
\newdir{ (}{{}*!/-3mm/@^{(}}

\theoremstyle{plain}
\newtheorem{proposition}{Proposition}[section]
\newtheorem{lemma}[proposition]{Lemma}
\newtheorem{corollary}[proposition]{Corollary}
\newtheorem{theorem}[proposition]{Theorem}

\theoremstyle{definition}
\newtheorem{notation}[proposition]{Notation}
\newtheorem{definition}[proposition]{Definition}
\newtheorem{definition and remark}[proposition]{Definition and Remark}
\newtheorem{remark}[proposition]{Remark}
\newtheorem{convention}[proposition]{Convention}
\newtheorem{convention and remark}[proposition]{Convention and Remark}
\newtheorem{notation and remark}[proposition]{Notation and Remark}
\newtheorem{notation, definition and remark}[proposition]{Notation, Definition and Remark}
\newtheorem{notation and definition}[proposition]{Notation and Definition}
\newtheorem{remark and definition}[proposition]{Remark and Definition}
\newtheorem{remark, notation and definition}[proposition]{Remark, Notation and Definition}
\newtheorem{remark and notation}[proposition]{Remark and Notation}
\newtheorem{reminder and remark}[proposition]{Reminder and Remark}
\newtheorem{remark and problems}[proposition]{Remark and Problems}
\newtheorem{example and remark}[proposition]{Example and Remark}
\newtheorem{example}[proposition]{Example}
\newtheorem{open question}[proposition]{Question}

\allowdisplaybreaks

\DeclareMathOperator{\Mor}{Mor}
\DeclareMathOperator{\Ass}{Ass}
\DeclareMathOperator{\Proj}{Proj}
\DeclareMathOperator{\Spec}{Spec}
\DeclareMathOperator{\Gl}{Gl}
\DeclareMathOperator{\Glf}{\underline{Gl}}
\DeclareMathOperator{\Hilb}{Hilb}
\DeclareMathOperator{\barHilb}{\bar{Hilb}}
\DeclareMathOperator{\Hilbf}{\underline{Hilb}}
\DeclareMathOperator{\barHilbf}{\underline{\bar{Hilb}}}
\newcommand{\Xf}{\underline{X}}
\DeclareMathOperator{\reg}{reg}
\DeclareMathOperator{\End}{end}
\renewcommand{\Im}{\imago} \DeclareMathOperator{\imago}{Im}
\DeclareMathOperator{\Ker}{Ker}
\DeclareMathOperator{\Char}{char}
\DeclareMathOperator{\id}{id}
\DeclareMathOperator{\satu}{sat} \newcommand{\sat}{{\satu}}
\DeclareMathOperator{\redu}{red} \newcommand{\red}{{\redu}}
\DeclareMathOperator{\Bin}{Bin}
\DeclareMathOperator{\supp}{supp}
\DeclareMathOperator{\Lt}{Lt}
\DeclareMathOperator{\init}{in}
\DeclareMathOperator{\Gin}{Gin}

\DeclareMathOperator{\gv}{gv}
\DeclareMathOperator{\hv}{hv}
\DeclareMathOperator{\lexu}{lex} \newcommand{\lex}{{\lexu}}
\DeclareMathOperator{\closed}{m}

\newcommand{\sk}{\medskip}
\renewcommand{\phi}{\varphi}
\renewcommand{\subsetneq}{\varsubsetneq}
\renewcommand{\iff}{\ \Longleftrightarrow\ }
\newcommand{\ra}{\rightarrow}
\newcommand{\Ra}{\Rightarrow}

\newcommand{\defiff}{:\Leftrightarrow}
\newcommand{\N}{\mathbb{N}}
\newcommand{\Z}{\mathbb{Z}}
\newcommand{\Q}{\mathbb{Q}}
\newcommand{\unipotent}{\mathcal{U}}
\newcommand{\Ideals}{\mathbb{I}}
\newcommand{\affline}{\mathbb{A}^1}
\newcommand{\afflinef}{\underline{\mathbb{A}}^1}
\newcommand{\etaf}{\overline{\eta}}
\newcommand{\terms}[1][]{\mathbb{T}_{#1}}
\newcommand{\ideala}{{\mathfrak{a}}}
\newcommand{\idealb}{{\mathfrak{b}}}
\newcommand{\idealc}{{\mathfrak{c}}}
\newcommand{\ideald}{{\mathfrak{d}}}
\newcommand{\ideall}{{\mathfrak{l}}}
\newcommand{\idealp}{{\mathfrak{p}}}
\newcommand{\idealq}{{\mathfrak{q}}}
\newcommand{\proj}[1][K]{{\mathbb{P}^n_{\!\!#1}}}
\newcommand{\alg}{{\mathcal{A}lg}}
\newcommand{\cat}{{\mathcal{C}}}
\newcommand{\sch}{{\mathcal{S}ch}}
\newcommand{\sets}{{\mathcal{S}et}}
\newcommand{\sheafo}{{\mathcal{O}}}
\newcommand{\sheaff}{{\mathcal{F}}}
\newcommand{\sheafg}{{\mathcal{G}}}
\newcommand{\sheafh}{{\mathcal{H}}}
\newcommand{\sheafhilb}{{\mathcal{H}ilb}}
\newcommand{\sheafi}{{\mathcal{I}}}
\newcommand{\sheafj}{{\mathcal{J}}}
\newcommand{\restr}{\!\!\upharpoonright}
\newcommand{\fbar}{{\boldsymbol{f}}}

\newcommand{\borelg}{>_\textup{\tiny{Bor}}}
\newcommand{\borelgeq}{\geq_\textup{\tiny{Bor}}}
\newcommand{\deglex}{\textup{hlex}}
\newcommand{\deglexg}{>_\textup{\tiny{hlex}}}
\newcommand{\deglexl}{<_\textup{\tiny{hlex}}}
\newcommand{\degrevlex}{\textup{rlex}}
\newcommand{\degrevlexg}{>_\textup{\tiny{rlex}}}

\begin{document}

\title{Connectedness of Hilbert scheme strata defined by bounding cohomology}

\author{Stefan Fumasoli}
\date{Z\"urich, \today}

\pagestyle{empty}
\begin{titlepage}

	\begin{center}
		\newcommand{\abgr}{\vspace{7ex}}
		\newcommand{\abkl}{\vspace{1.8ex}}

		\begin{bfseries}\begin{Large}
			Connectedness of Hilbert Scheme Strata Defined by Bounding Cohomology\\
		\end{Large}\end{bfseries}

		\abgr\abkl\abkl

		Dissertation\\
		\abkl
		zur\\
		\abkl
		Erlangung der naturwissenschaftlichen Doktorw\"urde\\
		(Dr.\ sc.\ nat.)\\
		\abkl
		vorgelegt der\\
		\abkl
		Mathematisch-naturwissenschaftlichen Fakult\"at\\
		\abkl
		der\\
		\abkl
		Universit\"at Z\"urich\\
		\abkl
		von\\
		\abkl
		Stefan Fumasoli\\
		\abkl
		von\\
		\abkl
		Z\"urich und Cadro TI\\

		\abgr

		Begutachtet von\\
		Prof.\ Dr.\ Markus Brodmann\\
		Prof.\ Dr.\ J\"urgen Herzog\\
		Prof.\ Dr.\ Bernd Sturmfels\\
		
		\abgr

		Z\"urich 2005
	\end{center}

\end{titlepage}
\ \vfill

\noindent
Die vorliegende Arbeit wurde von der Mathematisch-naturwissenschaftlichen Fakult\"at der Universit\"at Z\"urich auf Antrag von Prof.\ Dr.\ Markus Brodmann und Prof.\ Dr.\ Christian Okonek als Dissertation angenommen.
\cleardoublepage

\section*{Dank}
Als erstes danke ich meinem Doktorvater, Prof.\ Markus Brodmann, ganz herzlich f\"ur die hervorragende Betreuung und die immer wieder ermunternde Unterst\"utzung. 

Ein besonderer Dank geht an Enrico Sbarra, der mich bei einem Besuch in Z\"urich f\"ur das Rechnen von Beispielen begeisterte. Ihm verdanke ich den Hinweis auf den Begriff der Cohen-Macaulay-Filtrierung, der sich f\"ur diese Arbeit als \"ausserst wichtig erwiesen hat. Ferner m\"ochte ich mich bei Prof.\ Peter Schenzel f\"ur die anregenden Diskussionen bedanken.

Ich danke Prof.\ Christian Okonek sowie den beiden Gutachtern Prof.\ J\"urgen Herzog und Prof.\ Bernd Sturmfels, dass sie sich die Zeit genommen haben, meine Arbeit anzuschauen.

Sehr dankbar bin ich meinem ehemaligen B\"urokollegen Mihai-Sorin Stupariu f\"ur seine M\"uhe, die Arbeit nach Druck- und Englischfehlern abzusu\-chen. Er hat mir viele Verbesserungen vorgeschlagen.

\vspace{2cm}
\noindent
Stefan Fumasoli, Z\"urich im Fr\"uhling 2005
\cleardoublepage

\renewcommand{\abstractname}{Zusammenfassung}
\begin{abstract}
\frenchspacing
Sei $\Hilb^p_K$ das Hilbertschema, das die abgeschlossenen Un\-ter\-sche\-ma\-ta von $\proj$ mit Hil\-bert\-po\-ly\-nom $p\in\Q[t]$ \"uber einem K\"orper $K$ mit $\Char K\!=\!0$ para\-metri\-siert. Durch Beschr\"ankung der ko\-ho\-mo\-lo\-gi\-schen Hil\-bert\-funk\-ti\-onen der Punkte von $\Hilb^p_K$ nach unten werden lokal abgeschlossene Unterr\"aume des Hilbertschemas definiert. In dieser Arbeit wird bewiesen, dass einige dieser  Unterr\"aume zusammenh\"angend sind. Dazu wird die Theorie der Binomialideale, die von D.~Mall in \cite{M2000} untersucht worden sind, weiterentwickelt. Es stellt sich heraus, dass die von Mall konstruierten Binomial\-ideale Cohen-Macaulay-filtriert sind und dass f\"ur diese Ideale das Initialideal und das generische Initialideal bez\"uglich jeglicher zul\"assiger Termordnung \"uber\-ein\-stim\-men.
\nonfrenchspacing
\end{abstract}

\sk
\renewcommand{\abstractname}{Abstract}
\begin{abstract}
Let $\Hilb^p_K$ be the Hilbert scheme parametrizing the closed subschemes of $\proj$ with Hilbert polynomial $p\in\Q[t]$ over a field $K$ of characteristic zero. By bounding below the cohomological Hilbert functions of the points of $\Hilb^p_K$ we define locally closed subspaces of the Hilbert scheme. The aim of this thesis is to show that some of these subspaces are connected. For this we exploit the binomial ideals constructed by D.~Mall in \cite{M2000}. It turns out that these binomial ideals are sequentially Cohen-Macaulay and that their initial ideals and their generic initial ideals coincide for any admissible term order.
\end{abstract}
\clearpage
\tableofcontents
\cleardoublepage
\pagestyle{headings}
\section{Introduction}
Let $\proj$ be the projective space of dimension $n\geq 1$ over a field $K$ and let $\Hilb^p_K$ be the Hilbert scheme which parametrizes the closed subschemes of $\proj$ with Hilbert polynomial $p\in\Q[t]$, thus the quotients $\sheafo_{\proj}\twoheadrightarrow\nolinebreak\sheaff$ with Hilbert polynomial $p$, i.\,e.\ the ideal sheaves $\sheafi\subset\sheafo_{\proj}$ with Hilbert polynomial $q(t):=\binom{t+n}{n}-\nolinebreak p(t)$. R.~Hartshorne proved in 1963 in his thesis \cite{H1966} that $\Hilb^p_K$ is linearly connected. This means that for any two points of the Hilbert scheme there is a sequence of deformations defined over $\mathbb P^1_K$ connecting these points. Since then his techniques have been developed further and several connectedness results concerning interesting topological subspaces of $\Hilb^p_K$ have been proved: 

Let $\sheafj\subset\sheafo_{\proj\times_K\Hilb^p_K}$ be the universal sheaf of ideals with Hilbert polynomial $q$. For each point $x\in\Hilb^p_K$ define the function $h_x:\Z\ra\N$ to be the Hilbert function of the ideal sheaf $\sheafj^{(x)}:=\sheafj\otimes_K\kappa(x)\subset\sheafo_{\proj[\kappa(x)]}$. If $f,g:\Z\ra\N$ are two numerical functions we write $f\geq g$ if $f(j)\geq g(j)$ for all $j\in\Z$. Let $f:\Z\ra\N$ be a numerical function. By the Semicontinuity Theorem
$$H_{\geq f}:=\{x\in\Hilb^p_K\mid h_x\geq f\}$$ 
is a closed and 
$$H_f:=\{x\in\Hilb^p_K\mid h_x=f\}$$ 
is a locally closed subspace of $\Hilb^p_K$. G.~Gotzmann \cite{Go1988} showed in 1988 that $H_{\geq f}$ is connected if $K$ is of characteristic zero. K.~Pardue \cite{P} showed in 1996 that $H_f$ is connected if $K$ is an infinite field of any characteristic. 

The $K$-points of $\Hilb^p_K$ are precisely the saturated homogeneous ideals with Hilbert polynomial $q$ of the polynomial ring $S:=K[X_0,\dots,X_n]$. Two such points are said to be connected by a Gr\"obner deformation if one of them is the initial ideal or the generic initial ideal of the other. In this case, if one of these ideals is generated by monomials and the other by monomials and binomials, the Gr\"obner deformation is called binomial. In his Habilitations\-schrift of 1997, D. Mall \cite{M2000} gave an algorithmic proof of the fact that over an algebraically closed field $K$ of characteristic zero, not only $H_f$ and $H_{\geq f}$, but also 
$$H_{\leq f}:=\{x\in\Hilb^p_K\mid h_x\leq f\}$$ 
are connected by a sequence of Gr\"obner deformations, all of them binomial excepting the first and the last one. Recently I. Peeva and M. Stillman \cite{PS} proved that in characteristic zero the set of all homogeneous (not necessarily saturated) ideals with Hilbert function $f$ is connected by a sequence of Gr\"obner deformations, using some slightly different binomial ideals.

The aim of this thesis is to show that some subspaces of $H_{\geq f}$ and $H_f$ defined by bounding below the cohomological Hilbert functions of the points of $\Hilb^p_K$ are connected. More precisely, for $x\in\Hilb^p_K$ let $\sheaff^{(x)}:=\sheafo_{\proj[\kappa(x)]}/\sheafj^{(x)}$. For $i\in\N$ and $x\in\Hilb^p_K$ define the cohomological Hilbert functions 
$$h^i_{x}:\Z\ra\N,\ j\mapsto\dim_{\kappa(x)}H^i\bigl(\proj[\kappa(x)],\sheaff^{(x)}(j)\bigr),$$
$$\bar h^i_{x}:\Z\ra\N,\ j\mapsto\dim_{\kappa(x)}H^i\bigl(\proj[\kappa(x)],\sheafj^{(x)}(j)\bigr).$$ 
Let $\fbar=(f_i)_{i\in\N}$ be a sequence of numerical functions $f_i:\Z\ra\N$. Then by the Semicontinuity Theorem 
$$H^{\geq\fbar}:=\{x\in\Hilb^p_K\mid h^i_x\geq f_i\ \forall\,i\in\N\},$$ 
$$\bar H^{\geq\fbar}:=\{x\in\Hilb^p_K\mid \bar h^i_x\geq f_i\ \forall\,i\geq 1\}$$ 
are closed subspaces of $\Hilb^p_K$. 

In Theorem \ref{main Theorem} we show that the spaces
$$H_f\cap H^{\geq\fbar},\ \bar H^{\geq\fbar},\ H_f\cap\bar H^{\geq\fbar}\text{ and }H_{\geq f}\cap\bar H^{\geq\fbar}$$ 
are connected if $\Char K=0$.

\sk
\emph{Sketch of the proof:}
Let $X_K$ be one of these subsets of the Hilbert scheme. Since $X_K$ is locally closed in $\Hilb^p_K$, we may endow it with the induced reduced scheme structure. By flat base change and general nonsense it follows that $X_k\cong(X_K\times_K k)_\red$ for any field extension $K\subset k$. So it suffices to show the connectedness of the set $\closed(X_K)$ of closed points of $X_K$ in the case when $K$ is algebraically closed.

Let $\Char K=0$. We use the techniques of D.~Mall \cite{M1997, M2000}, J.~Herzog and E.~Sbarra \cite{HS, S} to show that $\closed(X_K)$ is connected by Gr\"obner deformations:

By the Serre-Grothendieck Correspondence the cohomology groups of sheaves correspond to local cohomology groups. So the problem is translated into the language of commutative algebra. Given two homogeneous saturated ideals $\ideala$, $\idealb\subset S$ with the same Hilbert polynomial, the algorithm of Mall yields a sequence of Borel ideals $(\idealc_i)_{i=0}^r$ and a sequence of binomial ideals $(\ideald_i)_{i=1}^r$ such that $\idealc_0^\sat=\Gin_{\degrevlex}(\ideala)$ and $\idealc_r^\sat=\Gin_{\degrevlex}(\idealb)$ with respect to the reverse lexicographic term order $\leq_{\degrevlex}$ and such that for each $i\in\{1,\dots,r\}$ the set of all possible initial ideals of $\ideald_i$ equals $\{\idealc_{i-1},\idealc_i\}$. 

One of our main tools is a Theorem of \cite{S} and \cite{HS}: For $i\in\N$ and a graded $S$-module $M$ define the function $h^i_M:\Z\ra\N$, $j\mapsto\dim_K H^i_{S_+}(M)_j$, where $H^i_{S_+}(M)_j$ denotes the $j$th homogeneous part of the $i$th local cohomology module of $M$ with respect to the irrelevant ideal $S_+$.  If $\ideala\subset S$ is a homogeneous ideal, then $h^i_{S/\ideala}\leq h^i_{S/\init_\tau\ideala}$ for all $i\in\N$ and for any term order $\tau$. Moreover, $h^i_{S/\ideala}=h^i_{S/\Gin_\degrevlex\ideala}$ for all $i\in\N$ if and only if $S/\ideala$ is sequentially Cohen-Macaulay (cf.\ Proposition \ref{main Theorem of Herzog-Sbarra}). Therefore, there are essentially two facts which have to be proven, namely: One of the two initial ideals of a Mall binomial ideal $\ideald$ is the generic initial ideal of $\ideald$ with respect to $\leq_{\degrevlex}$\linebreak[1] (cf.~Theorem~\ref{Gin = in for good binomial ideals}). Moreover, $S/\ideald$ is sequentially Cohen-Macaulay (cf.~Theorem~\ref{binomial ideals and sCM}).

The connectedness of $\closed(X_K)$ now follows from the fact that $\closed(X_K)$ is closed under isomorphisms. This means that for every two closed points $x$, $y\in\Hilb^p_K$ such that $\sheaff^{(x)}$ and $\sheaff^{(y)}$ are isomorphic as $\sheafo_{\proj}$-modules it holds $x\in X_K$ if and only if $y\in X_K$ (cf.~section~\ref{subsection connectedness}).
\clearpage
\section{Preliminaries}
In this section we introduce the needed combinatorial and algebraic notions and collect some results about initial and generic initial ideals. The latter ones play an important role in the proof of connectedness of Hilbert schemes: A homogeneous ideal $\ideala$ in a polynomial ring $K[X_1,\dots,X_n]$ over a field $K$ and its generic initial ideal $\Gin\ideala$ have the same Hilbert function. Furthermore, by means of weight orders, they are connected in the Hilbert scheme by a sequence of lines (cf.~section~\ref{subsection connectedness}).

Generic initial ideals are Borel-fixed, which means that they are fixed under the action of upper triangular matrices. If $K$ has characteristic zero, they correspond to Borel sets which have a combinatorial behaviour which is rather easy to understand (cf.~section~\ref{subsection Borel sets}).  

First we want to fix some basic notations concerning essentially polynomials and term orders. 
\subsection{Notations and definitions}
\begin{convention}
Let $\N$ denote the set of nonnegative integers. Throughout this thesis $n$, $d\in\N\setminus\{0\}$ are two positive integers.

Let $K$ be a field. We assume that $\Char K=0$ throughout the sections \ref{subsection Borel sets} and \ref{section binomial ideals}. 

An expression like $1\leq i\leq n$ has to be read as $i\in\{1,\dots,n\}$. An expression like $1\leq i\leq j\leq n$ has to be read as 
$$(i,j)\in\{(k,l)\in\N\times\N\mid 1\leq k\leq l\leq n\}.$$

We fix a polynomial ring $S:=K[X_1,\dots,X_n]$. (In the last section $S$ will denote the polynomial ring $K[X_0,\dots,X_n]$ in one variable more.) Any polynomial ring is endowed with the standard $\Z$-grading. 

By a \emph{ring} we always mean a commutative ring with identity element. In particular, all $K$-algebras are supposed to be commutative.
\end{convention}

\begin{definition}
If $R$ is a graded ring, $M$ a graded $R$-module and $N$ a subset of $M$, we denote by $\langle N\rangle_R$ the submodule of $M$ generated by $N$. If $i\in\Z$, we denote the $i$th graded component of $M$ by $M_i$. For $k\in\Z$ let $M_{\geq k}:=\bigoplus_{i\geq k}M_i$ denote the \emph{$k$-truncation} of $M$. We say that $M$ is \emph{generated in degree $d$} if $\langle M_d\rangle_R=M_{\geq d}$.

Let $R=\bigoplus_{i\in\N}R_i$ be a homogeneous Noetherian ring with $R_0=K$ and let $R_+:=\bigoplus_{i>0}R_i$ denote the irrelevant ideal of $R$. Let $M$ be a graded $R$-module. For $i\in\N$ let $H^i_{R_+}(M)$ denote the $i$th \emph{local cohomology module} of $M$ with respect to $R_+$, endowed with its natural grading (cf.\ \cite[Chap.~12]{BS}). 

If $M$ is finitely generated, we denote by 
$$h_M:\Z\ra\N,\ i\mapsto\dim_K M_i$$ 
the \emph{Hilbert function} of $M$. There exists a polynomial $p\in\Q[t]$,  called the \emph{Hilbert polynomial} of $M$, such that 
$$p(i)=h_M(i)\text{ for all }i\gg 0.$$ 
A polynomial $p\in\Q[t]$ is called an \emph{admissible Hilbert polynomial} if there exists a homogeneous ideal $\ideala\subset S$ with Hilbert polynomial $p$.
\end{definition}

\begin{notation}
For any set $L$ and any integer $m\in\N$ let $L^{(n,m)}$ denote the set of all $n\times m$ matrices $[M_{ij}\mid 1\leq i\leq n,\:1\leq j\leq m]$ with entries $M_{ij}\in L$. 

Let $g=[g_{ij}\mid 1\leq i\leq n,\:1\leq j\leq n]\in K^{(n,n)}$ be a matrix. Then we denote by $g:S\ra S$ the homomorphism of $K$-algebras defined by $X_j\mapsto\sum_{i=1}^n g_{ij}X_i$ for $1\leq j\leq n$.

\sk
If $N$ is a set, $\#N$ denotes its cardinality.

For two functions $f,g\colon\Z\ra\N$ we write $f\geq g$ if $f(k)\geq g(k)$ for all $k\in \Z$.

For $1\leq i\leq n$, let $e_i\in\N^n$ denote the standard vector with $(e_i)_j=
1$ if $i=j$ and $(e_i)_j=0$ otherwise.

\sk
For $a=(a_1,\dots,a_n)$, $b=(b_1,\dots,b_n)\in\Z^n$ and for a subset $A\subset\N^n$ we introduce the following notations:

$|a|:=\sum_{i=1}^n a_i$,

$m(a):=\max{(\{1\leq i\leq n\mid a_i\neq 0\}\cup\{1\})}$,

$\mu(a):=\min{(\{1\leq i\leq n\mid a_i\neq 0\}\cup\{n\})}$,

$a+k b:=(a_1+k b_1,\dots,a_n+k b_n)$ for any $k\in\Z$,

$A+a:=\{c+a\mid c\in A\}$,

$a^*:=a-a_n e_n=(a_1,\dots,a_{n-1},0)$,

$A^*:=\{c^*\mid c\in A\}$,

$A^{(i)}:=\{c\in A\mid m(c)=i\}$ for $1\leq i\leq n$,

$A(i):=\{c\in A\mid c_n=i\}$ for $i\in\N$,

$a^+:=(\max{\{a_1,0\}},\dots,\max{\{a_n,0\}})\in\N^n$, 

$a^-:=a^+-a\in\N^n$,

$\N_d^n:=\{c\in\N^n\mid |c|=d\}$.
\end{notation}

\begin{definition}\label{def - deglex, degrevlex - def}
In $\N^n$ we define the \emph{homogeneous lexicographic order} $\deglexg$ and the \emph{reverse lexicographic order} $\degrevlexg$ as follows: Let $a,b\in\N^n$. 

$a\deglexg b\ \defiff\ |a|>|b|$ or ($|a|=|b|$ and $a_{\mu(a-b)}>b_{\mu(a-b)}$).

$a\degrevlexg b\ \defiff\ |a|>|b|$ or ($|a|=|b|$ and $a_{m(a-b)}<b_{m(a-b)}$).
\end{definition}

\begin{definition}\label{def - term order of N^n, admissible - def}
A \emph{term order} of $\N^n$ is a total order of $\N^n$ with the following properties:
\begin{enumerate}
\item $(0,\dots,0)\leq a$ for all $a\in\N^n$,
\item $a<b\ \Rightarrow\ a+c<b+c$ for all $a,b,c\in\N^n$.
\end{enumerate}

A term order of $\N^n$ is called \emph{admissible} if it has the further properties
\begin{enumerate}
\setcounter{enumi}{2}
\item $e_1>\dots>e_n$,
\item $|a|>|b|\ \Rightarrow\ a>b$ for all $a,b\in\N^n$.
\end{enumerate}
\end{definition}

\begin{remark}
The homogeneous lexicographic order and the reverse lexicographic order are admissible term orders of $\N^n$.
\end{remark}

\begin{notation}
For $a=(a_1,\dots,a_n)\in\Z^n$ write $X^a:=X_1^{a_1}\cdots X_n^{a_n}$.

For a subset $A\subset\N^n$ one puts $X^A:=\{X^a\in\Z[X_1,\dots,X_n]\mid a\in A\}$.

We consider the set of \emph{monomials} $\terms:=X^{\N^n}$ as a subset of $R[X_1,\dots,X_n]$ for any ring $R$, and we define $\log:\terms\ra\N^n$, $X^a\mapsto a$.

Let $R$ be a ring. For $f\in R[X_1,\dots,X_n]$ and $m\in\terms$ let $c_m^f\in R$ denote the coefficient of $f$ with respect to the monomial $m$. This means that we can write $f=\sum_{m\in\terms}c^f_m m$.

The set $\supp_R(f):=\{m\in\terms\mid c^f_m\neq 0\}$ is called the \emph{support} of $f\in R[X_1,\dots,X_n]$. For $f\in S$ we write $\supp(f):=\supp_K(f)$.
\end{notation}

\begin{definition}\label{def - term order of T(S), deglex, degrevlex - def}
An order of $\terms$ is called a \emph{term order} if it is induced by a term order of $\N^n$. A term order of $\terms$ is called \emph{admissible} if it is induced by an admissible term order of $\N^n$.

The term orders of $\terms$ which are induced by the homogeneous lexicographic order and by the reverse lexicographic order of $\N^n$ are called  \emph{homogeneous lexicographic order} of $\terms$ and \emph{reverse lexicographic order} of $\terms$ respectively.
\end{definition}

\begin{definition}
Let $B\subset\N^n_d$. A subset $A\subset B$ is called a \emph{lexicographic segment} of $B$ if for all $a,b\in B$ with $a\in A$ and $b\deglexg a$ we have $b\in A$.

The subset $B\subset\N^n_d$ is called \emph{lexicographic} if $B$ is a lexicographic segment of $\N^n_d$.

An ideal of $S$ is called a \emph{monomial ideal} if it is generated by elements of $\terms$.

A monomial ideal $\ideala\subset S$ is called a \emph{lex ideal} if $\log{(\ideala_i\cap\terms)}$ is lexicographic for all $i\in\N$.
\end{definition}

\begin{remark}
Since the homogeneous lexicographic order is total, there exists for each $0\leq m\leq\#\N^n_d$ a unique lexicographic set $B\subset\N^n_d$ with $m$ elements. It follows that lex ideals are uniquely determined by their Hilbert function.

If $\ideala\subset S$ is a lex ideal, then $\ideala^\sat$ is also a lex ideal.

If $\idealb\subset S$ is a saturated homogeneous ideal and $k\in\N$, then $(\idealb_{\geq k})^\sat=\idealb$.
\end{remark}

\begin{lemma}\label{existence of lex ideals}
a) Let $\ideala\subset S$ be a homogeneous ideal. Then there exists a unique lex ideal $\ideala^\lex\subset S$ with $h_\ideala=h_{\ideala^\lex}$.

b) Let $p\in\Q[t]$ be an admissible Hilbert polynomial. Then there exists a unique saturated lex ideal $\ideall_p\subset S$ with Hilbert polynomial $p$.
\end{lemma}

\begin{proof}
a) See for example \cite[4.]{Sp}.

b) Let $\ideala\subset S$\/ be a homogeneous ideal with Hilbert polynomial $p$. Set $\ideall_p:=(\ideala^\lex)^\sat$. This is a saturated lex ideal with Hilbert polynomial $p$. Let $\idealb\subset S$ be another saturated lex ideal with Hilbert polynomial $p$. Then there exists $k\in\N$ such that $h_\idealb(i)=p(i)=h_{\ideall_p}(i)$ for all $i\geq k$. It follows that $\idealb_{\geq k}=(\ideall_p)_{\geq k}$ and therefore $\idealb=(\idealb_{\geq k})^\sat=((\ideall_p)_{\geq k})^\sat=\ideall_p$.
\end{proof}

\begin{definition}\label{def - lt, in - def}
Let $\tau$ be a term order of $\terms$ and let $f\in S$. If $f\neq 0$, the \emph{leading term} of $f$ is defined by $\Lt_\tau(f):=\max_\tau\supp(f)$; furthermore we define $\Lt_\tau(0):=0$. If $M\subset S$ is a subset, let  $\Lt_\tau M:=\{\Lt_\tau(f)\mid f\in M\}$.

Let $\ideala\subset S$ be an ideal. Then the \emph{initial ideal} of $\ideala$ with respect to $\tau$ is defined to be  $\init_\tau\ideala:=\langle\Lt_\tau\ideala\rangle_S$.
\end{definition}

\begin{lemma}[{\cite[15.26]{E}}]\label{ideal and initial ideal have same Hilbert function}
Let $\ideala\subset S$ be a homogeneous ideal and let $\tau$ be a term order of\/ $\terms$. Then $h_\ideala=h_{\init_\tau\ideala}$ and $h_{S/\ideala}=h_{S/\init_\tau\ideala}$.\hfill$\square$
\end{lemma}
\sk 
\subsection{Weight orders}
One has a standard procedure to connect an ideal $\ideala\subset S$ with its initial ideal with respect to any term order $\tau$ in the Hilbert scheme by an affine line (cf.~Proposition~\ref{proposition about the affine line}). The idea is to find a flat family of $K[z]$-algebras whose fiber over $1$ is $S/\ideala$ and whose fiber over $0$ is $S/\init_\tau\ideala$.

\begin{notation}\label{notation Weight orders}
Let $R$ be a ring and $\omega\in\Z^n$. Define 
$$\alpha^R_\omega:R[X_1,\dots,X_n]\setminus\{0\}\ra\Z,\;f\mapsto\sup{\{\sum_{i=1}^n(\omega_i\,(\log m)_i)\mid m\in\supp_R(f)\}}.$$ 

Let $\phi:R\ra R'$ be a homomorphism of rings. For $a\in R'$ define 
$$\beta^{\phi,a}_\omega:R[X_1,\dots,X_n]\ra R'[X_1,\dots,X_n],\;f\mapsto\hspace{-0.3cm}\smash{\sum_{m\in\supp_R(f)}}\!\!a^{\alpha^R_\omega(f)-\alpha^R_\omega(m)}\phi(c^f_m) m.$$
If $R'=R$ and $\phi=\id_R$ we write $\beta^a_\omega:=\beta^{\phi,a}_\omega$.
\end{notation}

\begin{lemma}\label{lemma on weighted ideals}
Let $\ideala\subset S$ be a homogeneous ideal, $a\in K$ and $\omega\in\Z^n$. Let $K[z]$ be a polynomial ring in one variable and $\iota:K\ra K[z]$ the canonical inclusion. Let $\phi:K[z]\ra K$ denote the homomorphism of $K$-algebras defined by $\phi(z)=a$. Let $h:K[z][X_1,\dots,X_n]\otimes_{K[z]}K\ra S$ denote the canonical isomorphism and $j:\langle\beta^{\iota,z}_\omega(\ideala)\rangle_{K[z][X_1,\dots,X_n]}\hookrightarrow K[z][X_1,\dots,X_n]$ the inclusion map. Then there exists an isomorphism
$$\langle\beta^{\iota,z}_\omega(\ideala)\rangle_{K[z][X_1,\dots,X_n]}\otimes_{K[z]}K\cong\langle\beta^a_\omega(\ideala)\rangle_S$$
such that the diagram
$$\xymatrix{
K[z][X_1,\dots,X_n]\otimes_{K[z]}K\ar[r]^-{h}_-{\cong}&S\\
\langle\beta^{\iota,z}_\omega(\ideala)\rangle_{K[z][X_1,\dots,X_n]}\otimes_{K[z]}K\ar[r]^-{\cong}\ar@<-11mm>[u]^{j\otimes_{K[z]} K}&\langle\beta^a_\omega(\ideala)\rangle_S\ar@{ (->}[u]
}$$
commutes.
\end{lemma}

\begin{proof}
Let $\bar\phi:K[z][X_1,\dots,X_n]\ra S$ be the homomorphism of $S$-algebras induced by $\phi$. Define $h':=h\circ(j\otimes_{K[z]}K)$. Then $h'(g\otimes k)=\bar\phi(k\,g)$ for all $k\in K$ and all $g\in\langle\beta^{\iota,z}_\omega(\ideala)\rangle_{K[z][X_1,\dots,X_n]}$. 

For $f\in\ideala$ it holds 
\begin{align*}
h'(\beta^{\iota,z}_\omega(f)\otimes 1)&=\bar\phi(\beta^{\iota,z}_\omega(f)) =\bar\phi(\sum_{m\in\supp_{K[z]}(f)}z^{\alpha^K_\omega(f)-\alpha^K_\omega(m)}c^f_m m)\\ &=\sum_{m\in\supp_{K[z]}(f)}a^{\alpha^K_\omega(f)-\alpha^K_\omega(m)}c^f_m m =\beta^a_\omega(f).
\end{align*}
We deduce that $\Im h'=\langle\beta^a_\omega(\ideala)\rangle_S$. 

Let $g\in\Ker h'$. Write $g=\sum_{i=1}^r g_i\otimes k_i$ with $g_i\in\langle\beta^{\iota,z}_\omega(\ideala)\rangle_{K[z][X_1,\dots,X_n]}$ and $k_i\in K$ and set $g':=\sum_{i=1}^r k_i\,g_i$. Then 
\begin{align*}g&=\sum_{i=1}^r g_i\otimes k_i\ =\ g'\otimes 1\ =\!\!\sum_{m\in\supp_{K[z]}(g')}\hspace{-0.4cm}c^{g'}_m m\otimes 1\ =\!\!\sum_{m\in\supp_{K[z]}(g')} m\otimes \phi(c^{g'}_m)\\
&=\sum_{m\in\supp_{K[z]}(g')}\hspace{-0.4cm}\phi(c^{g'}_m)\,m\otimes 1\ =\ \bar\phi\bigl(\sum_{m\in\supp_{K[z]}(g')}\hspace{-0.4cm}c^{g'}_m m\bigr)\otimes 1\ =\ \bar\phi(g')\otimes 1\\
&=\sum_{i=1}^r \bar\phi(k_i\,g_i)\otimes 1=\sum_{i=1}^r h'(g_i\otimes k_i)\otimes 1=h'(g)\otimes 1=0.
\end{align*}
This shows that $\Ker h'=0$ and hence the statement is proved.
\end{proof}

\begin{remark}\label{remark about weight orders}
Let $R$ be a ring, $a\in R$, $\ideala\subset R[X_1,\dots,X_n]$ a homogeneous ideal and $\omega\in\Z^n$. Then $\langle\beta^a_\omega(\ideala)\rangle_{R[X_1,\dots,X_n]}$ is homogeneous. 

Moreover, if $R=K$ and if $a\in K^*$ is a unit, then $S/\ideala$ and $S/\langle\beta^a_\omega(\ideala)\rangle_S$ are isomorphic as graded $S$-modules (cf.~\cite[15.17]{E}).
\end{remark}

\begin{definition}
For an ideal $\ideala\subset S$ and $\omega\in\Z^n$ define
$\init_\omega\ideala:=\langle\beta^0_\omega(\ideala)\rangle_S.$
\end{definition}

\begin{proposition}[{\cite[15.16, Ex.\ 15.12]{E}}]\label{existence of weight orders}
Let $\ideala\subset S$ be a homogeneous ideal and $\tau$ a termorder of\/ $\terms$. Then there exists $\omega\in\Z^n$ such that $\init_\tau\ideala=\init_\omega\ideala$.\hfill$\square$
\end{proposition}

\begin{proposition}[{\cite[15.17]{E}}]\label{a flat homomorphism}
Let $\ideala\subset S$ be a homogeneous ideal and $\omega\in\Z^n$. Let $K[z]$ be a polynomial ring in one variable and $\iota:K\ra K[z]$ the canonical inclusion. Then the canonical homomorphism of rings
$$K[z]\ra K[z][X_1,\dots,X_n]/\langle\beta^{\iota,z}_\omega(\ideala)\rangle_{K[z][X_1,\dots,X_n]}$$
is flat.\hfill$\square$
\end{proposition}
\sk
\subsection{Borel sets}\label{subsection Borel sets}
In this section we assume that $\Char K=0$. Then the Borel-fixed ideals $\ideala\subset S$ are monomial ideals which are characterized by the following property: If a monomial $m\in\ideala$ is divisible by an indeterminate $X_j$, then $\frac{X_i}{X_j}\,m\in\ideala$ for all $1\leq i\leq j$. In each homogeneous component they correspond to so called Borel sets. Borel sets are the Borel order analogue to lexicographic sets. 
\pagebreak[1]

\sk
There are several equivalent ways to define the Borel order. The most plausible is the following one: For all monomials $m\in\terms$ and for all $1\leq k<n$ set $X_k\,m\borelg X_{k+1}\,m$ and take the associative hull:
\begin{definition}\label{def - Borel order, set - def}
Define the \emph{Borel order} $\borelgeq$ of $\N^n_d$ as follows. Let $a,b\in\N_d^n$. 

$a\borelgeq b\ \defiff\ \forall\ 1\leq k<n\ \exists\alpha_k\in\N:a-b=\sum_{j=1}^{n-1}\alpha_j(e_j-e_{j+1})$.

A set $A\subset\N^n_d$ is called a \emph{Borel set} if for all $a,b\in\N^n_d$ with $a\in A$ and $b\borelgeq a$ we have $b\in A$.
\end{definition}

For technical reasons it will be more convenient to have a description of the Borel order by upper triangular integer matrices (s.\ Lemma \ref{properties of Borel order}). Remark that the Borel order is not total, so it is not a term order.

As a consequence of the fact that generic initial ideals are Borel-fixed irrespective of the admissible term order (Proposition \ref{existence of gin}) we have the following characterization of the Borel order:

\begin{lemma}[{\cite[2.2]{C}}]\label{Borel greater implies greater for any term order}
Let $a,b\in\N^n_d$. Then it holds $a\borelgeq b$ if and only if $a\geq_\tau b$ for all admissible term orders $\tau$ of\/ $\N^n$.\hfill$\square$
\end{lemma}

\begin{definition}\label{def - root - def}
An element of a Borel set $B\subset\N^n_d$ which is minimal with respect to the Borel order is called a \emph{root} of $B$. 
\end{definition}

\begin{definition}\label{def - monomial, Borel ideal - def}
A monomial ideal $\ideala\subset S$ is called a \emph{Borel ideal} if $\log(\ideala_i\cap\terms)$ is a Borel set for all $i\in\N$.
\end{definition}

\begin{remark}\label{Groebner bases of Borel ideals}
If $\ideala\subset S$ is a Borel ideal, then  $\ideala^\sat=(\ideala:X_n^\infty)$ (cf.~\cite[15.24]{E}). Hence, if $B\subset\N^n_d$ is a Borel set, then the set $X^{B^*}$ generates the ideal $\langle X^B\rangle_S^\sat$. Since this ideal is monomial, $X^{B^*}$ is a Gr\"obner basis of $\langle X^B\rangle_S^\sat$ with respect to any admissible term order of\/ $\terms$.

Let $B\subset\N^n_d$ be a Borel set. Since $\langle X^B\rangle_S$ is generated in degree $d$, it is clear that $(\langle X^B\rangle_S^\sat)_{\geq d}\supset\langle X^B\rangle_S$. On the other hand, as $\langle X^B\rangle_S$ is a Borel ideal, it follows that $(\langle X^{B^*}\rangle_S)_{\geq d}\subset\langle X^B\rangle_S$. We conclude that $(\langle X^B\rangle_S^\sat)_{\geq d}=\langle X^B\rangle_S$.  
\end{remark}

\begin{notation}
Set $U(n):=\{M\in \N^{(n,n)}\mid M_{ij}=0\ \ \forall\ 1\leq j<i\leq n\}$. 
\end{notation}

\begin{lemma}\label{properties of Borel order}
Let $a,b\in\N_d^n$. Then the following are equivalent:
\begin{align*}
&\textup{(i)}&&a\borelgeq b,\\
&\textup{(ii)}&&\exists M\in U(n):\ \sum_{i=1}^n M_{ij}=b_j\ \ \forall\ 1\leq j\leq n,\ \sum_{j=1}^n M_{ij}=a_i\ \ \forall\ 1\leq i\leq n,\\
&\textup{(iii)}&&\smash[t]{\sum_{k=1}^i a_k\geq\sum_{k=1}^i b_k}\ \ \forall\ 1\leq i\leq n.
\end{align*}
\end{lemma}

\begin{proof}
(i)$\Ra$(ii): Let $a\borelgeq b$. Choose $\alpha_k\in\N$ for all $1\leq k<n$  such that $a-b=\sum_{k=1}^{n-1}\alpha_k(e_k-e_{k+1})$. Set $\alpha_0:=0$ and 
$m:=-\sum_{i=1}^n \min{\{0,b_i-\alpha_{i-1}\}}.$ We will construct a sequence of matrices $M(0)$, \dots, $M(m)\in\Z^{(n,n)}$ such that for all $0\leq k\leq m$ the following properties hold:
\begin{align*}
&\text{(1)}\hspace{12pt}M(k)_{ij}=0\text{ for all }1\leq j<i\leq n,\\
&\text{(2)}\hspace{12pt}M(k)_{ij}\geq 0\text{ for all }1\leq i<j\leq n,\\
&\text{(3)}\hspace{12pt}\sum_{i=1}^n M(k)_{ij}=b_j\text{ for all }1\leq j\leq n,\\
&\text{(4)}\hspace{12pt}\sum_{j=1}^n M(k)_{ij}=a_i\text{ for all }1\leq i\leq n,\\
&\text{(5)}\hspace{12pt}\sum_{i=1}^n \min{\{0,M(k)_{ii}\}}=k-m.
\end{align*}
Property (5) implies that $M(m)_{ii}\geq 0$ for all $1\leq i \leq n$, hence $M(m)$ is the requested matrix in $U(n)$.

Define $M(0)\in\Z^{(n,n)}$ by $M(0)_{ij}:=
\begin{cases}
\alpha_{i},&\text{if $i=j-1$;}\\
b_i-\alpha_{i-1},&\text{if $i=j$;}\\
0,&\text{otherwise.}
\end{cases}$\\
It is clear that $M(0)$ has the properties (1), (2), (3) and (5). Setting $\alpha_n:=0$, we have for all $1\leq i\leq n$
$$\sum_{j=1}^n M(0)_{ij}=b_i-\alpha_{i-1}+\alpha_i=b_i+\bigl(\sum_{k=1}^{n-1}\alpha_k(e_k-e_{k+1})\bigl)_i=a_i,$$
whence $M(0)$ has also property (4).

If $m>0$, we construct $M(1)$, \dots, $M(m)$ recursively. Let $0\leq k<m$ and assume that $M(k)\in\Z^{(n,n)}$ with the required properties is constructed already. Since $k<m$, by property (5) there exists 
$$l:=\min{\{1\leq j\leq n\mid M(k)_{jj}<0\}}.$$ 
Properties (1), (3) and (4) imply that there exist $p,q\in\N$ with $1\leq p<l<q\leq n$ and $M(k)_{pl}$, $M(k)_{lq}>0$. Now define $M(k+1)\in\Z^{(n,n)}$ by $$M(k+1)_{ij}:=
\begin{cases}
M(k)_{ij}+1,&\text{if $(i,j)\in\{(l,l),(p,q)\}$;}\\
M(k)_{ij}-1,&\text{if $(i,j)\in\{(p,l),(l,q)\}$;}\\
M(k)_{ij},&\text{otherwise.}
\end{cases}$$
It is clear that $M(k+1)$ has the required properties.

\sk
(ii)$\Ra$(iii): Let $M\in U(n)$ be such that $\sum_{i=1}^n M_{ij}=b_j$ for all $1\leq j\leq n$ and $\sum_{j=1}^n M_{ij}=a_i$ for all $1\leq i\leq n$. Then for all $1\leq i\leq n$ it holds
$$\sum_{k=1}^i a_k=\sum_{k=1}^i\sum_{j=1}^n M_{kj}=\sum_{j=1}^n\sum_{k=1}^i M_{kj}\geq\sum_{j=1}^i\sum_{k=1}^i M_{kj}=\sum_{j=1}^i\sum_{k=1}^n M_{kj}=\sum_{j=1}^i b_j.$$

\sk
(iii)$\Ra$(i): Let $a,b\in\N_d^n$ be such that $\sum_{i=1}^k a_i\geq\sum_{i=1}^k b_i$ for all $1\leq k\leq n$. For $1\leq k<n$ set $\alpha_k:=\sum_{i=1}^k a_i-\sum_{i=1}^k b_i$. Then 
\begin{align*}
a-b&=\sum_{k=1}^n(a_k-b_k)e_k=(a_1-b_1)e_1+\bigl(d-\sum_{i=1}^{n-1}a_i-\bigl(d-\sum_{i=1}^{n-1}b_i\bigr)\bigr)e_n\\*
&\quad +\sum_{k=2}^{n-1}\bigl(\sum_{i=1}^k(a_i-b_i)-\sum_{i=1}^{k-1}(a_i-b_i)\bigr)e_k \\
&=\sum_{k=1}^{n-1}\sum_{i=1}^k(a_i-b_i)e_k-\sum_{k=2}^n\sum_{i=1}^{k-1}(a_i-b_i)e_k\\
&=\sum_{k=1}^{n-1}\sum_{i=1}^k(a_i-b_i)(e_k-e_{k+1})\\
&=\sum_{k=1}^{n-1}\alpha_k(e_k-e_{k+1}).
\end{align*}
\end{proof}

It follows immediately from Lemma \ref{Borel greater implies greater for any term order} that lexicographic sets are Borel sets. The idea of D. Mall's proof of the connectedness of certain Hilbert scheme strata is to find a sequence of Borel sets which are ``more and more lexicographic''. For this way of proceeding the following definition is crucial: 

\begin{definition}
Let $B\subset\N^n_d$ be a Borel set. Then $B$ is called \emph{growth-height-lexicographic} if the following conditions hold:
\begin{enumerate}
\item $B^{(i)}$ is a lexicographic segment of $(\N^n_d)^{(i)}$ for all $1\leq i<n$,
\item $B(i)$ is a lexicographic segment of $\N^{n}_d(i)$ for all $1\leq i\leq d$.
\end{enumerate}

Define the \emph{growth vector} of $B$ by $\gv(B):=(\#B^{(1)},\dots,\#B^{(n)})$ and the \emph{height vector} of $B$ by $\hv(B):=(\#B(1),\dots,\#B(d))$.
\end{definition}

\begin{remark}
Let $A$ be a lexicographic segment of $B\subset\N^n_d$. Then $A$ is uniquely determined by its cardinality and by $B$. Thus a growth-height-lexicographic Borel set $B\subset\N^n_d$ is uniquely determined by its growth and height vectors.

It holds $\left|\hv(B)\right|=\gv(B)_n$.
\end{remark}

\begin{proposition}[{\cite[2.17]{M1997}}]\label{existence of growth-height-lexicographic normal form}
Let $B\subset\N^n_d$ be a Borel set. Then there exists a unique growth-height-lex\-i\-co\-graphic Borel set $L_{gh}(B)\subset\N^n_d$ with $\gv(L_{gh}(B))=\gv(B)$ and $\hv(L_{gh}(B))=\hv(B)$.

The set $L_{gh}(B)$ is called \emph{growth-height-lexicographic normal form} of $B$.\\
\hspace*{\fill}$\square$
\end{proposition}

\begin{proposition}[{\cite[2.9]{M2000}}]\label{computation of the Hilbert function}
Let $A$, $B\subset\N^n_d$ be Borel sets. Then we have

a) $h_{\langle X^A\rangle_S}=h_{\langle X^B\rangle_S}\iff\gv(A)=\gv(B)$.

b) $h_{\langle X^B\rangle_S^\sat}(i)=
\begin{cases}
{\displaystyle\sum_{j=0}^{n-1}\gv(B)_{n-j}\binom{i-d+j}{j}},&\text{if $i\geq d$;}\\
{\displaystyle\sum_{j=0}^i\hv(B)_{d-j}},&\text{if $0\leq i<d$.}
\end{cases}$\\
\hspace*{\fill}$\square$
\end{proposition}

\begin{corollary}\label{Hilb poly determines growth-vector}
Let $A$, $B\subset\N^{n+1}_d$ be two Borel sets such that the ideals $\langle X^A\rangle_S$ and $\langle X^B\rangle_S$ have the same Hilbert polynomial. Then it holds $\gv(A)=\gv(B)$.
\end{corollary}

\begin{proof}
Let $p\in\Q[t]$ be the Hilbert polynomial of $\langle X^A\rangle_S$, resp.\ $\langle X^B\rangle_S$. It holds $\langle X^A\rangle_S=(\langle X^A\rangle_S)^\sat_{\geq d}$ and $\langle X^B\rangle_S=(\langle X^B\rangle_S)^\sat_{\geq d}$ (cf.~Remark~\ref{Groebner bases of Borel ideals}). By Proposition \ref{computation of the Hilbert function}~b) it follows that $h_{\langle X^A\rangle_S}(i)=p(i)=h_{\langle X^B\rangle_S}(i)$ for all $i\geq d$, whence $h_{\langle X^A\rangle_S}=h_{\langle X^B\rangle_S}$. The corollary now follows from part a) of Proposition \ref{computation of the Hilbert function}.
\end{proof}
\sk
\subsection{Generic initial ideals and reverse lexicographic order}\label{subsection Gin}
Generic initial ideals have a lot of nice properties. We already mentioned that they are Borel-fixed. In this subsection we show that their formation commutes with truncation, and if $\Char K=0$ and $S$ is endowed with the reverse lexicographic order, then their formation commutes with saturation.

\begin{definition}
The \emph{unipotent subgroup} $\unipotent\subset\Gl(n,K)$ is the group of all upper triangular matrices with ones on the diagonal.
\end{definition}

\begin{proposition}\label{existence of gin}
Let $\ideala\subset S$ be a homogeneous ideal and $\tau$ an admissible term order of\/ $\terms$. 

a) \cite[15.18]{E} There is a non-empty Zariski open set $U\subset\Gl(n,K)$ and a unique ideal\/ $\Gin_\tau\ideala\subset S$ such that $\Gin_\tau\ideala=\init_\tau g(\ideala)$ for all $g\in U$. 

Furthermore, the open set $U$ meets the unipotent group $\unipotent$.

The ideal\/ $\Gin_\tau\ideala$ is called the \emph{generic initial ideal} of $\ideala$ with respect to $\tau$.

b) \cite[15.20]{E} The generic initial ideal $\Gin_\tau\ideala$ is \emph{Borel-fixed}, i.\,e.\ for all upper triangular matrices $g\in\Gl(n,K)$ it holds $g(\Gin_\tau\ideala)=\Gin_\tau\ideala$. 

c) \cite[15.23]{E} If $\Char K=0$, then an ideal $\idealb\subset S$ is Borel-fixed if and only if it is a Borel ideal.\hfill$\square$
\end{proposition}

\begin{corollary}\label{Gin=in for unipotent fixed ideals}
If a homogeneous ideal of $S$ remains fixed under the action of $\unipotent$, then its generic initial ideal and its initial ideal with respect to any admissible term order coincide.\hfill$\square$
\end{corollary}

\begin{lemma}\label{gin commutes with truncation}
Let $\ideala\subset S$ be a homogeneous ideal and $\tau$ an admissible term order of\/ $\terms$. Then $\init_\tau(\ideala_{\geq d})=(\init_\tau\ideala)_{\geq d}$ and $\Gin_\tau(\ideala_{\geq d})=(\Gin_\tau\ideala)_{\geq d}.$
\end{lemma}

\begin{proof}
The first claim is obvious. Choose $g\in\Gl(n,K)$ such that $\Gin_\tau(\ideala_{\geq d})=\init_\tau g(\ideala_{\geq d})$ and $\Gin_\tau\ideala=\init_\tau g(\ideala)$. It follows that
$\Gin_\tau(\ideala_{\geq d})=\init_\tau g(\ideala_{\geq d})=\init_\tau (g(\ideala)_{\geq d})=(\init_\tau g(\ideala))_{\geq d}=(\Gin_\tau\ideala)_{\geq d}.$
\end{proof}

\begin{proposition}\label{Gin and sat commute}
Assume that $\Char K=0$. Let $\ideala\subset S$ be a homogeneous ideal. Then
$$\Gin_\degrevlex(\ideala^\sat)=(\Gin_\degrevlex\ideala)^\sat.$$
\end{proposition}

\begin{proof}
Let $P:=\bigcup_{\idealp\in\Ass(S/\ideala)\setminus\{S_+\}}\idealp$ be the set of all elements of $S$ contained in some associated prime of $S/\ideala$ excepting possibly the irrelevant ideal. We first show that $(\ideala: u^\infty)=\ideala^\sat$ for all $u\in S_1\setminus P$. Let $u\in S_1\setminus P$. Let $r\in\N$ be such that $\ideala^\sat=(\ideala:_S S_+^r)$. Since $u^r\notin P$, it holds $(0:_{S/\ideala}u^r)=H^0_{S_+}(0:_{S/\ideala}u^r)$ (cf.\ \cite[18.3.8\,(iii)]{BS}). Hence, for any $f\in(\ideala:_S u^r)$ there exists $m\in\N$ such that $S_+^m f\subset\ideala$, which means that $f\in\ideala^\sat$. Since $u^r\in S_+^r$, it is clear that $\ideala^\sat=(\ideala:_S\nolinebreak S_+^r)\subset(\ideala:_S u^r)$. It follows $(\ideala:_S u^r)=\ideala^\sat$.

We next prove that there is a non-empty Zariski open set $U\subset\Gl(n,K)$ such that $g(\ideala^\sat)=(g(\ideala):X_n^\infty)$ for all $g\in U$. Since $\#K=\infty$, the open subset $S_1\setminus P\subset S_1$ is not empty, whence $U:=\{g\in\Gl(n,K)\mid g^{-1}(X_n)\notin P\}$ 
is a non-empty open subset of $\Gl(n,K)$. Let $g\in U$. Then it holds $g(\ideala^\sat)=g(\ideala:\nolinebreak g^{-1}(X_n)^\infty)=(g(\ideala):X_n^\infty)$.

Now, by Proposition \ref{existence of gin}, we may choose $g\in U$ such that $\Gin_\degrevlex\ideala=\init_\degrevlex g(\ideala)$ and $\Gin_\degrevlex(\ideala^\sat)=\init_\degrevlex g(\ideala^\sat)$. Since $\Gin_\degrevlex(\ideala)$ is a Borel ideal it holds $(\Gin_\degrevlex\ideala)^\sat=(\Gin_\degrevlex \ideala:X_n^\infty)$ (cf.\ Remark \ref{Groebner bases of Borel ideals}); and it follows from \cite[15.12]{E} that $\init_\degrevlex(g(\ideala):\nolinebreak X_n^\infty)=(\init_\degrevlex g(\ideala):X_n^\infty)$.
Altogether we obtain
$\Gin_\degrevlex(\ideala^\sat)=\init_\degrevlex g(\ideala^\sat)=\init_\degrevlex(g(\ideala):\nolinebreak X_n^\infty)=(\init_\degrevlex g(\ideala):X_n^\infty)
=(\Gin_\degrevlex \ideala:X_n^\infty)=(\Gin_\degrevlex\ideala)^\sat.$
\end{proof}

The following Lemma states a property of the reverse lexicographic order which is needed to prove Proposition \ref{Groebner bases of (saturated) binomial ideals}.
\begin{lemma}[{\cite[12.1]{Stu}}]\label{a result of Bayer Mumford}
Let $f_1$, \dots, $f_r\in S$ be homogeneous polynomials such that $X_n\nmid f_i$ for all $1\leq i\leq r$, let $a_1$, \dots, $a_r\in\N$ and let $\ideala:=\langle\{X_n^{a_i}f_i\mid 1\leq\nolinebreak i\leq\nolinebreak r\}\rangle_S$. If\/ $\{X_n^{a_i}f_i\mid 1\leq i\leq r\}$ is a Gr\"obner basis of $\ideala$ with respect to the reverse lexicographic order, then $\{f_i\mid 1\leq i\leq r\}$ is a Gr\"obner basis of\/ $(\ideala:X_n^\infty)$ with respect to the reverse lexicographic order. \hfill$\square$
\end{lemma}

\clearpage
\section{Binomial ideals}\label{section binomial ideals}
Throughout this section $K$ is a field of characteristic zero.

D. Mall showed in \cite{M2000} that for any Borel set $B\subset\N^n_d$ there exists a sequence of Borel sets $(B_i)_{i=0}^r$, beginning with $B_0=B$ and ending with its growth-height-lexicographic normal form $L_{gh}(B)=B_r$, such that all $B_i$ have the same growth and height vectors and such that for each $1\leq i\leq r$ there is a $\rho_i\in\Z^n$ with $(B_{i-1}\!\setminus\! B_i)+\rho_i=B_i\!\setminus\!B_{i-1}$ (cf.~Proposition \ref{rebuilding systems are Mall}). The triples $(B_{i-1}\cap B_i$, $B_{i-1}\setminus B_i$, $\rho_i)$ are called binomial systems. There is an analogous result, saying that the growth-height-lexicographic normal form $L_{gh}(B)$ and the lexicographic set $L$ with $\#L=\#B$ are connected by a sequence of binomial systems (cf.~Proposition \ref{existence of Mall bin syst to a lexicographic segment}). This sequences yield a sequence of Groebner deformations, connecting any point of the Hilbert scheme $\Hilb^p_K$ with the unique lexicographic point.

We want to use this sequences of Groebner deformations in order to prove that the subsets of the Hilbert scheme defined by bounding cohomology are connected.

Let $(B_i)_{i=0}^r$ be a sequence of Borel sets provided by Mall's algorithm, let $1\leq i\leq r$, and let $\idealb_{i-1}$, $\idealb_i\subset S$ be two ideals generated by $B_{i-1}$ and $B_i$ respectively. Let $\idealc_i\subset S$ be a binomial ideal given by the binomial system $(B_{i-1}\cap B_i$, $B_{i-1}\setminus B_i$, $\rho_i)$ (cf.~\ref{not-binomial ideal-not}). Then $\init_\degrevlex\idealc_i=\idealb_{i-1}$ and $\init_\deglex\idealc_i=\idealb_i$. Our aim is to compare the local cohomology modules of $S/\idealb_{i-1}$ and $S/\idealb_i$. This is possible by a Theorem of Herzog and Sbarra (s.~Proposition~\ref{main Theorem of Herzog-Sbarra}). In order to apply this Theorem in section \ref{subsection sequences}, we have to show two facts: In section~\ref{subsection generic initial ideal} we prove that $\init_\degrevlex\idealc_i=\Gin_\degrevlex\idealc_i$ (s.~Theorem~\ref{Gin = in for good binomial ideals}) and in section~\ref{subsection seq. CM} we prove that $S/\idealc_i$ is sequentially Cohen-Macaulay (s.~Theorem~\ref{binomial ideals and sCM}). Both results are a consequence of certain properties of the binomial system $(B_{i-1}\cap B_i$, $B_{i-1}\setminus B_i$, $\rho_i)$. The first section~\ref{subsection bin syst} is devoted to prove that all occurring binomial systems have several common properties providing the subsequent results. Since Mall did not state these properties explicitely, we have to go into the technical details of \cite{M2000}:

\subsection{Binomial systems}\label{subsection bin syst}
\begin{definition}\label{def - binomial system - def}
A triple $(A,C,\rho)$ consisting of two subsets $A,C\subset \N_d^n$ and of an $n$-tuple $\rho\in\Z^n$ is called a \emph{binomial system} (of degree $d$ in $n$ indeterminates) if the following conditions hold:
\begin{enumerate}
\item  $C+\rho\subset\N^n_d$,
\item $A\cap C=A\cap(C+\rho)=C\cap(C+\rho)=\emptyset$,
\item $A\cup C$ and $A\cup(C+\rho)$ are Borel sets.
\end{enumerate}
\end{definition}

\begin{remark}
If $(A,C,\rho)$ is a binomial system, we always assume that it is of degree $d$ in $n$ indeterminates unless otherwise stated.

If $(A,C,\rho)$ is a binomial system, then $A$ is a Borel set.

If $(A,C,\rho)$ is a binomial system, then for any term order of $\N^n$ we have: If $c<c+\rho$ for some $c\in C$, then $c<c+\rho$ for all $c\in C$.

If $(A,C,\rho)$ is a binomial system, then $(A,C+\rho,-\rho)$ is also a binomial system.
\end{remark}

\begin{definition}\label{def - admissible, good, Mall - def}
A binomial system $(A,C,\rho)$ is \emph{admissible} if $C$ is empty or if $\rho_{m(\rho)}>0$.

A binomial system $(A,C,\rho)$ is \emph{good} if it is admissible and if $b_i=c_i$ for all $1\leq i<m(\rho)$ and for all $b,c\in C$.

A binomial system $(A,C,\rho)$ is \emph{Mall} if it is good and if $\rho_{\mu(\rho)}>0$ and $m(c)=m(c+\rho)$ for all $c\in C$.
\end{definition}

\begin{remark}\label{remark on admissible binomial systems}
If $(A,C,\rho)$ is a binomial system which is not admissible, then $(A,C+\rho,-\rho)$ is an admissible binomial system.

If a binomial system $(A,C,\rho)$ is Mall, then $c\degrevlexg c+\rho$ and $c\deglexl c+\rho$ for all $c\in C$.
\end{remark}

\begin{notation}\label{not-binomial ideal-not}
If $C\subset\N^n$ and $\rho\in\Z^n$ are such that $C+\rho\subset\N^n$, set 
$$\Bin(C,\rho):=\{X^c-X^{c+\rho}\in S\mid c\in C\}.$$
If $(A,C,\rho)$ is a binomial system, set 
$$F(A,C,\rho):=\langle X^A\cup\Bin(C,\rho)\rangle_S.$$
\end{notation}

In the following Proposition we give some important properties of binomial ideals:

\begin{proposition}\label{Groebner bases of (saturated) binomial ideals}
Let $(A,C,\rho)$ be a binomial system.

a) If\/ $(A,C,\rho)$ is admissible, then $X^A\cup\Bin(C,\rho)$ is a Gr\"obner basis of $F(A,C,\rho)$ with respect to the reverse lexicographic order.

b) $F(A,C,\rho)^\sat=(F(A,C,\rho):X_n^\infty)$.

c) If\/ $(A,C,\rho)$ is admissible, then $X^{A^*}\cup\Bin(C^*,\rho)$ is a Gr\"obner basis of $F(A,C,\rho)^\sat$ with respect to the reverse lexicographic order.

d) If\/ $(A,C,\rho)$ is Mall, then $\init_\deglex F(A,C,\rho)=\langle X^{A\cup(C+\rho)}\rangle_S$.

e) If\/ $(A,C,\rho)$ is Mall and $\rho_n=0$, then  $$\init_\deglex(F(A,C,\rho)^\sat)=\langle X^{A^*\cup(C+\rho)^*}\rangle_S.$$

f) $(F(A,C,\rho)^\sat)_{\geq d}=F(A,C,\rho)$.
\end{proposition}

\begin{proof}
The first two properties are shown in \cite[3.7(1)]{M2000} and \cite[3.8]{M2000}. Assertion c) follows from a), b) and Lemma \ref{a result of Bayer Mumford}. Statement d) is a consequence of \cite[3.7(2)]{M2000} and of Remark \ref{remark on admissible binomial systems}, statement e) is proved in \cite[3.10]{M2000}.

\sk
f) We may assume that $(A,C,\rho)$ is admissible (cf.\ Remark \ref{remark on admissible binomial systems}). Then by statement c) it holds $F(A,C,\rho)^\sat=\langle X^{A^*}\cup\Bin(C^*,\rho)\rangle_S$. If $a\in A$ and $u\in\N^n_{a_n}$, then $u+a^*=a+u-a_n e_n\in A$ since $A$ is a Borel set, whence $X^u X^{a^*}\in F(A,C,\rho)_d$. If $c\in C$ and $v\in\N^n_{c_n}$, then $v+c^*\in A\cup C$ and $v+c^*+\rho\in A\cup(C+\rho)$, since $A\cup C$ and $A\cup(C+\rho)$ are Borel sets. Furthermore, it holds $v+c^*\in C$ if and only if $v+c^*+\rho\in C+\rho$, whence $X^v(X^{c^*}-X^{c^*+\rho})\in F(A,C,\rho)_d$. 

Now let $f\in(F(A,C,\rho)^\sat)_d$. Then we may write 
$$f=\sum_{a\in A} f_a X^{a^*}+\sum_{c\in C}g_c(X^{c^*}-X^{c^*+\rho})$$
with some homogeneous polynomials $f_a$, $g_c\in S$ with $\deg f_a=a_n$ and $\deg g_c=c_n$ for all $a\in A$, $c\in C$. By what has been shown above it is clear that all summands of $f$ are in $F(A,C,\rho)_d$, and hence $f\in F(A,C,\rho)_d$. Since the other inclusion is obvious, we have proven that $(F(A,C,\rho)^\sat)_d=F(A,C,\rho)_d$. Hence, our statement follows.
\end{proof}

\begin{proposition}\label{rebuilding systems are Mall}
Let $B\subset\N^n_d$ be a Borel set and let $L_{gh}(B)$ be its growth-height-lex\-i\-co\-graphic normal form. If $B\neq L_{gh}(B)$ then there exists a finite sequence of Mall binomial systems $((A_i,C_i,\rho_i))_{i=1}^r$ with the following properties:
\begin{enumerate}
\item $B=A_1\cup C_1$,
\item $A_i\cup (C_i+\rho_i)=A_{i+1}\cup C_{i+1}$ for all $1\leq i<r$,
\item $L_{gh}(B)=A_r\cup(C_r+\rho_r)$,
\item $m(\rho_i)<n$ for all $1\leq i\leq r$.
\end{enumerate}
\end{proposition}

\begin{proof}
Mall \cite{M2000} gives an algorithm (\cite[4.41]{M2000}) to find a sequence of binomial systems with the properties (i) -- (iii). Furthermore, these binomial systems satisfy the equations $\hv(A_i\cup C_i)=\hv(A_i\cup (C_i+\rho_i))$ for all $1\leq i\leq r$, and this implies the property (iv).

It remains to show that these binomial systems are Mall. 

It is enough to show that $(A,C,\rho):=(A_1,C_1,\rho_1)$ is Mall, then one can use an inductive argument to prove the general case. At this point it is necessary to enter deeper in the details of Mall's paper:

If $c$ is a root of $B$ and $1\leq i\leq n$, let $T_c(i)$ denote the lexicographically maximal element in the set 
$$\{a\in\N_d^n\setminus B\mid a\deglexg c,\:m(a)=m(c),\:a_j=c_j\ \forall\ i\leq j\leq n\}\cup\{c\}.$$
Since $B$ is not growth-height-lexicographic, there exists $2<m<n$ with the property that there is a root $c\in B$ such that $T_c(m+1)\neq c$ (cf.~Proposition \cite[4.31]{M2000}). Let now $2<m<n$ be the minimal integer with this property. (Then $B$ is said to be \emph{$m$-trivial, but not $(m+1)$-trivial}, cf.\ Definition \cite[4.27]{M2000} and Lemma \cite[4.28]{M2000}). 

In step (7) the algorithm of Mall chooses $\rho$ to be the lexicographically maximal element in the set $\{T_c(m+1)-c\in\Z^n\mid c \text{ is a root of }B\}$. Since $B$ is not $(m+1)$-trivial, we may select, in step (8) of the algorithm, a root $c_\rho$ of $B$ such that $\rho=T_{c_\rho}(m+1)-c_\rho\deglexg (0,\dots,0)$. Therefore, $\rho_{\mu(\rho)}>0$. 

On the other hand, since $B$ is $m$-trivial, but not $(m+1)$-trivial, it is shown by Mall that $m(\rho)=m$ (Lemma \cite[4.34]{M2000}) and that $(A,C,\rho)$ is admissible (Lemma \cite[4.35]{M2000}). 

In steps (9) and (10) of the algorithm, the set $C$ is constructed by adding some elements of $B$ to the set $\{c_\rho\}$. This is done in a finite loop by building so-called \emph{completion-adjunction pairs of width at most $n-m$} (cf.\ Lemma \cite[4.21]{M2000} and Lemma \cite[4.40]{M2000}). At the end of this procedure, in particular $(C,C)$ is a completion-adjunction pair of width at most $n-m$. So, by definition, we get $(c_\rho)_i=c_i$ for all $1\leq i<m=m(\rho)$ and for all $c\in C$ (Definition \cite[4.36(3)]{M2000}). Therefore, $(A,C,\rho)$ is good.

By the definition of $T_{c_\rho}(m+1)$ it holds $m(c_\rho)=m(T_{c_\rho}(m+1))$. It follows from $T_{c_\rho}(m+\nolinebreak 1)=c_\rho+\rho$ and $\rho_{m(\rho)}>0$ that $m(c_\rho)=m(c_\rho+\rho)\geq m(\rho)$. Let $c\in C$. Since $(A,C,\rho)$ is good, we have $(c_\rho)_i=c_i$ for all $1\leq i<m(\rho)$. Using the relations $\sum_{i=m(\rho)}^n c_i=\sum_{i=m(\rho)}^n (c_\rho)_i>0$, we deduce that $m(c)\geq m(\rho)$. Now $m(c+\rho)=m(c)$ follows from $\rho_{m(\rho)}>0$, and we conclude that $(A,C,\rho)$ is Mall.
\end{proof}

\begin{proposition}\label{existence of Mall bin syst to a lexicographic segment}
Let $B\subset\N^n_d$ be a growth-height-lexico\-graphic Borel set and let $L\subset\N^n_d$ be a lexicographic set such that $\gv(B)=\gv(L)$. If $L\neq B$, then there exists a finite sequence of Mall binomial systems $((A_i,C_i,\rho_i))_{i=1}^r$ with the following properties:
\begin{enumerate}
\item $B=A_1\cup C_1$,
\item $A_i\cup (C_i+\rho_i)=A_{i+1}\cup C_{i+1}$ for all $1\leq i<r$,
\item $L=A_r\cup(C_r+\rho_r)$,
\item $\gv(A_i\cup C_i)=\gv(L)$ for all $1\leq i\leq r$,
\item $h_{\langle X^{A_i\cup C_i}\rangle_S^\sat}<h_{\langle X^{A_i\cup (C_i+\rho_i)}\rangle_S^\sat}$ for all $1\leq i\leq r$.
\end{enumerate}
\end{proposition}

\begin{proof}
Assume that $B\neq L$. The condition $\gv(B)=\gv(L)$ implies $\#B=\#L$. Let $c:=\min_\deglex B$ and $c':=\max_\deglex\N^n_d\setminus B$. Then 
\begin{enumerate}
\item $c'>_\deglex c$,
\item $c'_n>c_n$,
\item $c\in B\setminus L$,
\item $c'\in L\setminus B$,
\item $m(c)=m(c')=n$.
\end{enumerate}

Proof of (i): If $c>_\deglex c'$, then $B$ is a lexicographic set and hence equals $L$.

Proof of (ii): Assume that $c'_n\leq c_n$. Let $b:=c+(c_n-c'_n)(e_{n-1}-e_n)$. Since $B$ is a Borel set, we have $b\in B(c'_n)$. Thus, because $B(c'_n)$ is a lexicographic segment of $\N^n_d(c'_n)$ and $c'\notin B(c'_n)$, it holds $c'<_\deglex b$. From $b_n=c'_n$ and $|b|=|c'|$ it follows that $\mu:=\mu(b-c')<n-1$. Therefore, we have $c_i=b_i=c'_i$ for all $1\leq i<\mu$ and $c_\mu=b_\mu>c'_\mu$, whence $c>_\deglex c'$, which contradicts the previous statement.

Proof of (iii): If $c\in L$, then $b\in L$ for all $b\in B$, because $L$ is lexicographic and $c\in B$ is minimal with respect to the lexicographic order. But $B\subset L$ contradicts $\#B=\#L$ and $B\neq L$.

Proof of (iv): The set $B':=\{b\in B\mid b>_\deglex c'\}$ is lexicographic by the definition of $c'$. If $c'\notin L$, then $c'<_\deglex l$ for all $l\in L$, since $L$ is lexicographic. It follows that $L\subset B'\subset B$, which is a contradiction.

Proof of (v): Since $B$ and $L$ are both growth-height-lexico\-graphic and since $\gv(B)=\gv(L)$, it holds $\{b\in B\mid m(b)<n\}=\{l\in L\mid m(l)<n\}$. By the previous three statements we have $c'\in L\setminus(B\cup\{l\in L\mid m(l)<n\})$ and $c\in B\setminus L$, whence $c\in B\setminus\{b\in B\mid m(b)<n\}$.

\sk
Set $A:=B\setminus\{c\}$, $C:=\{c\}$ and $\rho:=c'-c$. We claim that $(A,C,\rho)$ is a binomial system. We have to show that $A\cup\{c'\}$ is a Borel set. By Lemma \ref{Borel greater implies greater for any term order} it is clear that $A$ is a Borel set. Let now $b\in\N^n_d$ be such that $b\borelg c'$. By Lemma \ref{Borel greater implies greater for any term order} it holds $b>_\deglex c'$. Hence, by the definition of $c'$, we have $b\in B$. We know that $c'>_\deglex c$, whence $b\neq c$, and the claim is proved. Moreover, $A\cup(C+\rho)=A\cup\{c'\}$ is growth-height-lexico\-graphic.

Now it is easy to see that $(A,C,\rho)$ is Mall. In fact, from $c'_n>c_n$ it follows that $(A,C,\rho)$ is admissible. So $(A,C,\rho)$ is good, as $C$ has only one element. From $c'>_\deglex c$ it follows that $c'_{\mu(c'-c)}>c_{\mu(c'-c)}$, whence $\rho_{\mu(\rho)}=c'_{\mu(c'-c)}-c_{\mu(c'-c)}>0$. Furthermore,  $m(c)=m(c')=m(c+\rho)$. Thus, $(A,C,\rho)$ is Mall.

We have $\gv(A\cup(C+\rho))=\gv(A\cup C)$, and for $1\leq i\leq d$
$$\hv(A\cup(C+\rho))_i=
\begin{cases}
\hv(A\cup C)_i-1,&\text{if $i=c_n$;}\\
\hv(A\cup C)_i+1,&\text{if $i=c'_n$;}\\
\hv(A\cup C)_i,&\text{otherwise.}
\end{cases}$$
Applying the formula of Proposition \ref{computation of the Hilbert function} b) we compute
$$h_{\langle X^{A\cup (C+\rho)}\rangle_S^\sat}(i)=
\begin{cases}
h_{\langle X^{A\cup C}\rangle_S^\sat}(i),&\text{if $0\leq i<d-c'_n$;}\\
h_{\langle X^{A\cup C}\rangle_S^\sat}(i)+1,&\text{if $d-c'_n\leq i<d-c_n$;}\\
h_{\langle X^{A\cup C}\rangle_S^\sat}(i),&\text{if $i\geq d-c_n$.}\\
\end{cases}$$
This shows that $h_{\langle X^{A\cup C}\rangle_S^\sat}<h_{\langle X^{A\cup (C+\rho)}\rangle_S^\sat}.$

Recall that $c\in B\setminus L$ and $c'\in L\setminus B$. Hence, setting $(A_1,C_1,\rho_1):=(A,C,\rho)$, we conclude by induction on the number $\#(L\setminus B)$.
\end{proof}
\sk
\subsection{The generic initial ideal of a binomial ideal}\label{subsection generic initial ideal}
Let $(A,C,\rho)$ be a good binomial system and $F(A,C,\rho)$ the induced binomial ideal. Then the generic initial ideal of $F(A,C,\rho)$ and the initial ideal of $F(A,C,\rho)$ with respect to any admissible term order coincide (Theorem \ref{Gin = in for good binomial ideals}). This conclusion is not at all trivial: It does not hold in general if $(A,C,\rho)$ is only admissible, but not good (s.~Example \ref{counterexample}). The crucial point is that $F(A,C,\rho)$ is fixed under the action of the unipotent group if $(A,C,\rho)$ is good (Proposition \ref{good binomial ideals are fixed by unipotents}). 

In case the binomial system  $(A,C,\rho)$ is good, we want to compute $g(f)$ if $f$ is a generator of the binomial ideal $F(A,C,\rho)$ and $g\in\unipotent$ is unipotent. To do this, we introduce generic coordinates for $g$: 
\begin{notation}
Set $T:=K[Y_{ij}\mid 1\leq i\leq j\leq n]$.

Define the automorphism of $T$-algebras 
$$\phi:T[X_1,\dots,X_n]\ra T[X_1,\dots,X_n]$$ 
by $\phi(X_j):=\sum_{i=1}^j Y_{ij}X_i$.

Let $g=[g_{ij}\mid 1\leq i\leq n,\:1\leq j\leq n]\in K^{(n,n)}$ be a matrix. Then we denote by $\bar{g}:T[X_1,\dots,X_n]\ra S$ the homomorphism of $S$-algebras defined by $\bar{g}(Y_{ij}):=g_{ij}$ for $1\leq i\leq j\leq n$.

For $a\in\N^n$ set $g^a:=\prod_{i=1}^n g_{ii}^{a_i}$ and $Y^a:=\prod_{i=1}^n Y_{ii}^{a_i}$.

For $M\in\N^{(n,n)}$ set $Y^M:=\prod_{1\leq i\leq j\leq n}Y_{ij}^{M_{ij}}$.

For $\rho\in\Z^n$ and $M\in\Z^{(n,n)}$ let $M+\rho\in\Z^{(n,n)}$ denote the matrix which is defined by
$$(M+\rho)_{ij}:=
\begin{cases}
M_{ij},&\text{if $i\neq j$;}\\
M_{ii}+\rho_i,&\text{otherwise.}
\end{cases}$$

For $a,b\in\N^n$ set 
$$U(a,b):=\{M\in U(n)\mid \sum_{i=1}^n M_{ji}=a_j,\sum_{i=1}^n M_{ij}=b_j\ \forall\ 1\leq j\leq n\}.$$
\end{notation}

To introduce more notations, we need the following Lemma:

\begin{lemma}\label{lemma on Borel matrices}
Let $b,c\in\N^n_d$ and $\rho\in\Z^n$ be such that $b+\rho$, $c+\rho\in\N^n_d$ and such that $b_i=c_i$ for all $1\leq i<m(\rho)$.

a) Let $b\borelgeq c$ and $M\in U(b,c)$. Then $M_{jj}=c_j$ for all $1\leq j\leq m(\rho)$.

b) Let $b\borelgeq c$ and $M\in U(b,c)$. Then $M_{jj}+\rho_j\geq 0$ for all $1\leq j\leq n$.

c) Let $M\in\Z^{(n,n)}$. Then $M\in U(b,c)$ if and only if $M+\rho\in U(b+\rho,c+\rho)$.
\end{lemma}

\begin{proof}
a) We use an inductive argument. It is clear that $M_{11}=\sum_{i=1}^n M_{i1}=c_1$. Let $1<j\leq m(\rho)$ and assume that $c_k=M_{kk}$ for all $1\leq k<j$. We have $M_{kk}=c_k=b_k=\sum_{i=1}^n M_{ki}$ for all $1\leq k<j$ and therefore $M_{kj}=0$ for all $1\leq k<j$. Since $M_{ij}=0$ for all $i>j$, it holds $c_j=\sum_{i=1}^n M_{ij}=M_{jj}$.

\sk
b) The statement follows from part a) and from the condition $c_j+\rho_j\geq 0$ for all $1\leq j\leq n$.

\sk
c) Let $M\in U(b,c)$. Then part b) states that $M+\rho\in U(n)$. Since
$$(b+\rho)_j=b_j+\rho_j=\sum_{i=1}^n M_{ji}+\rho_j=\sum_{i=1}^n(M+\rho)_{ji}$$
and
$$(c+\rho)_j=c_j+\rho_j=\sum_{i=1}^n M_{ij}+\rho_j=\sum_{i=1}^n(M+\rho)_{ij}$$
for all $1\leq j\leq n$, it holds $M+\rho\in U(b+\rho,c+\rho)$.

In order to prove the converse implication, one just has to replace $b$, $c$, $\rho$ with $b+\rho$, $c+\rho$, $-\rho$ respectively.
\end{proof}

\begin{remark}\label{remark on the diagonal of Borel matrices}
Part b) and hence part c) of the previous Lemma are also true if $b$, $c$ and $\rho$ miss the preconditioned property $b_i=c_i$ for all $1\leq i<m(\rho)$. To prove this, one can use the following fact:

\nopagebreak
\it{Let $b,c\in \N^n_d$ with $b\borelgeq c$ and let $M\in U(b,c)$. Then
$$\forall\ 1\leq j\leq n\ \exists a\in\N^n_d:b\borelgeq a\borelgeq c\text{ and }M_{jj}=a_j.$$}
\end{remark}

\begin{notation}
For $m$ and $m_1,\dots,m_n\in\N$ with $m=\sum_{i=1}^n m_i$ set 
$$\binom{m}{m_1,\dots,m_n}:=\frac{m!}{m_1!\cdots m_n!}.$$

For $M\in\N^{(n,n)}$ set
$$\mu_M:=\prod_{j=1}^n\binom{\sum_{i=1}^n M_{ij}}{M_{1j},\dots,M_{nj}}.$$

For $a,b\in\N^n_d$ let $\alpha_a^b:=c^{\phi(X^b)}_{X^a}\in T$ be the coefficient of $X^a$ in the polynomial $\phi(X^b)$. (We have $\phi(X^b)=\sum_{a\in\N^n_d}\alpha^b_a X^a=\sum_{a\in\supp_T(\phi(X^b))}\alpha^b_a X^a$).

For $b,c\in\N^n_d$ and $\rho\in\Z^n$ such that $b+\rho$, $c+\rho\in\N^n_d$ and $b_i=c_i$ for all $1\leq i<m(\rho)$ set 
$$p^\rho_{b,c}:=\sum_{M\in U(b,c)}\mu_M Y^{M-\rho^-}.$$
(By Lemma \ref{properties of Borel order} and Lemma \ref{lemma on Borel matrices} b) this is a polynomial in $T[X_1,\dots,X_n]$).
\end{notation}

In the proof of Proposition \ref{good binomial ideals are fixed by unipotents} we shall make use of the following two Lemmata:

\begin{lemma}\label{computation of alpha^b_a}
Let $a,b\in\N^n_d$. 

a) Then $\alpha^b_a=\sum_{M\in U(a,b)}\mu_M Y^M$.

b) Then $\alpha^b_b=Y^b$.

c) Let $\rho\in\Z^n$ be such that $b+\rho\in\N^n_d$. Then  $p^\rho_{b,b}=Y^{b-\rho^-}$.

d) Then $a\borelgeq b\iff\alpha^b_a\neq 0\iff a\in\supp_T(\phi(X^b))$.
\end{lemma}

\begin{proof}
a) From the multinomial formula it follows
\begin{align*}
\sum_{a\in\supp_T(\phi(X^b))}\alpha^b_a X^a&=\phi(X^b)=\prod_{j=1}^n \bigl(\sum_{i=1}^j Y_{ij}X_i\bigr)^{b_j}\\
&=\prod_{j=1}^n\sum_{\substack{k_1,\,\dots,\,k_j\in\N\\k_1+\dots+k_j=b_j}}\binom{b_j}{k_1,\dots,k_j}\prod_{i=1}^j (Y_{ij}X_i)^{k_i}\\
&=\sum_{\substack{k_{11},\,k_{12},\,k_{22},\,\dots,\,k_{1n},\,\dots,\,k_{nn}\in\N\\\sum_{i=1}^j k_{ij}=b_j\ \forall\ 1\leq j\leq n}}\ \prod_{j=1}^n\binom{b_j}{k_{1j},\dots,k_{jj}}\prod_{i=1}^j (Y_{ij}X_i)^{k_{ij}}\\
&=\sum_{\substack{M\in U(n)\\\sum_{i=1}^n M_{ij}=b_j\ \forall\ 1\leq j\leq n}}\prod_{j=1}^n\binom{b_j}{M_{1j},\dots,M_{nj}}\prod_{i=1}^j (Y_{ij}X_i)^{M_{ij}}\\
&=\sum_{\substack{M\in U(n)\\\sum_{i=1}^n M_{ij}=b_j\ \forall\ 1\leq j\leq n}}\mu_M\prod_{1\leq i\leq j\leq n} Y_{ij}^{M_{ij}}X_i^{M_{ij}}\\
&=\sum_{\substack{M\in U(n)\\\sum_{i=1}^n M_{ij}=b_j\ \forall\ 1\leq j\leq n}}\mu_M Y^M\prod_{i=1}^n X_i^{\sum_{j=1}^n M_{ij}}\\
&=\sum_{\substack{M\in U(n)\\\sum_{i=1}^n M_{ij}=b_j\ \forall\ 1\leq j\leq n}}\mu_M Y^M X^{(\sum_{j=1}^n M_{1j},\dots,\sum_{j=1}^n M_{nj})}.
\end{align*}
Hence we get 
$$\alpha^b_a=\sum_{\substack{M\in U(n)\\\sum_{i=1}^n M_{ij}=b_j\ \forall\ 1\leq j\leq n\\\sum_{j=1}^n M_{ij}=a_i\ \forall\ 1\leq i\leq n}} \mu_M Y^M=\sum_{M\in U(a,b)}\mu_M Y^M.$$

\sk
b) and c): The set $U(b,b)$ contains only one matrix $M\in\N^{(n,n)}$, defined by 
$$M_{ij}:=
\begin{cases}
b_j,&\text{if $i=j$;}\\
0,&\text{otherwise.}
\end{cases}$$
It holds $\mu_M=\prod_{j=1}^n\binom{b_j}{M_{1j},\dots,M_{nj}}=\prod_{j=1}^n\binom{b_j}{b_j}=1$ and $Y^M=\prod_{1\leq i\leq j\leq n} Y_{ij}^{M_{ij}}=\prod_{j=1}^n Y_{jj}^{b_j}=Y^b$. Therefore, it follows that $p^\rho_{b,b}=\mu_M Y^{M-\rho^-}=Y^{b-\rho^-}$
and, by part a) of the Lemma, $\alpha^b_b=\mu_M Y^M=Y^b$.

\sk
d) We notice that $\mu_M>0$ for all $M\in\N^{(n,n)}$. By part a) and by Lemma \ref{properties of Borel order} we have $$a\in\supp_T(\phi(X^b))\iff \alpha^b_a\neq 0\iff U(a,b)\neq\emptyset\iff a\borelgeq b.$$
\end{proof}

\begin{lemma}\label{computation of alpha^c_b}
Let $b,c\in\N^n_d$ and $\rho\in\Z^n$ be such that $b+\rho$, $c+\rho\in\N^n_d$ and $b_i=c_i$ for all $1\leq i<m(\rho)$. Then  $\alpha^c_b=p^\rho_{b,c}\,Y^{\rho^-}$\! and $\alpha^{c+\rho}_{b+\rho}=p^\rho_{b,c}\,Y^{\rho^+}$.
\end{lemma}

\begin{proof}
If $b\not\borelgeq c$, then $\alpha^c_b=\alpha^{c+\rho}_{b+\rho}=0$ by Lemma \ref{computation of alpha^b_a} d) and $p^\rho_{b,c}=0$ by Lemma \ref{properties of Borel order}. Hence, we can assume that $b\borelgeq c$.
The first equation follows immediately from Lemma \ref{computation of alpha^b_a} a). To prove the second one, we first show that $\mu_M=\mu_{M+\rho}$ for all $M\in U(b,c)$. Let $M\in U(b,c)$. Let $1\leq j\leq n$. If $j\leq m(\rho)$ we have by part a) of Lemma \ref{lemma on Borel matrices} that
$$\binom{c_j}{M_{1j},\dots,M_{nj}}=\frac{c_j!}{M_{jj}!}=1=\frac{(c_j+\rho_j)!}{(M_{jj}+\rho_j)!}=\binom{c_j+\rho_j}{(M+\rho)_{1j},\dots,(M+\rho)_{nj}}.$$
If $j>m(\rho)$ we have $\binom{c_j}{M_{1j},\dots,M_{nj}}=\binom{c_j+\rho_j}{(M+\rho)_{1j},\dots,(M+\rho)_{nj}}$. The claim now follows from the definition of $\mu_M$ and $\mu_{M+\rho}$.

By Lemma \ref{computation of alpha^b_a} a) and \ref{lemma on Borel matrices}~c) we get
$$\alpha ^{c+\rho}_{b+\rho}=\sum_{M\in U(b+\rho,c+\rho)}\mu_{M+\rho}Y^{M+\rho}=\sum_{M\in U(b,c)}\mu_M Y^{M-\rho^-+\rho^+}=p^\rho_{b,c}\,Y^{\rho^+}.$$
\end{proof}

\begin{proposition}\label{good binomial ideals are fixed by unipotents}
Let $(A,C,\rho)$ be a good binomial system and $g\in\unipotent$ an unipotent matrix. Then $g(F(A,C,\rho))=F(A,C,\rho).$
\end{proposition}

\begin{proof}
If $C=\emptyset$, the statement follows from the fact that the ideal $\langle X^A\rangle_S$ is Borel-fixed (s.~ Proposition~\ref{existence of gin}). 

Now let $C\neq\emptyset$ and $c:=\min_\degrevlex C$. Then $(A,C\setminus\{c\},\rho)$ is again a good binomial system (cf.\ Lemma \ref{Borel greater implies greater for any term order}). Since $\langle g(G)\rangle_S=g(\langle G\rangle_S)$ for any set $G\subset S$, it is enough to show that
$$\langle g(X^A\cup\Bin(C,\rho))\rangle_K=\langle X^A\cup\Bin(C,\rho)\rangle_K.$$ 
Hence it is enough to show that 
$$g(X^A\cup\Bin(C,\rho))\subset\langle X^A\cup\Bin(C,\rho)\rangle_K.$$
By induction we may assume that $g(F(A,{C\setminus\{c\}},\rho))=F(A,{C\setminus\{c\}},\rho)$, whence 
$$g(X^A\cup\Bin(C\setminus\{c\},\rho))\subset\langle X^A\cup\Bin(C,\rho)\rangle_K.$$ 
Therefore, it is enough to show that 
$$g(X^c-X^{c+\rho})\in\langle X^A\cup\Bin(C,\rho)\rangle_K.$$
Observe that $g^a=1$ for all $a\in\N^n$ as $g$ is unipotent.
Using Lemma \ref{computation of alpha^b_a} and \ref{computation of alpha^c_b} we compute
\begin{align*}
g(X^c-X^{c+\rho})=&\bar{g}(\phi(X^c-X^{c+\rho}))\\
=&\bar{g}\bigl(\sum_{a\in\supp_T(\phi(X^c))}\hspace{-12pt}\alpha^c_a X^a - \hspace{-12pt}\sum_{a\in\supp_T(\phi(X^{c+\rho}))}\hspace{-12pt}\alpha^{c+\rho}_a X^a\bigr)\\
=&\bar{g}\bigl(\sum_{a\borelgeq c}\alpha^c_a X^a - \hspace{-8pt}\sum_{a\borelgeq c+\rho}\hspace{-8pt} \alpha^{c+\rho}_a X^a\bigr)\\
=&\bar{g}\bigl(\!\sum_{\substack{a\in A\\a\borelg c}}\!\alpha^c_a X^a -\hspace{-8pt}\sum_{\substack{a\in A\\a\borelg c+\rho}}\hspace{-8pt} \alpha^{c+\rho}_a X^a\bigr)+\bar{g}\bigl(\!\sum_{\substack{b\in C\\b\borelgeq c}}\!\alpha^c_b X^b - \hspace{-8pt}\sum_{\substack{b\in C+\rho\\b\borelgeq c+\rho}}\hspace{-8pt} \alpha^{c+\rho}_b X^b\bigr)\\
=&\bar{g}\bigl(\sum_{a\in A}\alpha^c_a X^a - \sum_{a\in A}\alpha^{c+\rho}_a X^a\bigr)+\bar{g}\bigl(\sum_{b\in C}(\alpha^c_b X^b - \alpha^{c+\rho}_{b+\rho} X^{b+\rho})\bigr)\\
=&\bar{g}\bigl(\sum_{a\in A}(\alpha^c_a-\alpha^{c+\rho}_a)X^a\bigr)+\bar{g}\bigl(\sum_{b\in C} (p^\rho_{b,c}\,Y^{\rho^-}X^b -p^\rho_{b,c}\,Y^{\rho^+} X^{b+\rho})\bigr)\\
=&\sum_{a\in A}\bar{g}(\alpha^c_a-\alpha^{c+\rho}_a)X^a+\sum_{b\in C} \bar{g}(p^\rho_{b,c})\bigl(\bar{g}(Y^{\rho^-})X^b -\bar{g}(Y^{\rho^+}) X^{b+\rho}\bigr)\\
=&\sum_{a\in A}\bar{g}(\alpha^c_a-\alpha^{c+\rho}_a)X^a+ \sum_{b\in C} \bar{g}(p^\rho_{b,c})(X^b - X^{b+\rho}).
\end{align*}
Hence our Proposition is proved.
\end{proof}

Now, it follows immediately from Proposition \ref{good binomial ideals are fixed by unipotents} and Corollary \ref{Gin=in for unipotent fixed ideals}:

\begin{theorem}\label{Gin = in for good binomial ideals}
Let $(A,C,\rho)$ be a good binomial system and $\tau$ an admissible term order. Then $\Gin_\tau F(A,C,\rho)=\init_\tau F(A,C,\rho).$\hfill$\square$
\end{theorem}

As a further consequence we have by Proposition \ref{Gin and sat commute} and \ref{Groebner bases of (saturated) binomial ideals}:

\begin{corollary}\label{gin and sat commute for binomial ideals}
Let $(A,C,\rho)$ be a good binomial system. Then 
$$\Gin_\degrevlex (F(A,C,\rho)^\sat)=\langle X^{A\cup C}\rangle_S^\sat=\langle X^{A^*\cup C^*}\rangle_S.$$
Furthermore, if $(A,C,\rho)$ is Mall and $\rho_n=0$, then 
\begin{equation*}\init_\deglex (F(A,C,\rho)^\sat)=\langle X^{A\cup(C+\rho)}\rangle_S^\sat=\langle X^{A^*\cup(C+\rho)^*}\rangle_S.\tag*{$\square$}\end{equation*}
\end{corollary}

\pagebreak[2]
We conclude this subsection with an example which shows that Theorem \ref{Gin = in for good binomial ideals} fails if the binomial system $(A,C,\rho)$ is not good:
\begin{example}\label{counterexample}
Consider the ring $R:=K[X_1,\dots,X_5]=K[x,y,z,t,u]$. Let $\rho:=(1,-2,2,-2,1)$,  $b:=(0,2,0,3,0)$, $c:=(0,2,0,2,1)$ and $C:=\{b,c\}$. Let $B\subset\N^5_5$ be the smallest Borel set containing $D:=C\cup(C+\rho)$ and set $A:=B\setminus D$. Then $(A,C,\rho)$ is an admissible binomial system. Let $\ideala:=F(A,C,\rho)$. We then compute:\pagebreak[1]
\begin{align*}
\Gin_\degrevlex\ideala=\langle\{&
x^5, x^4y, x^3y^2, x^2y^3, xy^4, y^5, x^4z, x^3yz, x^2y^2z, xy^3z, y^4z, x^3z^2,\\* &x^2yz^2, xy^2z^2, y^3z^2, x^2z^3, xyz^3, y^2z^3, xz^4, x^4t, x^3yt, x^2y^2t,\\* &xy^3t, y^4t, x^3zt, x^2yzt, xy^2zt, y^3zt, x^2z^2t, xyz^2t, y^2z^2t, xz^3t,\\* &x^3t^2, x^2yt^2, xy^2t^2, y^3t^2, x^2zt^2, xyzt^2, y^2zt^2, xz^2t^2, x^2t^3, xyt^3, y^2t^3,\\* &x^4u, x^3yu, x^2y^2u, xy^3u, y^4u, x^3zu, x^2yzu, xy^2zu, y^3zu, x^2z^2u,\\* & xyz^2u, y^2z^2u, xz^3u, x^3tu, x^2ytu, xy^2tu, y^3tu, x^2ztu, xyztu,\\* &y^2ztu, \underline{xz^2tu}, x^2t^2u, xyt^2u, x^3u^2, x^2yu^2, xy^2u^2, x^2zu^2, xyzu^2
\}\rangle_{R},\\
\ &\\
\init_\degrevlex\ideala=\langle\{&
x^5, x^4y, x^3y^2, x^2y^3, xy^4, y^5, x^4z, x^3yz, x^2y^2z, xy^3z, y^4z, x^3z^2,\\* &x^2yz^2, xy^2z^2, y^3z^2, x^2z^3, xyz^3, y^2z^3, xz^4, x^4t, x^3yt, x^2y^2t,\\* &xy^3t, y^4t, x^3zt, x^2yzt, xy^2zt, y^3zt, x^2z^2t, xyz^2t, y^2z^2t, xz^3t,\\* &x^3t^2, x^2yt^2, xy^2t^2, y^3t^2, x^2zt^2, xyzt^2, y^2zt^2, xz^2t^2, x^2t^3, xyt^3, y^2t^3,\\* &x^4u, x^3yu, x^2y^2u, xy^3u, y^4u, x^3zu, x^2yzu, xy^2zu, y^3zu, x^2z^2u,\\ &xyz^2u, y^2z^2u, xz^3u, x^3tu, x^2ytu, xy^2tu, y^3tu, x^2ztu, xyztu,\\* &y^2ztu, x^2t^2u, xyt^2u, \underline{y^2t^2u}, x^3u^2, x^2yu^2, xy^2u^2, x^2zu^2, xyzu^2
\}\rangle_{R}.
\end{align*}
These ideals are not equal, as is indicated by the underlined generators. The reason is the following: Let $M\in\N^{(5,5)}$ be the unique element of $U(b,c)$. Then $M_{55}=0\neq c_5$, but $m(\rho)=5$ (counterexample to part (a) of Lemma~\ref{lemma on Borel matrices}). We cannot conclude that $\mu_M$ equals $\mu_{M+\rho}$; indeed $\mu_M=1$ and $\mu_{M+\rho}=2$. It follows that $\alpha_{b+\rho}^{c+\rho}=2p^\rho_{b,c}Y^{\rho^+}$ (counterexample to Lemma~\ref{computation of alpha^c_b}). Let $g\in\Gl(5,K)$ be unipotent. We then compute 
\begin{multline*}
g(X^c-X^{c+\rho})=X^c - X^{c+\rho} + \bar{g}(p^\rho_{b,c})(X^b -2X^{b+\rho})+ \sum_{a\in A}\bar{g}(\alpha^c_a-\alpha^{c+\rho}_a)X^a\\*
\notin\langle X^A\cup\Bin(C,\rho)\rangle_K
\end{multline*}
(counterexample to Proposition \ref{good binomial ideals are fixed by unipotents}). A further computation yields 
$$g(X^b-X^{b+\rho})=X^b - X^{b+\rho} + \sum_{a\in A}\bar{g}(\alpha^b_a-\alpha^{b+\rho}_a)X^a.$$ 
Since $A$ is a Borel set, it is clear that $X^A\subset g(F(A,C,\rho))$, whence $X^b-X^{b+\rho}$, $X^b -2X^{b+\rho} + \bar{g}(p^\rho_{b,c})^{-1}(X^c - X^{c+\rho})\in g(F(A,C,\rho))$. Since $b>_\degrevlex c$ we get $X^{b+\rho}\in\Gin_\degrevlex F(A,C,\rho)\setminus\init_\degrevlex F(A,C,\rho)$ and $X^c\in\init_\degrevlex F(A,C,\rho)\setminus\Gin_\degrevlex F(A,C,\rho)$.
\end{example}
\sk
\subsection{Binomial ideals and sequentially Cohen-Macaulay\-ness}\label{subsection seq. CM}
In \cite[2.2]{HS} J.~Herzog and E.~Sbarra showed that in characteristic zero the $S$-module $S/\idealb$ is sequentially Cohen-Macaulay if $\idealb$ is a Borel ideal. Later J.~Herzog, D.~Popescu, and M.~Vladoiu \cite{HPV} generalized this result to monomial ideals of Borel type in any characteristic of $K$. An ideal $\idealb\subset S$ is of Borel type if $(\idealb:X_j^\infty)=(\idealb:_S\langle X_1,\dots,X_j\rangle_S^\infty)$ for all $1\leq j\leq n$. It is well known that Borel-fixed ideals are of Borel type (\cite[15.24]{E}). Since a binomial ideal $F(A,C,\rho)$ is fixed under the action of the unipotent group, if the binomial system $(A,C,\rho)$ is good (Proposition \ref{good binomial ideals are fixed by unipotents}), it is natural to ask whether $S/F(A,C,\rho)$ is sequentially Cohen-Macaulay. In this section we prove that in fact $S/F(A,C,\rho)$ and $S/F(A,C,\rho)^\sat$ are sequentially Cohen-Macaulay if $(A,C,\rho)$ is an admissible binomial system (Theorem \ref{binomial ideals and sCM}).

\begin{definition}
Let $R$ be a graded ring. A finitely generated graded $R$-module $M$ is \emph{sequentially Cohen-Macaulay} if there exists a finite filtration
$$0=M_0\subsetneq M_1\subsetneq\dots\subsetneq M_r=M$$ of $M$ by graded submodules such that
\begin{enumerate}
\item $M_i/M_{i-1}$ is Cohen-Macaulay for all $1\leq i\leq r$,
\item $\dim(M_i/M_{i-1})<\dim(M_{i+1}/M_i)$ for all $1\leq i<r$.
\end{enumerate}
\end{definition}

\sk
Provided that $(A,C,\rho)$ is an admissible binomial system, we construct a filtration $F(A,C,\rho)^\sat=F_0\subset\dots\subset F_{n-1}=S$ of ideals such that the quotients $F_i/F_{i-1}$ are zero or Cohen-Macaulay of dimension $i$ for all $1\leq i<n$ (Proposition \ref{Binomial quotients are CM}). There is a natural way to define the ideals $F_i$:
\begin{notation}
For $i\in\N$ let $S_{(i)}$ denote the polynomial ring $K[X_1,\dots,X_i]$ and let $\terms[i]$ be the set of monomials of $S_{(i)}$. In particular,  $S_{(n)}=S$ and $\terms[n]=\terms$. For $i<j$ one has a canonical inclusion $S_{(i)}\subset S_{(j)}$. 

For $A\subset\N^n$ and $1\leq i\leq n$ set 
$$A_i:=\{(a_1,\dots,a_i)\in\N^i\mid(a_1,\dots,a_i,0,\dots,0)\in A\}.$$

Let $(A,C,\rho)$ be a binomial system. For $0\leq i\leq n-m(\rho)$ set
$$F_i(A,C,\rho):=(F(A,C,\rho)\cap S_{(n-i)})^\sat S.$$
For $n-m(\rho)<i<n$ set
$$F_i(A,C,\rho):=\bigl((\init_\degrevlex F(A,C,\rho))\cap S_{(n-i)}\bigr)^\sat S.$$
\end{notation}

\begin{lemma}\label{generators of retractions}
Let $1\leq l\leq n$ and let $G\subset S$ be a set of polynomials such that $g\in S_{(l)}$ or\/ $\supp(g)\cap S_{(l)}=\emptyset$ for all $g\in G$. Then $\langle G\cap S_{(l)}\rangle_{S_{(l)}}=\langle G\rangle_S\cap S_{(l)}$.
\end{lemma}

\begin{proof}
The inclusion ``$\subset$'' is obvious. We have to show that $\langle G\rangle_S\cap S_{(l)}\subset\langle G\cap S_{(l)}\rangle_{S_{(l)}}$. Let $f\in\langle G\rangle_S\cap S_{(l)}$. Then there are monomials $u_1$, \dots, $u_r\in\terms$ and polynomials $g_1$, \dots, $g_r\in G$ such that $f=\sum_{i=1}^r u_i g_i$. For all $1\leq i\leq r$ by our assumption $u_i g_i\in S_{(l)}$ or $\supp(u_i g_i)\cap S_{(l)}=\emptyset$. Hence it follows that
$$f=\sum_{\substack{1\leq i\leq r\\ u_i g_i\in S_{(l)}}} u_i g_i\in\langle G\cap S_{(l)}\rangle_{S_{(l)}}.\vspace{-6mm}$$
\end{proof}

The following Lemma is crucial, because it allows us to compute the generators and the initial ideals of all ideals $F_i(A,C,\rho)$.
\begin{lemma}\label{Groebner basis of F_i(A,C,rho)}
Let $(A,C,\rho)$ be an admissible binomial system and $m(\rho)\leq i\leq n$. Then 
$$X^{(A_i)^*}\cup\Bin((C_i)^*,(\rho_1,\dots,\rho_i))$$ 
is a Gr\"obner basis of $(F(A,C,\rho)\cap S_{(i)})^\sat$ with respect to the reverse lexicographic order of\/ $\terms[i]$.
\end{lemma}

\begin{proof}
Since $(A_i,C_i,(\rho_1,\dots,\rho_i))$ is an admissible binomial system in $i$ indeterminates, by Proposition \ref{Groebner bases of (saturated) binomial ideals}~c) the set $X^{(A_i)^*}\cup\Bin((C_i)^*,(\rho_1,\dots,\rho_i))$ is a Gr\"obner basis of $F(A_i,C_i,(\rho_1,\dots,\rho_i))^\sat$ with respect to the reverse lexicographic order of $\terms[i]$. To complete the proof it is enough to show that $F(A,C,\rho)\cap S_{(i)}=F(A_i,C_i,(\rho_1,\dots,\rho_i))$.

Since $m(\rho)\leq i$, it holds $g\in S_{(i)}$ or $\supp(g)\cap\nolinebreak S_{(i)}=\emptyset$ for all $g\in X^A\cup\Bin(C,\rho)$. By Lemma \ref{generators of retractions} we get 
\begin{align*}
F(A,C,\rho)\cap S_{(i)}&=\langle X^A\cup\Bin(C,\rho)\rangle_S\cap S_{(i)}
=\langle (X^A\cup\Bin(C,\rho))\cap S_{(i)}\rangle_{S_{(i)}}\\
&=\langle X^{A_i}\cup\Bin(C_i,(\rho_1,\dots,\rho_i))\rangle_{S_{(i)}}
=F(A_i,C_i,(\rho_1,\dots,\rho_i)).
\end{align*}
\end{proof}

\begin{lemma}\label{lemma about F_i(A,C,rho)}
Let $(A,C,\rho)$ be an admissible binomial system and put $m:=m(\rho)$. Then

a) $F_i(A,C,\rho)=(F_i(A,C,\rho)\cap S_{(n-i-1)})S\ \  \forall\,i\in\{0,\dots,n-1\}\setminus\{n-m\}$,

b) $F_{i+1}(A,C,\rho)=(F_i(A,C,\rho)\cap S_{(n-i-1)})^\sat S\ \  \forall\,i\in\{0,\dots,n-2\}\!\setminus\!\{n-m\}$,

c) $(F(A,C,\rho)\cap S_{(m)})^\sat\subset\bigl((\init_\degrevlex F(A,C,\rho))\cap S_{(m-1)}\bigr)^\sat S_{(m)},$

d) $F_i(A,C,\rho)\subset F_{i+1}(A,C,\rho)$ for all\/ $0\leq i<n-1$. 

e) $\{X_1,\dots,X_{m-1}\}\subset$
$$\sqrt{\Bigl((F(A,C,\rho)\cap S_{(m)})^\sat\smash{\underset{S_{(m)}}{:}}\bigl((\init_\degrevlex F(A,C,\rho))\cap S_{(m-1)}\bigr)^\sat S_{(m)}\Bigr)}.$$

f) Let $u\in\terms[m]$. Then 
$u\notin\init_\degrevlex(F(A,C,\rho)\cap S_{(m)})^\sat$ implies $X_m u\notin\init_\degrevlex(F(A,C,\rho)\cap S_{(m)})^\sat.$
\end{lemma}

\smallskip
\begin{proof}
a) Let $i\in\{0,\dots,n-1\}\setminus\{n-m\}$. If $i>n-m$, then by Remark \ref{Groebner bases of Borel ideals} $F_i(A,C,\rho)=\bigl((\init_\degrevlex F(A,C,\rho))\cap S_{(n-i)}\bigr)^\sat S$ is generated in $S_{(n-i-1)}$ since $(\init_\degrevlex F(A,C,\rho))\cap S_{(n-i)}\subset S_{(n-i)}$ is a Borel ideal by Proposition \ref{Groebner bases of (saturated) binomial ideals}~a). On the other hand, if $n-i>m$, then by Lemma \ref{Groebner basis of F_i(A,C,rho)} we have
\begin{align*}
F_i(A,C,\rho)&=(F(A,C,\rho)\cap S_{(n-i)})^\sat S\\*
&=\bigl\langle X^{(A_{n-i})^*}\cup\Bin((C_{n-i})^*,(\rho_1,\dots,\rho_{n-i}))\bigr\rangle_{S_{(n-i)}} S.
\end{align*}
Since $\rho_{n-i}=0$, therefore $F_i(A,C,\rho)$ is also generated in $S_{(n-i-1)}$.

\sk
b) Let $i\in\{0,\dots,n-2\}\setminus\{n-m\}$. Set 
$$\ideala:=
\begin{cases}
F(A,C,\rho),&\text{if $i<n-m$;}\\
\init_\degrevlex F(A,C,\rho),&\text{otherwise.}
\end{cases}$$
Since $i\neq n-m$ we have $F_i(A,C,\rho)=(\ideala\cap S_{(n-i)})^\sat S$ and $F_{i+1}(A,C,\rho)=(\ideala\cap S_{(n-i-1)})^\sat S$. It follows that
\begin{align*}
F_{i+1}(A,C,\rho)&=(\ideala\cap S_{(n-i-1)})^\sat S=((\ideala\cap S_{(n-i)})\cap S_{(n-i-1)})^\sat S\\
&=((\ideala\cap S_{(n-i)})^\sat\cap S_{(n-i-1)})^\sat S=(F_i(A,C,\rho)\cap S_{(n-i-1)})^\sat S.
\end{align*}

\sk
c) Set $\idealb:=\bigl((\init_\degrevlex F(A,C,\rho))\cap S_{(m-1)}\bigr)^\sat S_{(m)}.$ We first prove the following claim 
$$X^{(A_m)^*}\cup X^{(C_m)^*}\cup X^{(C_m)^*+(\rho_1,\dots,\rho_m)}\subset\idealb.$$
Let $b\in(A_m)^*\cup(C_m)^*\cup ((C_m)^*+(\rho_1,\dots,\rho_m))$. 

\emph{Case 1:} $b\in(A_m)^*\cup(C_m)^*$. Then there exists $a\in A\cup C$ such that 
$b=(a_1,\dots,a_{m-1},0)$ and $m(a)\leq m$. Let $a':=a+a_m(e_{m-1}-e_m)$. Since $A\cup C$ is a Borel set, we have $a'\in A\cup C$. Hence $(a_1,\dots,a_{m-2},a_{m-1}+\nolinebreak a_m)\in (A\cup C)_{m-1}$ and $(a_1,\dots,a_{m-2},0)\in ((A\cup C)_{m-1})^*$. By Remark~\ref{Groebner bases of Borel ideals} and Proposition \ref{Groebner bases of (saturated) binomial ideals}~a) it follows that 
$$X^{(a_1,\dots,a_{m-2},0)}\in\langle X^{(A\cup C)_{m-1}}\rangle_{S_{(m-1)}}^\sat\subset\bigl((\init_\degrevlex F(A,C,\rho))\cap S_{(m-1)}\bigr)^\sat,$$ and therefore, $X^b=X_{m-1}^{a_{m-1}}X^{(a_1,\dots,a_{m-2},0,0)}\in\idealb.$

\emph{Case 2:} $b\in(C_m)^*+(\rho_1,\dots,\rho_m)$. Then there exists $c\in C$ such that $b=(c_1,\dots,c_{m-1},0)+(\rho_1,\dots,\rho_m)$. By similar arguments to those used in case 1 we obtain $a':=c+\rho+(c+\rho)_m(e_{m-1}-e_m)\in A\cup(C+\rho)$. Since $m(a')=m-1$ it holds $a'\in A$, and therefore 
$$X^{(c_1,\dots,c_{m-2},0)+(\rho_1,\dots,\rho_{m-2},0)}\in\bigl((\init_\degrevlex F(A,C,\rho))\cap S_{(m-1)}\bigr)^\sat.$$ 
It follows again that $X^b=X_{m-1}^{(c+\rho)_{m-1}}X_m^{\rho_m}X^{(c_1,\dots,c_{m-2},0,0)+(\rho_1,\dots,\rho_{m-2},0,0)} \in\idealb,$ and our claim is proved.

Now we get our statement by means of Lemma \ref{Groebner basis of F_i(A,C,rho)}:
\begin{align*}
(F(A,C,\rho)\cap S_{(m)})^\sat&=\langle X^{(A_m)^*}\cup\Bin((C_m)^*,(\rho_1,\dots,\rho_m))\rangle_{S_{(m)}}\\*
&\subset\langle X^{(A_m)^*}\cup X^{(C_m)^*}\cup X^{(C_m)^*+(\rho_1,\dots,\rho_m)}\rangle_{S_{(m)}}\subset\idealb.
\end{align*}

\sk
d) Let $0\leq i<n-1$. If $i\neq n-m$ we have by parts a) and b) that 
\begin{align*}
F_i(A,C,\rho)&=(F_i(A,C,\rho)\cap S_{(n-i-1)})S\\*
&\subset(F_i(A,C,\rho)\cap S_{(n-i-1)})^\sat S=F_{i+1}(A,C,\rho).
\end{align*}

If $i=n-m$ we have by part c) that
$$F_{n-m}(A,C,\rho)=(F(A,C,\rho)\cap S_{(m)})^\sat S
\subset\bigl((\init_\degrevlex F(A,C,\rho))\cap S_{(m-1)}\bigr)^\sat S.$$
By definition, the last ideal equals $F_{n-m+1}(A,C,\rho).$

\sk
e) By Proposition \ref{Groebner bases of (saturated) binomial ideals}~a) and Remark \ref{Groebner bases of Borel ideals} it holds 
\begin{align*}
\bigl((\init_\degrevlex F(A,C,\rho))\cap S_{(m-1)}\bigr)^\sat&=(\langle X^{A\cup C}\rangle_S\cap S_{(m-1)})^\sat\\*
&=\langle X^{((A\cup C)_{m-1})^*}\rangle_{S_{(m-1)}}.
\end{align*}
Without loss of generality we may assume that $A_{m-1}\neq\emptyset$ and set 
$$r:=\max{\{a_{m-1}\in\N\mid a\in (A\cup C)_{m-1}\}}+\rho_m.$$ Let $1\leq i<m$. It is enough to show that 
$$X_i^r\langle X^{((A\cup C)_{m-1})^*}\rangle_{S_{(m-1)}}S_{(m)}\subset(F(A,C,\rho)\cap S_{(m)})^\sat.$$
Let $a\in(A\cup C)_{m-1}$. Since $(A\cup C)_{m-1}$ is a Borel set, we have 
$$c:=a_{m-1}e_i+a^*=a+a_{m-1}(e_i-e_{m-1})\in(A\cup C)_{m-1}.$$ 

\emph{Case 1:} $c\in A_{m-1}$. Then 
$$X_i^{a_{m-1}}X^{a^*}\in F(A,C,\rho)\cap S_{(m-1)}\subset(F(A,C,\rho)\cap S_{(m)})^\sat.$$

\emph{Case 2:} $c\in C_{m-1}$. Then $c':=(c_1,\dots,c_{m-1},0,\dots,0)\in C$. Since $(A\cup(C+\rho))$ is a Borel set, we have $b:=c'+\rho+\rho_m(e_i-e_{m})\in A\cup(C+\rho).$ It is clear that $m(b)<m$, whence $b\in A$. It follows that $X_i^{a_{m-1}+\rho_m}X^{a^*}=x_i^{\rho_m}(X^{c'}-X^{c'+\rho})+X_m^{\rho_m}X^c\in F(A,C,\rho)\cap S_{(m)}\subset(F(A,C,\rho)\cap S_{(m)})^\sat$.

In both cases our claim follows.

\sk
f) From Lemma \ref{Groebner basis of F_i(A,C,rho)} it follows that 
$$\init_\degrevlex(F(A,C,\rho)\cap S_{(m)})^\sat=\langle X^{(A_m\cup C_m)^*}\rangle_{S_{(m)}}.$$ 
Thus, if $X_m u\in\init_\degrevlex(F(A,C,\rho)\cap S_{(m)})^\sat$, there exists $b\in(A_m\cup C_m)^*$ such that $X^b$ divides $X_m u$. Since $X^b$ is not divisible by $X_m$, it has to divide $u$. Hence  $u\in\init_\degrevlex(F(A,C,\rho)\cap S_{(m)})^\sat$.
\end{proof}

\begin{proposition}\label{Binomial quotients are CM}
Let $(A,C,\rho)$ be an admissible binomial system and let\/ $1\leq i<n$. Then $F_i(A,C,\rho)/F_{i-1}(A,C,\rho)$ is zero or Cohen-Macaulay of dimension $i$.
\end{proposition}

\begin{proof}
Let $m:=m(\rho)$. Assume first that $i\neq n-m+1$. Set 
$$\ideala:=F_{i-1}(A,C,\rho)\cap S_{(n-i)}.$$ Then by Lemma \ref{lemma about F_i(A,C,rho)}~a) and~b) we have $F_{i-1}(A,C,\rho)=\ideala\,S$ and $F_i(A,C,\rho)=\ideala^\sat S$. It follows that
\begin{align*}
F_i(A,C,\rho)/F_{i-1}(A,C,\rho)&=\ideala^\sat S/\ideala\,S=(\ideala^\sat/\ideala)\otimes_{S_{(n-i)}}S\\*
&=H^1_{(S_{(n-i)})_+}(\ideala)\otimes_{S_{(n-i)}}S.
\end{align*}
The (finitely generated) $S_{(n-i)}$-module $H^1_{(S_{(n-i)})_+}\!(\ideala)$ is Artinian (cf.\ \cite[7.1.4]{BS}), therefore it is zero or Cohen-Macaulay and zero-dimensional. Thus, the $S$-module 
$$F_i(A,C,\rho)/F_{i-1}(A,C,\rho)$$ 
is zero or Cohen-Macaulay of dimension $i$ (cf.\ \cite[2.1.9]{BH}).

\sk
We now prove our statement for $i=n-m+1$. Set
$$\ideala:=\bigl((\init_\degrevlex F(A,C,\rho))\cap S_{(m-1)}\bigr)^\sat S_{(m)}$$
and
$$\idealb:=(F(A,C,\rho)\cap S_{(m)})^\sat.$$
By Lemma \ref{lemma about F_i(A,C,rho)}~c) it holds $\idealb\subset\ideala$.
Since 
\begin{align*}
F_{n-m+1}&(A,C,\rho)/F_{n-m}(A,C,\rho)\\*
&=\bigl((\init_\degrevlex F(A,C,\rho))\cap S_{(m-1)}\bigr)^\sat S/(F(A,C,\rho)\cap S_{(m)})^\sat S\\*
&=(\ideala/\idealb)\otimes_{S_{(m)}}S,
\end{align*}
it is enough to show that the $S_{(m)}$-module $\ideala/\idealb$ is zero or one-dimensional and Cohen-Macaulay.

By Lemma \ref{lemma about F_i(A,C,rho)}~e) it holds
\begin{align*}
\dim \ideala/\idealb&=\dim S_{(m)}/(0\underset{S_{(m)}}{:}\ideala/\idealb)=\dim S_{(m)}/(\idealb\underset{S_{(m)}}{:}\ideala)=\dim S_{(m)}/\sqrt{(\idealb\smash{\underset{S_{(m)}}{:}}\ideala)}\\
&\leq \dim S_{(m)}/\langle X_1,\dots,X_{m-1}\rangle_{S_{(m)}}=1.
\end{align*}
Thus, it is enough to show that $X_m$ is a non-zerodivisor of $\ideala/\idealb$. 

Lemma \ref{lemma about F_i(A,C,rho)}~f) states that $X_m$ is a non-zerodivisor of $S_{(m)}/\init_\degrevlex\idealb$. Since any non-zerodivisor of $S_{(m)}/\init_\degrevlex\idealb$ is a non-zerodivisor of $S_{(m)}/\idealb$ by \cite[15.15]{E} and hence of $\ideala/\idealb$, our proof is finished.
\end{proof}

\begin{theorem}\label{binomial ideals and sCM}
Let $(A,C,\rho)$ be an admissible binomial system. Then the $S$-modules $S/F(A,C,\rho)^\sat$ and $S/F(A,C,\rho)$ are sequentially Cohen-Macaulay.
\end{theorem}

\begin{proof}
Assume that $F(A,C,\rho)\neq 0$. For $0\leq i<n$ set $$M_i:=F_i(A,C,\rho)/F(A,C,\rho)^\sat.$$ Then, by Lemma \ref{lemma about F_i(A,C,rho)}~d), one has a filtration 
$$0=M_0\subset M_1\subset\dots\subset M_{n-1}=S/F(A,C,\rho)^\sat.$$
By Proposition \ref{Binomial quotients are CM} the quotients $M_i/M_{i-1}=F_i(A,C,\rho)/F_{i-1}(A,C,\rho)$ are zero or Cohen-Macaulay of dimension $i$ for all $1\leq i<n$. So $S/F(A,C,\rho)^\sat$ is sequentially Cohen-Macaulay.

Furthermore, set $\bar M_i:=F_{i-1}(A,C,\rho)/F(A,C,\rho)$ for $1\leq i\leq n$ and $\bar M_0:=0$. Then $0=\bar M_0\subset \bar M_1\subset\dots\subset\bar M_n=S/F(A,C,\rho)$ is a filtration of $S/F(A,C,\rho)$ such that the quotients $\bar M_{i+1}/\bar M_i=F_{i}(A,C,\rho)/F_{i-1}(A,C,\rho)$ are zero or Cohen-Macaulay of dimension $i$ for all $1\leq i<n$. Using the fact that 
$$\bar M_1=F(A,C,\rho)^\sat/F(A,C,\rho)=H^1_{S_+}(F(A,C,\rho))$$ 
is Artinian, the quotient $\bar M_1/\bar M_0=F(A,C,\rho)^\sat/F(A,C,\rho)$ is zero or Cohen-Mac\-aulay and zero-dimensional. Hence $S/F(A,C,\rho)$ is also sequentially Cohen-Macaulay.
\end{proof}
\sk
\subsection{Sequences of Gr\"obner deformations}\label{subsection sequences}
It is well known that for any polynomial $p\in\Q[t]$ there exists an integer $b\in\N$ such that every saturated homogeneous ideal of $S$ with Hilbert polynomial $p$ is generated in degree $b$ (s.~Proposition \ref{admissible Hilbert polynomials}~c)). Consider the set $\Ideals$ of all homogeneous ideals of $S$ with fixed Hilbert polynomial $p$ and with generators all sitting in degree $b$. Let $f:\Z\ra\N$ be a numerical function. Define $\Ideals_{=f}:=\{\ideala\in\Ideals\mid h_{\ideala^\sat}=f\}$ and $\Ideals_{\geq f}:=\{\ideala\in\Ideals\mid h_{\ideala^\sat}\geq f\}$. Then \cite{M2000} states that any two points in $\Ideals_{=f}$ resp.\ in $\Ideals_{\geq f}$ are connected by a sequence of Gr\"obner deformations. 

Now let $\fbar:=(f_i)_{i=1}^\infty$ be a sequence of numerical functions $f_i:\Z\ra\N$ and define 
$$\Ideals^{\geq\fbar}_{=f}:=\{\ideala\in\Ideals_{=f}\mid h^i_{S/\ideala}\geq f_i\text{ for all }i\geq 1\},$$ 
$$\Ideals^{\geq\fbar}_{\geq f}:=\{\ideala\in\Ideals_{\geq f}\mid h^i_{S/\ideala}\geq f_i\text{ for all }i\geq 1\}.$$ 
Then, as we shall show below, any two points in $\Ideals^{\geq\fbar}_{=f}$ resp.\ in $\Ideals^{\geq\fbar}_{\geq f}$ are also connected by a sequence of Gr\"obner deformations (Corollary \ref{corollary of connection sequences (= f)}, Theorem \ref{theorem of connection sequences (<= f)}).

\begin{definition}
We say that two homogeneous ideals $\ideala$, $\idealb\subset S$ are \emph{connected by a Gr\"obner deformation} if one of the following equalities hold:
$$
\Gin_\degrevlex\ideala=\idealb,\quad
\init_\deglex\ideala=\idealb,\quad
\ideala=\Gin_\degrevlex\idealb\text{\quad or \quad}
\ideala=\init_\deglex\idealb.$$

Let $M$ be a set of homogeneous ideals of $S$, and let $\ideala$, $\idealb\in M$. A \emph{connecting sequence} in $M$ from $\ideala$ to $\idealb$ is a sequence $(\idealc_i)_{i=0}^r$ of ideals $\idealc_i\in M$ with the properties:
\begin{enumerate}
\item $\idealc_0=\ideala$, $\idealc_r=\idealb$;
\item $\idealc_{i-1}$ and $\idealc_i$ are connected by a Gr\"obner deformation for all $1\leq i\leq r$.
\end{enumerate}
\end{definition}

\begin{remark}\label{remark on connecting sequences}
Let $M$ be a set of homogeneous ideals of $S$, and let $(\idealc_i)_{i=0}^r$ be a connecting sequence in $M$ from $\ideala\in M$ to $\idealb\in M$. Then all ideals $\idealc_0$, \dots, $\idealc_r$ have the same Hilbert function (cf.\ Lemma \ref{ideal and initial ideal have same Hilbert function}).

Let $\ideala$, $\idealb$, $\idealc\in M\!$, and let $(\idealc_i)_{i=0}^r$ and $(\idealc_{i+r})_{i=0}^{r'}$ be two connecting sequences in $M$ from $\ideala$ to $\idealb$, respectively from $\idealb$ to $\idealc$. Then $(\idealc_i)_{i=0}^{r+r'}$ is a connecting sequence in $M$ from $\ideala$ to $\idealc$.
\end{remark}

\begin{proposition}\label{main Theorem of Herzog-Sbarra}
Let $\ideala\subset S$ be a homogeneous ideal. 

a) \cite[2.4]{S} Then $h^i_{S/\ideala}\leq h^i_{S/\init_\tau\ideala}$ for all $i\in\N$ and for any term order $\tau$ of\/ $\terms$.

b) \cite[2.2]{HS} If $\ideala\subset S$ is a Borel ideal, then $S/\ideala$ is sequentially Cohen-Macaulay.

c) \cite[3.1]{HS} The following are equivalent:
\begin{enumerate}\setlength{\labelsep}{-0.7cm}
\item \hspace{1cm}$S/\ideala$ is sequentially Cohen-Macaulay;
\item \hspace{1cm}$h^i_{S/\ideala}=h^i_{S/\Gin_\degrevlex\ideala}$ for all $i\in\N$.\vspace{-0.7cm}
\end{enumerate}\hfill$\square$
\end{proposition}

\begin{theorem}\label{theorem of connection sequences (= f)}
Let $(f_i)_{i\in\N}$ be a sequence of functions $f_i:\Z\ra\N$. Let $M$ be the set of all saturated homogeneous ideals $\ideala$ of $S$ such that $S/\ideala$ has Hilbert function $f_0$ and such that $h^i_{S/\ideala}\geq f_i$ for all $i\geq 1$. If $B\subset\N^n_d$ is a Borel set such that $\langle X^B\rangle_S^\sat\in M$ and if $L_{gh}(B)$ is the growth-height-lexico\-graphic normal form of $B$, then $\langle X^{L_{gh}(B)}\rangle_S^\sat$ is in $M$, too. Moreover, for all $\ideala$, $\idealb\in M$ there is a connecting sequence in $M$ from $\ideala$ to $\idealb$.
\end{theorem}

\begin{proof}
Let $\ideala$, $\idealb\in M$. Then $\idealc:=\Gin_\degrevlex\ideala$ is a saturated Borel ideal such that $S/\idealc$ has Hilbert function $f_0$ and such that $h^i_{S/\idealc}\geq h^i_{S/\ideala}\geq f_i$ for all $i\geq 1$ (cf.\ Proposition \ref{existence of gin}, Proposition \ref{Gin and sat commute}, Remark \ref{remark on connecting sequences} and Proposition \ref{main Theorem of Herzog-Sbarra}).

Assume that $\idealc$ and $\Gin_\degrevlex\idealb$ are both generated in degree $d$. Set $B:=\log\left(\idealc_d\cap\terms\right)$. This is a Borel set in $\N^n_d$ and it holds $\idealc=\langle X^B\rangle_S^\sat$. Let $L_{gh}(B)$ be the growth-height-lexico\-graphic normal form of $B$. 

We claim that there is a connecting sequence in $M$ from $\ideala$ to $\langle X^{L_{gh}(B)}\rangle_S^\sat$. Set $\idealc_0:=\ideala$ and $\idealc_1:=\idealc$. If $B$ is growth-height-lexicographic, we are done. Therefore, assume that $B\neq L_{gh}(B)$. By Proposition \ref{rebuilding systems are Mall} there exists a sequence of Mall binomial systems $((A_i,C_i,\rho_i))_{i=1}^r$ with the following properties:
\begin{enumerate}
\item $B=A_1\cup C_1$,
\item $A_i\cup (C_i+\rho_i)=A_{i+1}\cup C_{i+1}$ for all $1\leq i<r$,
\item $L_{gh}(B)=A_r\cup(C_r+\rho_r)$,
\item $m(\rho_i)< n$ for all $1\leq i\leq r$.
\end{enumerate}
It suffices to show that there is a connecting sequence in $M$ from $\langle X^{A_1\cup C_1}\rangle_S^\sat$ to $\langle X^{A_1\cup (C_1+\rho_1)}\rangle_S^\sat$, then one can use an inductive argument to get the claim. Set $\idealc_2:=F(A_1,C_1,\rho_1)^\sat$ and $\idealc_3:=\langle X^{A_1\cup (C_1+\rho_1)}\rangle_S^\sat$. By property (iv) we have $(\rho_1)_n=0$. Since $(A_1,C_1,\rho_1)$ is Mall, by Corollary \ref{gin and sat commute for binomial ideals} we have $\Gin_\degrevlex\idealc_2=\idealc_1$ and $\init_\deglex\idealc_2=\idealc_3$. It remains to show that $\idealc_2$ and $\idealc_3$ are in $M$. 

It is clear that $S/\idealc_2$ and $S/\idealc_3$ have the Hilbert function $f_0$ (cf.\ Remark \ref{remark on connecting sequences}). Since $(A_1,C_1,\rho_1)$ is Mall, by Theorem \ref{binomial ideals and sCM} the $S$-module $S/\idealc_2$ is sequentially Cohen-Macaulay. Thus, $S/\idealc_1$ and $S/\idealc_2$ have the same cohomological Hilbert functions (Proposition~\ref{main Theorem of Herzog-Sbarra}~c)). Finally we have $h^i_{S/\idealc_3}=h^i_{S/\init_\deglex\idealc_2}\geq h^i_{S/\idealc_2}$ for all $i\in\N$ by Proposition \ref{main Theorem of Herzog-Sbarra}~a). Thus, our claim is proved. 

Now, by Proposition \ref{computation of the Hilbert function}, the Borel set $B':=\log{((\Gin_\degrevlex\idealb)_d\cap\terms)}$ has the same growth and height vectors as $B$. Hence there exists a connecting sequence in $M$ from $\idealb$ to $\langle X^{L_{gh}(B)}\rangle_S^\sat=\langle X^{L_{gh}(B')}\rangle_S^\sat$, and our Theorem is proved. 
\end{proof}

\begin{corollary}\label{corollary of connection sequences (= f)}
Let $(f_i)_{i\in\N}$ be a sequence of functions $f_i:\Z\ra\N$. Let $M$ be the set of all homogeneous ideals $\ideala$ of $S$ with the following properties
\begin{enumerate}
\item $\ideala=(\ideala^\sat)_{\geq d}$,
\item $S/\ideala^\sat$ has Hilbert function $f_0$,
\item $h^i_{S/\ideala}\geq f_i$ for all $i\geq 1$.
\end{enumerate}
Then for all $\ideala$, $\idealb\in M$ there exists a connecting sequence in $M$ from $\ideala$ to $\idealb$.
\end{corollary}

\begin{proof}
Let $\idealc$, $\idealc'\in M$. By Theorem \ref{theorem of connection sequences (= f)} one has a connecting sequence $(\idealc_i)_{i=0}^r$ from $\idealc^\sat$ to $(\idealc')^\sat$ in the set of all saturated homogeneous ideals $\ideala$ of $S$ such that $S/\ideala$ has Hilbert function $f_0$ and such that $h^i_{S/\ideala}\geq f_i$ for all $i\geq 1$. We claim that $((\idealc_i)_{\geq d})_{i=0}^r$ is a connecting sequence from $\idealc$ to $\idealc'$ in $M$. It is clear that $(\idealc_i)_{\geq d}\in M$ for all $0\leq i\leq r$. Furthermore, by Lemma \ref{gin commutes with truncation} it holds $\Gin_\degrevlex(\ideala_{\geq d})=(\Gin_\degrevlex\ideala)_{\geq d}$ and $\init_\deglex(\ideala_{\geq d})=(\init_\deglex\ideala)_{\geq d}$ for all homogeneous ideals $\ideala\subset S$. Hence, our claim is proved.
\end{proof}

\begin{theorem}\label{theorem of connection sequences (<= f)}
Let $p\in\Q[t]$. Assume that all saturated homogeneous ideals of $S$ with Hilbert polynomial $p$ are generated in degree $d$. Let $(f_i)_{i\in\N}$ be a sequence of functions $f_i:\Z\ra\N$. Let $M$ be the set of all homogeneous ideals $\ideala$ of $S$ with the following properties
\begin{enumerate}
\item $\ideala=(\ideala^\sat)_{\geq d}$,
\item $\ideala$ has Hilbert polynomial $p$,
\item $h_{S/\ideala^\sat}\leq f_0$ ,
\item $h^i_{S/\ideala}\geq f_i$ for all $i\geq 1$.
\end{enumerate}
Then for all $\ideala$, $\idealb\in M$ there exists a connecting sequence in $M$ from $\ideala$ to $\idealb$.
\end{theorem}

\begin{proof}
Let $\ideala\in M$. Then, as in the proof of Theorem~\ref{theorem of connection sequences (= f)}, we see that $\idealc:=\Gin_\degrevlex\ideala^\sat$ is a saturated Borel ideal such that $h_{S/\idealc}=h_{S/\ideala^\sat}$ and such that $h^i_{S/\idealc}\geq f_i$ for all $i\geq 1$. Set $B:=\log\left(\idealc_d\cap\terms\right)$ and let $L_{gh}(B)$ be its growth-height-lexico\-graphic normal form. By Theorem \ref{theorem of connection sequences (= f)} the ideal $\idealc':=\langle X^{L_{gh}(B)}\rangle_S^\sat$ has the same Hilbert function as $\idealc$ and moreover, $h^i_{S/\idealc'}\geq f_i$ for all $i\geq 1$. Hence, by Corollary \ref{corollary of connection sequences (= f)}, there is a connecting sequence in $M$ from $\ideala$ to $\langle X^{L_{gh}(B)}\rangle_S$. 

\sk
Let $\ideall$ be the unique saturated lexicographic ideal with Hilbert polynomial $p$ (cf.~Lemma~\ref{existence of lex ideals}). We claim that there exists a connecting sequence in $M$ from $\langle X^{L_{gh}(B)}\rangle_S$ to $\ideall_{\geq d}$. Set $L:=\log\left(\ideall_d\cap\terms\right)$. By our assumptions, $\ideall$ is generated in degree $d$. So we have $\ideall_{\geq d}=\langle X^L\rangle_S$, and Corollary \ref{Hilb poly determines growth-vector} tells us that $\gv(L_{gh}(B))=\gv(L)$. Assume that $L_{gh}(B)\neq L$. Then, by Proposition \ref{existence of Mall bin syst to a lexicographic segment}, there exists a finite sequence of Mall binomial systems $((A_i,C_i,\rho_i))_{i=1}^r$ with the following properties:
\begin{enumerate}
\item $L_{gh}(B)=A_1\cup C_1$,
\item $A_i\cup (C_i+\rho_i)=A_{i+1}\cup C_{i+1}$ for all $1\leq i<r$,
\item $L=A_r\cup(C_r+\rho_r)$,
\item $\gv(A_i\cup C_i)=\gv(L)$ for all $1\leq i\leq r$,
\item $h_{\langle X^{A_i\cup C_i}\rangle_S^\sat}<h_{\langle X^{A_i\cup (C_i+\rho_i)}\rangle_S^\sat}$ for all $1\leq i\leq r$.
\end{enumerate}
It suffices to show that there is a connecting sequence in $M$ from $\langle X^{A_1\cup C_1}\rangle_S$ to $\langle X^{A_1\cup (C_1+\rho_1)}\rangle_S$, then one can use an inductive argument to get the claim. Set $\idealc_0:=\langle X^{A_1\cup C_1}\rangle_S$, $\idealc_1:=F(A_1,C_1,\rho_1)$ and $\idealc_2:=\langle X^{A_1\cup (C_1+\rho_1)}\rangle_S$. Since $(A_1,C_1,\rho_1)$ is Mall, by Proposition \ref{Groebner bases of (saturated) binomial ideals} and Theorem \ref{Gin = in for good binomial ideals} we have $\Gin_\degrevlex\idealc_1=\idealc_0$ and $\init_\deglex\idealc_1=\idealc_2$. It remains to show that $\idealc_1$ and $\idealc_2$ are in $M$. By Proposition \ref{Groebner bases of (saturated) binomial ideals} f) and Remark \ref{Groebner bases of Borel ideals} we have $(\idealc_1^\sat)_{\geq d}=\idealc_1$ and $(\idealc_2^\sat)_{\geq d}=\idealc_2$. As in the proof of Theorem \ref{theorem of connection sequences (= f)} one sees that $h_{S/\idealc_0^\sat}=h_{S/\idealc_1^\sat}$ and $h^i_{S/\idealc_2}\geq h^i_{S/\idealc_1}=h^i_{S/\idealc_0}$ for all $i\in\N$. From the above properties (ii) and (iv) it follows that $\gv(A_1\cup (C_1+\rho_1))=\gv(L)$, hence by Proposition \ref{computation of the Hilbert function} a), the ideal $\idealc_2$ has Hilbert polynomial $p$. From property (v) it follows that $h_{S/\idealc_2^\sat}<h_{S/\idealc_0^\sat}\leq f_0$. Thus, $\idealc_1$, $\idealc_2\in M$, and our claim is proved.

\sk
Altogether we have constructed a connecting sequence in $M$ from $\ideala$ to $\ideall_{\geq d}$. Connecting another ideal $\idealb\in M$ with $\ideall_{\geq d}$ in $M$ yields a connecting sequence in $M$ from $\ideala$ to $\idealb$.
\end{proof}
\clearpage
\section{Hilbert function strata}
Let $p\in\Q[t]$ be a polynomial. The Hilbert scheme $\Hilb^p_K$ is defined as the representing scheme of the Hilbert functor $\Hilbf^p_K:\sch_K\ra\sets$ which assigns to each locally Noetherian scheme $T$ over $K$ the set of all coherent quotient sheaves $\sheaff$ of $\sheafo_{\proj[T]}$ which are flat over $T$ with Hilbert polynomial $p_{\sheaff,x}=p$ for all $x\in T$ (s.~section \ref{subsection Hilbert scheme}). 

To each point $x\in\Hilb^p_K$ is associated a quotient sheaf $$\smash{\sheaff^{(x)}\in\Hilbf^p_K(\Spec(\kappa(x))).}$$ Hence it is possible to define subsets of $\Hilb^p_K$ by bounding the vector space dimensions of the Serre cohomology groups $H^i(\proj[\kappa(x)],\sheaff^{(x)}(j))$ of $\sheaff^{(x)}$ for $x\in\Hilb^p_K$, $i\in\N$ and $j\in\Z$ (s.~section~\ref{subsection main theorem}). 

Let $X_K\subset\Hilb^p_K$ be such a topological subspace defined by bounding cohomology. How to prove that it is connected? We proceed as follows: By use of the Serre-Grothendieck Correspondence, the bounding conditions are translated into the language of local cohomology. Hence, if $\Char K=0$, we may try to apply the previous results about connecting sequences of Gr\"obner deformations. To prove that connectedness by Gr\"obner deformations yields topological connectedness (cf.~section~\ref{subsection connectedness}), it is crucial that $X_K$ is closed under isomorphisms: If $x$, $y\in\Hilb^p_K$ are two closed points and $x\in X_K$, then $\sheaff^{(x)}\cong\sheaff^{(y)}$ implies $y\in X$. 

To get the connectedness of $X_K$ in case the field $K$ is not algebraically closed, we prove that $X_K$ carries a reduced scheme structure such that $(X_K\times_{\Spec(K)}\Spec(k))_\red=X_k$ for any field extension $K\subset k$. This is done in two steps: First by showing by means of the semicontinuity Theorem that $X_K$ is locally closed in $\Hilb^p_K$. Then by constructing a subfunctor $\Xf_K:\sch_K\ra\sets$ of the Hilbert functor $\Hilbf^p_K$, defined by the bounding conditions imposed on $X_K$, and by showing that $\Xf_K$ is representable by $X_K$ (cf.~section~\ref{subsection main theorem}). 

As before $n$ denotes a positive integer. But now we redefine the polynomial ring $S$ by $S:=K[X_0,\dots,X_n]$.

\subsection{Cohomological Hilbert functions and  Hil\-bert polynomials}
We first describe basic properties of graded local cohomology and cohomology of sheaves on a projective scheme. As an introduction to local cohomology we refer to \cite{BS, B}. As an introduction to sheaves and schemes we recommend \cite{H, EGA I}.

\begin{notation, definition and remark}\label{CM-regularity}
Let $R:=\bigoplus_{i\in\N}R_i$ be a positively graded homogeneous Noetherian ring with $R_0=K$. 

A) Let $M$ be a graded $R$-module. For $i\in\N$ we consider the $i$th local cohomology module $H^i_{R_+}(M)$ of $M$ with respect to the irrelevant ideal $R_+$ endowed with its natural grading. If $M$ is finitely generated, the $K$-vector space $H^i_{R_+}(M)_k$ has finite dimension for any $i\in\N$ and $k\in\Z$ (\cite[15.1.5(i)]{BS}). For $i\in\N$ we denote by 
$$h^i_M:\Z\ra\N,\;k\mapsto\dim_K H^i_{R_+}(M)_k$$ 
the $i$th \emph{cohomological Hilbert function} of $M$.

The \emph{end} of $M$ is defined by $\End(M):=\sup\,\{i\in\Z\mid M_i\neq 0\}$. (Use the following convention: $\sup A=\infty$ if $A\subset\Z$ is not bounded above; $\sup\emptyset=-\infty$; $-\infty+i=-\infty$ for all $i\in\Z$. For infima use analogous conventions).

The \emph{(Castelnuovo-Mumford) regularity} of $M$ is defined by 
$$\reg M:=\sup{\{\End{(H^i_{R_+}(M))}+i\mid i\in\N\}}.$$ 
If $M$ is finitely generated, it holds $\reg M<\infty$ (cf.\ \cite[15.1.5(ii), 6.1.2]{BS}).

If $M$ is finitely generated and non zero, then $M$ is generated in degree $\reg M$ (\cite[15.3.1]{BS}).

\sk
B) Let $X:=\Proj(R)$ be the projective scheme induced by $R$. If $M$ is a $R$-module, we denote the induced quasi-coherent sheaf of $\sheafo_X$-modules by $\widetilde M$. 

For $i\in\Z$ let $R(i)$ be the shifted graded $R$-module with $R(i)_k=R_{i+k}$ for all $k\in\Z$. Define $\sheafo_X(i):=\widetilde{R(i)}$. Let $\sheaff$ be a sheaf of $\sheafo_X$-modules. For $i\in\Z$ let $\sheaff(i):=\sheaff\otimes_{\sheafo_X}\sheafo_X(i)$ denote the $i$th \emph{twist} of $\sheaff$. 

For $i\in\Z$ let $H^i(X,\sheaff)$ denote the $i$th \emph{Serre cohomology group} of $X$ with coefficients in $\sheaff\!$, endowed with its natural $K$-vector space structure. 

If $\sheaff$ is coherent, then the $K$-vector space $H^i(X,\sheaff(k))$ has finite dimension for any  $i\in\N$ and $k\in\Z$ (cf.~\cite[III.\ 5.2(a)]{H}). For $i\in\N$ we denote by 
$$h^i_{X,\sheaff}:\Z\ra\N,\;k\mapsto\dim_K H^i(X,\sheaff(k))$$
the $i$th \emph{cohomological Hilbert function} of $(X,\sheaff)$.

Let $m\in\Z$. The sheaf $\sheaff$ is said to be \emph{$m$-regular} if $H^i(X,\sheaff(k-i))=0$ for all $i>0$ and $k\geq m$.

The \emph{(Castelnuovo-Mumford) regularity} of $\sheaff$ is defined by 
$$\reg\sheaff:=\inf\,\{m\in\Z\mid \sheaff \text{ is $m$-regular}\}.$$ 
If $\sheaff$ is coherent, it holds $\reg \sheaff<\infty$ (cf.\ \cite[III.\ 2.7, 5.2(b)]{H}).

If $A$ is a commutative ring, set $\proj[A]:=\Proj(A[X_0,\dots,X_n])$. If $T$ is a scheme over $K$, let $\proj[T]$ denote the fiber product $\proj\times_{\Spec(K)}T$. 

If $\ideala\subset S$ is a homogeneous ideal, then 
$$\widetilde\ideala=\widetilde{\ideala^\sat}\text{\quad and\quad}\bigoplus_{k\in\N}H^0(\proj,\widetilde\ideala(k))\cong\ideala^\sat$$ 
(cf.\ \cite[II.\ Ex.\ 5.10]{H}). 
\end{notation, definition and remark}

\begin{proposition}[{\cite[12.4.2, 20.4.4]{BS}, Serre-Grothendieck Correspondence}]\label{Serre-Grothendieck}
Let $R:=\bigoplus_{i\in\N}R_i$ be a homogeneous Noetherian ring and $M$ a finitely generated $R$-module. Then there is an exact sequence of $R$-modules
\begin{equation*}\label{exact sequence}
0\ra H^0_{R_+}(M)\ra M\ra \bigoplus_{k\in\N}H^0(X,\widetilde M(k))\ra H^1_{R_+}(M)\ra 0
\end{equation*}
where all homomorphisms are homogeneous, and there are homogeneous isomorphisms of $R$-modules
\begin{equation*}\label{isomorphisms}
\bigoplus_{k\in\N}H^i(X,\widetilde M(k))\overset{\cong}{\ra}H^{i+1}_{R_+}(M)
\end{equation*}
for all $i\geq1$. \hfill$\square$
\end{proposition}

\begin{proposition}[{\cite[13.1.8]{BS}, Graded flat base change}]
Let $R:=\bigoplus_{i\in\N}R_i$ be a homogeneous Noetherian ring, $R_0'$ a flat Noetherian $R_0$-algebra, $R':=R\otimes_{R_0}R_0'$ and $M$ a graded $R$-module. Then there are  homogeneous isomorphisms  of $R'$-modules
$$H^i_{R_+}(M)\otimes_R R'\overset{\cong}{\ra}H^i_{R_+'}(M\otimes_R R')$$
for all $i\in\N$. \hfill$\square$
\end{proposition}

\begin{remark}\label{consequences of Serre-Grothendieck Correspondence}
We collect some consequences of the Serre-Grothendieck Correspondence and graded flat base change. Let $R:=\bigoplus_{i\in\N}R_i$ be a homogeneous Noetherian ring with $R_0=K$. 

A) Let $M$ be a finitely generated graded $R$-module. Then 
$$\reg M=\sup{\{\End{(H^0_{R_+}(M))},\:\End{(H^1_{R_+}(M))}+1,\:\reg\sheaff\}}.$$

B) Let $\ideala\subset S$ be a homogeneous ideal. Then $H^1_{S_+}(\ideala)=\ideala^\sat/\ideala$ and therefore, $H^0_{S_+}(S/\ideala^\sat)=0$. Hence for each $k\in\Z$ there is an exact sequence of $K$-vector spaces
$$0\ra(S/\ideala^\sat)_k\ra H^0(\proj,\widetilde{S/\ideala}(k))\ra H^1_{S_+}(S/\ideala^\sat)_k\ra 0.$$
It holds $H^i_{S_+}(S)=0$ for $i\neq n+1$; and $H^{n+1}_{S_+}(S/\ideala)\neq 0$ if and only if $\ideala\neq 0$ (cf.\ \cite[6.1.2, 6.2.7]{BS}). The assumption $n>0$ implies $H^1_{S_+}(S)=0$, hence there is an exact sequence
$$0\ra H^1_{S_+}(S/\ideala)\ra H^2_{S_+}(\ideala)\ra H^2_{S_+}(S).$$
It follows that $h^1_{S/\ideala}=h^1_{S/\ideala^\sat}=h^0_{\proj,\widetilde{S/\ideala}}-h_{S/\ideala^\sat}$. 

C) Let $\sheaff$ a coherent sheaf of $\sheafo_{\Proj(R)}$-modules. Let $R_0'$ be a field extension of $K$, let $R':=R\otimes_{R_0}R_0'$, and let $f:\Proj(R')\ra\Proj(R)$ be the canonical morphism. Then  $h^i_{\Proj(R),\sheaff}=h^i_{\Proj(R'),f^*\sheaff}$ for all $i\in\N$.
\end{remark}

\begin{definition and remark}\label{Hilbert polynomial of sheaves}
Let $R:=\bigoplus_{i\in\N}R_i$ be a homogeneous Noetherian ring with $R_0=K$, $X:=\Proj(R)$, and let $\sheaff$ be a coherent sheaf of $\sheafo_X$-modules. Then there exists a polynomial $p\in\Q[t]$, the \emph{Hilbert polynomial} of $\sheaff$, such that $p(i)=h^0_{X,\sheaff}(i)$ for all $i\gg 0$. More precisely, since 
$$p(m)=\sum_{i\in\N}(-1)^i\dim_K H^i(X,\sheaff(m))$$ 
for all $m\in\Z$, it holds $p(m)=\dim_K H^0(X,\sheaff(m))=h^0_{X,\sheaff}(m)$ for all $m\geq\reg(\sheaff)$ (cf.\ \cite[III.\ Ex.\ 5.2]{H}).

If $R=S$ and if $\sheaff=\widetilde M$ is induced by a finitely generated graded $S$-module $M$, then the Hilbert polynomial of $\sheaff$ is the same as the Hilbert polynomial of $M$.
\end{definition and remark}

\begin{proposition}\label{admissible Hilbert polynomials} Let $p\in\Q[t]$ be an admissible Hilbert polynomial.

a) \cite[6.]{Sp} There exist $s\in\N$ with $0\leq s<n$ and $b_0$,~\dots, $b_s\in\N$ with $0\leq b_0\leq\dots\leq b_s$ such that 
$$p=\sum_{i=0}^s\binom{t+n-b_i-i}{n-i}.$$

b) \cite[2.9]{Go1978} Let $s$,~$b_0$,~\dots, $b_s\in\N$ be as in part a). Let $\sheafi\subset\sheafo_{\proj}$ be a coherent ideal sheaf with Hilbert polynomial $p$. Then $\sheafi$ is $b_s$-regular. 

c) Let $s$,~$b_0$,~\dots, $b_s\in\N$ be as in part a). Then every saturated homogeneous ideal in $S$ with Hilbert polynomial $p$ is generated in degree $b_s$.
\end{proposition}

\begin{proof}
c) Let $\ideala\subset S$ be a saturated homogeneous ideal with Hilbert polynomial $p$, and let $\sheafi:=\widetilde\ideala\subset\sheafo_{\proj}$ be the ideal sheaf induced by $\ideala$. Then $\sheafi$ is $b_s$-regular. Since $\ideala$ is saturated, $H^0_{R_+}(\ideala)=H^1_{R_+}(\ideala)=0$. Thus the regularity of $\ideala$ is $\reg\ideala=\sup\,\{\End(H^0_{R_+}(\ideala)),\End(H^1_{R_+}(\ideala))+1,\reg\sheafi\}=\reg\sheafi\leq b_s$. It follows that the ideal $\ideala$ is generated in degree $b_s$ (cf.\ \ref{CM-regularity} A)).
\end{proof}
\sk
\subsection[The Hilbert scheme $\Hilb^p_K$]{The Hilbert scheme $\boldsymbol{\Hilb^p_K}$}\label{subsection Hilbert scheme}
An introduction to Hilbert schemes in a general setting can be found in \cite{Gr}. As a more elementary introduction we recommend \cite{St}.

\begin{notation}
Let $\sets$ denote the category of sets, $\alg_K$ the category of Noetherian $K$-algebras and $\sch_K$ the category of locally Noetherian schemes over $K$.

If $T$ is a scheme, $\iota:U\hookrightarrow T$ an open subscheme and $\sheaff$ a sheaf on $T$, we write $\sheaff\restr_U:=\iota^*\sheaff$ for the restriction of $\sheaff$ on $U$.

If $T$ is a $K$-scheme, define the twisting sheaf $\sheafo_{\proj[T]}(1):=\pi^*(\sheafo_{\proj}(1))$, where $\pi:\proj[T]\ra\proj$ is the canonical morphism.

Let $R$ be a graded ring, and $M$ a graded $R$-module. If $\idealp\subset R$ is a prime ideal, set $M_{[\idealp]}:=(\bigcup_{i\in\Z}R_i\setminus\idealp)^{-1}M$. The $0$th component of this graded module is denoted by $M_{(\idealp)}:=(M_{[\idealp]})_0.$
\end{notation}

\begin{remark}\label{remark about computing inverse images}
Let $f:A\ra B$ be a homomorphism of $K$-algebras and let $\iota:\proj[B]\ra\proj[A]$ denote the induced morphism of schemes. Let $M$ be a graded $A[X_0,\dots,X_n]$-module and $m\in\Z$. Then 
$$\iota^*(\widetilde M(m))\cong\widetilde{M\otimes_A B}(m)$$
(cf.\ \cite[2.8.11]{EGA II}).
\end{remark}

\begin{lemma}\label{lemma about flatness}
a) (\cite[0.5.7.4]{EGA I}) Let $X\ra Y$ be a morphism of schemes and $0\ra\sheaff\ra\sheafg\ra\sheafh\ra 0$
an exact sequence of $\sheafo_X$-modules. If $\sheafh$ is flat over $Y$, then $\sheafg$ is flat over $Y$ if and only if $\sheaff$ is flat over $Y$.

b) Let $Y$ be a Noetherian scheme and $\sheaff$ a coherent sheaf of $\sheafo_{\proj[Y]}$-modules, flat over $Y$. Then for each $m\in\Z$ the twisted sheaf $\sheaff(m)$ is flat over $Y$.

c) Let $R:=\bigoplus_{i\in\N}R_i$ be a graded ring and $M$ a graded $R$-module. Let $X:=\Proj(R)$, $Y:=\Spec(R_0)$ and $f:X\ra Y$ the canonical morphism of schemes. If $M$ is flat over $R_0$, then $\widetilde M$ is flat over $Y$.
\end{lemma}

\begin{proof}
b) Let $f:\proj[Y]\ra Y$ be the canonical morphism. By \cite[7.9.14]{EGA III} the sheaf $\sheaff$ is flat over $Y$ if and only if the direct image sheaf $f_*(\sheaff(m))$ is a locally free $\sheafo_Y$-module for all $m\gg 0$.

c) A sheaf $\sheaff$ of $\sheafo_X$-modules is flat over $Y$ if the stalks $\sheaff_x$ are flat over $\sheafo_{Y\!,\,f(x)}$ for all $x\in X$. Hence the induced sheaf of $\sheafo_X$-modules $\widetilde M$ is flat over $Y$ if and only if $M_{(\idealp)}$ is flat over $(R_0)_{\idealp\cap R_0}$ for all $\idealp\in\Proj(R)$ (cf.~\cite[II.~5.11]{H}). Now assume that $M$ is flat over $R_0$. Let $\idealp\in\Proj(R)$ and $\idealq:=\idealp\cap R_0$. Then $M_\idealp$ is flat over $(R_0)_{\idealq}$. Observe that there are inclusions of $(R_0)_\idealq$-modules $N_{(\idealp)}\subset N_{[\idealp]}\subset N_\idealp$ for any graded $R$-module $N$. Let $V$ be a $(R_0)_\idealq$-module. Then  $M_{[\idealp]}\otimes_{(R_0)_\idealq}\!V=R_{[\idealp]}\otimes_R M\otimes_{(R_0)_\idealq}\!V=(M\otimes_{(R_0)_\idealq}\!V)_{[\idealp]}$ and $M_{(\idealp)}\otimes_{(R_0)_\idealq}\!V=(M_{[\idealp]}\otimes_{(R_0)_\idealq}\!V)_0$. We deduce that $M_{(\idealp)}\otimes_{(R_0)_\idealq}\!V\,\subset\, M_{[\idealp]}\otimes_{(R_0)_\idealq}\!V\,\subset\, M_\idealp\otimes_{(R_0)_\idealq}\!V$\!, and it follows easily that $M_{(\idealp)}$ is flat over $(R_0)_{\idealq}$.
\end{proof}

\begin{notation and remark}\label{Hilb poly of sheaves}
Let $X\ra Y$ be a morphism of schemes over $K$, let $y\in Y$, and let $\sheaff$ be a sheaf of $\sheafo_X$-modules. Let $\kappa(y)$ denote the residue field of $y$ on $Y$. Let $\kappa(y)\subset k$ be a field extension and $\iota:\Spec(k)\ra Y$ the canonical morphism of schemes. We write $\sheaff_{k}:=(\iota')^*\sheaff$ for the inverse image of $\sheaff$ by the induced morphism $\iota':X\times_Y\Spec(k)\ra X$. 

Now let $X:=\proj[Y]$. If $Y$ is Noetherian, then $X$ is also Noetherian. If $Y$ is locally Noetherian and $\sheaff$ is coherent on $X$, then $X\times_Y\Spec(k)=\proj[k]$, and $\sheaff_{k}$ is a coherent $\sheafo_{\proj[k]}$-module (cf.\ \cite[3.3.1]{EGA I}). By flat base change it holds $h^i_{\proj[k]\!,\,\sheaff_k}=h^i_{\proj[\kappa(x)],\,\sheaff_{\kappa(x)}}$ for all $i\in\N$ (cf.\ Remark \ref{consequences of Serre-Grothendieck Correspondence}~C)). In particular the Hilbert polynomials of $\sheaff_{\kappa(y)}$ and of $\sheaff_k$ are identical. The Hilbert polynomial of $\sheaff_{\kappa(y)}$ is denoted by $p_{\sheaff\!,\,y}$.
\end{notation and remark}

\begin{proposition}[{\cite[1.2]{H1966}}]\label{flatness and Hilbert polynomial}
Let $Y$ be a locally Noetherian scheme and $\sheaff$ a coherent sheaf of $\sheafo_{\proj[Y]}$-modules. If $\sheaff$ is flat over $Y$, then the function $Y\ra\Q[t]$, $y\mapsto p_{\sheaff\!,\,y}$ is locally constant on $Y$. The converse is true if $Y$ is integral.\hfill\qedsymbol
\end{proposition}

\begin{definition and remark}\label{def - Hilbert scheme - def}
If $T$ is a locally Noetherian scheme over $K$ and if $p\in\Q[t]$ is a polynomial, let $\Hilbf_K^p(T)$ denote the set of all coherent quotient sheaves $\sheaff$ of $\sheafo_{\proj[T]}$ which are flat over $T$ such that $p_{\sheaff,x}=p$ for all $x\in T$. 

Let $p\in\Q[t]$ be a polynomial. Define the polynomial $q\in\Q[t]$ by $q(t):=\binom{t+n}{n}-p(t)$. If $\ideala$ is a homogeneous ideal of $S$ with Hilbert polynomial $q$, then the quotient $\sheafo_{\proj[K]}/\widetilde\ideala$ has Hilbert polynomial $p$ and is flat over $\Spec(K)$. On the other hand the kernel of a morphism of coherent sheaves over a Noetherian scheme is coherent (cf.\ \cite[II. 5.7]{H}). Hence there is a bijection
$$\{\ideala\subset S\mid\begin{subarray}{l}\ideala\text{ is a saturated homogeneous ideal}\\ 
\text{with Hilbert polynomial }q\end{subarray}\}\cong\Hilbf^p_K(\Spec(K)).$$ 

If $f:T'\ra T$ is a morphism of locally Noetherian schemes over $K$, $f':\proj[T']\ra\proj[T]$ the induced morphism and $\sheaff$ a coherent quotient sheaf of $\sheafo_{\proj[T]}$ which is flat over $T$, then $(f')^*\sheaff$ is a coherent quotient sheaf of $\sheafo_{\proj[T']}$, flat over $T'$, and it holds $p_{(f')^*\sheaff\!,\,x}=p_{\sheaff\!,\,f(x)}$ for all $x\in T'$ (cf.\ \cite[II, 5.8; III, 9.2]{H}, \cite[3.3.1]{EGA I} and \ref{Hilb poly of sheaves}). Therefore, there is a map 
$$\Hilbf_K^p(f):\Hilbf_K^p(T)\ra\Hilbf_K^p(T'),\;\sheaff\mapsto(f')^*\sheaff.$$

Thus, we have defined a contravariant functor $\Hilbf_K^p:\sch_K\ra\sets$, the \emph{Hilbert functor}.

Grothendieck showed in \cite{Gr} that the Hilbert functor is representable by a projective scheme $\Hilb^p_K$, the \emph{Hilbert scheme}. We denote the representing natural equivalence of functors by $\psi_{\Hilb}:\Mor_{\sch_K}(\bullet,\,\Hilb^p_K)\ra\Hilbf^p_K$.

Define the coherent quotient sheaf of $\sheafo_{\proj[\Hilb_K^p]}$
$$\sheafhilb:=\psi_{\Hilb}(\Hilb_K^p)(\id_{\Hilb_K^p})\in\Hilbf^p_K(\Hilb^p_K),$$ 
the so-called \emph{universal sheaf}. Then for all morphisms $g\in\Mor_{\sch_K}(T,\,\Hilb^p_K)$ $$\psi_{\Hilb}(T)(g)=\Hilbf^p_K(g)(\sheafhilb).$$

The universal sheaf has the following property: For any locally Noetherian scheme $X$ over $K$ and for any $\sheaff\in\Hilbf^p_K(X)$ there is a unique morphism of schemes $f:X\ra\Hilb_K^p$ such that $\sheaff=\Hilbf^p_K(f)(\sheafhilb)$.
\end{definition and remark}
\sk
\subsection{Connectedness by Gr\"obner deformations}\label{subsection connectedness}
To prove the connectedness of the Hilbert scheme one uses the standard argument that two points, connected by a Gr\"obner deformation, are linearly connected. In this section we prove an analogue for certain subsets of the Hilbert scheme.

\begin{notation}
If $X$ is a scheme, let $\closed(X)$ denote the set of closed points of $X$ with the induced topology.
\end{notation}

\begin{remark and notation}\label{notations for points in Hilb}
Let $K$ be algebraically closed and $p\in\Q[t]$ a polynomial. 

A) Since $\Hilb^p_K$ is of finite type over $K$, each closed point of $\Hilb^p_K$ has residue field $K$. Hence, there is a bijection 
$$\iota:\closed(\Hilb^p_K)\ra\Mor_{\sch_K}(\Spec(K),\,\Hilb^p_K).$$

Let $\gamma:\closed(\Hilb^p_K)\ra\Hilbf^p_K(\Spec(K))$ denote the composition of $\iota$ and $\psi_{\Hilb}(\Spec(K)):\Mor_{\sch_K}(\Spec(K),\,\Hilb^p_K)\ra\Hilbf^p_K(\Spec(K))$.
If $x\in\closed(\Hilb^p_K)$, set
$$\sheaff^{(x)}:=\gamma(x).$$

B) Define the polynomial $q\in\Q[t]$ by $q(t):=\binom{t+n}{n}-p(t)$. Assume that $q$ is admissible. Let $\Ideals$ be the set of all saturated homogeneous ideals of $S$ with Hilbert polynomial $q$. By Proposition \ref{admissible Hilbert polynomials} and Remark \ref{Hilbert polynomial of sheaves} we may choose an integer $b\in\N$ such that any ideal $\ideala\in\Ideals$ has the following properties: $\ideala$ is generated in degree $b$ and $h_\ideala(i)=q(i)$ for all $i\geq b$.

Let $\Ideals_p:=\{\ideala_{\geq b}\mid\ideala\in\Ideals\}$. Then there is a bijection 
$$\delta:\Ideals_p\ra\closed(\Hilb^p_K),\;\ideala\mapsto\gamma^{-1}(\sheafo_{\proj}/\widetilde\ideala).$$ 

C) Let $\ideala\in\Ideals_p$, let $g\in\Gl(n+1,K)$ and let $\tau$ be a term order of $\terms$. Choose $\idealb\in\{g(\ideala),\init_\tau\ideala\}$. Then $h_\ideala=h_\idealb$ (cf.~Lemma \ref{ideal and initial ideal have same Hilbert function}) and therefore, $\idealb^\sat\in\Ideals$. By our choice of $b$ we have $h_{\idealb^\sat}(i)=q(i)=h_\ideala(i)=h_\idealb(i)$ for all $i\geq b$, whence $(\idealb^\sat)_{\geq b}=\idealb_{\geq b}=\idealb$. This shows that $\idealb\in\Ideals_p$.
\end{remark and notation}

\begin{definition}
Let $K$ be algebraically closed and $p\in\Q[t]$ a polynomial. 
A subset $W\subset\closed(\Hilb^p_K)$ is called \emph{closed under isomorphisms} if for each $w\in W$ and for each $x\in\closed(\Hilb^p_K)$ the following holds:
\vspace{-1mm}
\begin{center}
If $\sheaff^{(x)}$ and $\sheaff^{(w)}$ are isomorphic as $\sheafo_{\proj}$-modules, then $x\in W$.

\end{center}
\end{definition}

\begin{remark}\label{remark about natural transformations}
Let $F$, $G:\sch_K\ra\sets$ be two contravariant representable functors such that there exists a natural transformation $\eta:F'\ra G'$, where the covariant functors $F'$, $G':\alg_K\ra\sets$ are defined by $F':=F\circ\Spec$ and $G':=G\circ\Spec$. Then there exists a natural transformation $\bar\eta:F\ra G$ such that $\bar\eta(\Spec(A))=\eta(A)$ for all Noetherian $K$-algebras $A$.
\end{remark}

\begin{definition and remark}\label{remark about the general linear group}
The \emph{general linear group} over $K$ is defined by $\Gl^n_K:=\Spec(K[Y_{ij}\mid 0\leq i< n,\: 0\leq j<n]_{\det([Y_{ij}\mid 0\leq i<n,\: 0\leq j<n])}).$ The set of all invertible $n\times n$ matrices with entries in a $K$-algebra $A$ is denoted by $\Gl(n,A)$. We consider the contravariant functor $\Glf^n_K:\sch_K\ra\sets$, $T\rightsquigarrow\Gl(n,\sheafo_T(T))$ which sends each $K$-scheme to the set of all invertible $n\times n$ matrices with entries in the global sections of its structure sheaf. We identify $\Glf^n_K(\Spec(A))=\Gl(n,A)$ for any $K$-algebra $A$. The functor $\Glf^n_K$ is represented by the pair $(\Gl^n_K,[Y_{ij}\mid 0\leq i<n,\: 0\leq j<n])$. This means that there exists a natural equivalence of functors $$\psi_{\Gl}:\Mor_{\sch_K}(\bullet,\,\Gl^n_K)\ra\Glf^n_K$$ such that the map $\psi_{\Gl}(\Gl^n_K):\Mor_{\sch_K}(\Gl^n_K,\,\Gl^n_K)\ra\Glf^n_K(\Gl^n_K)$ sends the identity morphism $\id_{\Gl^n_K}$ to the matrix $[Y_{ij}\mid 0\leq i<n, 0\leq j<n]$.
\end{definition and remark}

\begin{definition and remark}
Let $K[z]$ be a polynomial ring over $K$ in one indeterminate. The \emph{affine line} over $K$ is defined by $\affline_K:=\Spec(K[z])$. We consider the contravariant functor $\afflinef_K:\sch_K\ra\sets$, $T\rightsquigarrow\sheafo_T(T)$ which sends each $K$-scheme to the set of global sections of its structure sheaf. We identify $\afflinef_K(\Spec(A))=A$ for any $K$-algebra $A$, in particular we have $\afflinef_K(\affline_K)=K[z]$. The functor $\afflinef_K$ is represented by the pair $(\affline_K,z)$. This means that there exists a natural equivalence of functors $$\psi_{\affline}:\Mor_{\sch_K}(\bullet,\,\affline_K)\ra\afflinef_K$$ such that the map $\psi_{\affline}(\affline_K):\Mor_{\sch_K}(\affline_K,\,\affline_K)\ra\afflinef_K(\affline_K)$ sends the identity morphism $\id_{\affline_K}$ to the indeterminate $z$.
\end{definition and remark}

\begin{proposition}\label{proposition about the affine line}
Let $K$ be algebraically closed, $p\in\Q[t]$ a polynomial, $W\subset\closed(\Hilb^p_K)$ a topological subspace, closed under isomorphisms, and $\ideala\in\Ideals_p$ such that $x:=\delta(\ideala)\in W$.

a) Let $g\in\Gl(n+1,K)$. Then $y:=\delta(g(\ideala))\in W$, and $x$ and $y$ lie in the same connected component of $W$.

b) Let $\tau$ be a term order of\/ $\terms$. If $y:=\delta(\init_\tau\ideala)\in W$, then $x$ and $y$ lie in the same connected component of $W$.
\end{proposition}

\begin{proof}
a) Since the general linear group $\Gl^{n+1}_K$ acts on the Hilbert scheme $\Hilb^p_K$ by linear transformation of coordinates, there is a natural transformation $\underline{\zeta}:\Glf^{n+1}_K\ra\Hilbf^p_K$ such that $\underline\zeta(\Spec(A))(h)=\sheafo_{\proj[A]}/(h(\ideala\otimes_K\!A))\sptilde$ for any $K$-algebra $A$ and all $h\in\Gl(n+1,A)$. In particular, we have $\underline\zeta(\Spec(K))(h)=\sheaff^{(\delta(h(\ideala)))}$ for all $h\in\Gl(n+\nolinebreak 1,K)$. Since $\sheaff^{(\delta(\ideala))}\cong\sheaff^{(\delta(h(\ideala)))}$ for all $h\in\Gl(n+\nolinebreak 1,K)$ and since $W$ is closed under isomorphisms, we have $\delta(h(\ideala))\in W$ for all $h\in\Gl(n+\nolinebreak 1,K)$.

By the Lemma of Yoneda (see for example \cite[A5.1]{E}) there is a unique morphism of $K$\!-schemes $\zeta:\Gl^{n+1}_K\ra\Hilb_K^p$ such that the diagram
$$\xymatrix@C+2cm{
\Mor_{\sch_K}(\bullet,\,\Gl^{n+1}_K)\ar[r]^-{\Mor_{\sch_K}(\bullet,\;\zeta)}\ar[d]_{\psi_{\Gl}}&\Mor_{\sch_K}(\bullet,\,\Hilb^p_K)\ar[d]^{\psi_{\Hilb}}\\
\Glf^{n+1}_K\ar[r]^-{\underline{\zeta}}&\Hilbf^p_K
}$$
commutes. 

Let $A:=K[Y_{ij}\mid 0\leq i\leq n,\: 0\leq j\leq n]_{\det([Y_{ij}\mid 0\leq i\leq n,\: 0\leq j\leq n])}$. For each $h=[h_{ij}\mid 0\leq i\leq n,\: 0\leq j\leq n]\in\Gl(n+\nolinebreak 1,K)$ define the maximal $A$-ideal $\idealp_h:=\langle\{Y_{ij}-h_{ij}\mid 0\leq i\leq n,\: 0\leq j\leq n\}\rangle_A$ and the morphism of $K$-algebras $\lambda_h:A\ra K$ with $\lambda_h(Y_{ij})=h_{ij}$ for all $0\leq i\leq n$ and $0\leq j\leq n$.

\emph{Claim:} $\zeta(\idealp_h)=\delta(h(\ideala))$ for all $h\in\Gl(n+\nolinebreak 1,K)$. 

Let $h\in\Gl(n+\nolinebreak 1,K)$. Since $\iota(\zeta(\idealp_h))=\zeta\circ\Spec(\lambda_h)$, we have to show that $\iota(\delta(h(\ideala)))=\zeta\circ\Spec(\lambda_h)$. In view of the commutative diagram
$$\xymatrix@C+0cm@R+5mm{
\Mor(\Gl^{n+1}_K,\Gl^{n+1}_K)\ar[r]^-{\raisebox{2mm}{$\scriptstyle\Mor(\Spec(\lambda_h),\,\Gl^{n+1}_K)$}}\ar[d]^{\cong}_{\psi_{\Gl}(\Gl^{n+1}_K)}&\Mor(\Spec(K),\Gl^{n+1}_K)\ar[r]^-{\raisebox{2mm}{$\scriptstyle\Mor(\Spec(K),\,\zeta)$}}\ar[d]^{\cong}_{\psi_{\Gl}(\Spec(K))}&\Mor(\Spec(K),\Hilb^p_K)\ar[d]_{\cong}^{\psi_{\Hilb}(\Spec(K))}\\
\Glf^{n+1}_K(\Gl^{n+1}_K)\ar[r]^-{\raisebox{2mm}{$\scriptstyle\Glf^{n+1}_K(\Spec(\lambda_h))$}}&\Glf^{n+1}_K(\Spec(K))\ar[r]^-{\underline{\zeta}(\Spec(K))}&\Hilbf^p_K(\Spec(K))
}$$
we verify indeed that
\begin{align*}
&\psi_{\Hilb}(\Spec(K))(\zeta\circ\Spec(\lambda_h))\\
&=(\psi_{\Hilb}(\Spec(K))\circ\Mor(\Spec(K),\,\zeta)\circ\Mor(\Spec(\lambda_h),\,\Gl^{n+1}_K))(\id_{\Gl^{n+1}_K})\\
&=(\underline{\zeta}(\Spec(K))\circ\Glf^{n+1}_K(\Spec(\lambda_h))\circ\psi_{\Gl}(\Gl^{n+1}_K))(\id_{\Gl^{n+1}_K})\\
&=\underline{\zeta}(\Spec(K))(\Glf^{n+1}_K(\Spec(\lambda_h))([Y_{ij}\mid 0\leq i\leq n,\: 0\leq j\leq n]))\\
&=\underline{\zeta}(\Spec(K))(h)=\sheaff^{(\delta(h(\ideala)))}=\psi_{\Hilb}(\Spec(K))(\iota(\delta(h(\ideala)))).
\end{align*}

Hence, the image of the closed points of $\Gl^{n+1}_K$ by the map $\zeta$ lies in $W$ and the points $\delta(\ideala)$ and $\delta(g(\ideala))$ lie in this image. Because $K$ is supposed to be algebraically closed, the set of closed points of $\Gl^{n+1}_K$ is homeomorphic to an irreducible algebraic variety (cf.\ \cite[II. 2.6]{H}), and the statement a) follows.

\sk
b) Choose $\omega\in\Z^{n+1}$ such that $\init_\tau\ideala=\init_\omega\ideala$ (cf.\ Proposition~\ref{existence of weight orders}). In the sequel we make use of the notations of \ref{notation Weight orders}. For $a\in K$ it holds 
$$\langle\beta^a_\omega(\ideala)\rangle_S=
\begin{cases}
\ideala,&\text{if $a=1$;}\\
{\init_\omega\ideala},&\text{if $a=0$.}
\end{cases}$$

Let $\iota:K\ra K[z]$ denote the canonical inclusion homomorphism and set  $\idealb:=\langle\beta^{\iota,z}_\omega(\ideala)\rangle_{K[z][X_0,\dots,X_n]}$. By Proposition \ref{a flat homomorphism}, the canonical homomorphism of rings $K[z]\ra K[z][X_0,\dots,X_n]/\idealb$ is flat. Hence, $\sheaff:=\sheafo_{\proj[{K[z]}]}/\widetilde\idealb$ is flat over $\Spec(K[z])$ (Lemma \ref{lemma about flatness}). By Proposition \ref{flatness and Hilbert polynomial} the Hilbert polynomial $p':=p_{\sheaff,y}\in\Q[t]$ of $\sheaff$ is constant in $y\in\Spec(K[z])$, whence $\sheaff\in\Hilbf^{p'}_K(\Spec(K[z]))$.

For any $K$-algebra $A$ and any $a\in A$ let $\phi^A_a:K[z]\ra A$ denote the homomorphism of $K$-algebras given by $\phi^A_a(z)=a$. For each Noetherian $K$-algebra $A$ we define the map $$\etaf(A):\afflinef_K(\Spec(A))\ra\Hilbf^{p'}_K(\Spec(A)),\;a\mapsto\Hilbf^{p'}_K(\Spec(\phi^A_a))(\sheaff).$$
Let $g:A\ra B$ be a homomorphism of Noetherian $K$-algebras and $a\in A$. Then we have $g\circ\phi^A_a=\phi^B_{g(a)}$. It follows that $$\Hilbf^{p'}_K(\Spec(g))\circ\Hilbf^{p'}_K(\Spec(\phi^A_a))=\Hilbf^{p'}_K(\Spec(\phi^B_{g(a)})),$$ and the diagram
$$\xymatrix{
\afflinef_K(\Spec(A))\ar[r]^-{\etaf(A)}\ar[d]_{\afflinef_K(\Spec(g))}&\Hilbf^{p'}_K(\Spec(A))\ar[d]^{\Hilbf^{p'}_K(\Spec(g))}\\
\afflinef_K(\Spec(B))\ar[r]^-{\etaf(B)}&\Hilbf^{p'}_K(\Spec(B))
}$$
commutes. Hence, $\etaf:\afflinef_K\circ\Spec\ra\Hilbf^{p'}_K\circ\Spec$ is a natural transformation, and by Remark \ref{remark about natural transformations} there is a natural transformation $\underline{\eta}:\afflinef_K\ra\Hilbf^{p'}_K$ such that $\underline{\eta}(\Spec(A))=\etaf(A)$ for any Noetherian $K$-algebra $A$. 

Let $a\in K$. By means of $\phi^K_a:K[z]\ra K$, $K$ is a $K[z]$-algebra. Using Lemma \ref{lemma on weighted ideals},  $(K[z][X_0,\dots,X_n]/\idealb)\otimes_{K[z]}K=S/\langle\beta^a_\omega(\ideala)\rangle_S$ as quotients of $S$. Bearing in mind Remark \ref{remark about computing inverse images}, we compute
\begin{align*}
\underline\eta(\Spec(K))(a)&=\etaf(K)(a)=\Hilbf^{p'}_K(\Spec(\phi^K_a))(\sheafo_{\proj[{K[z]}]}/\widetilde\idealb)\\
&=(K[z][X_0,\dots,X_n]/\idealb\otimes_{K[z]}K)\sptilde\\
&=(S/\langle\beta^a_\omega(\ideala)\rangle_S)\sptilde\\
&=\sheafo_{\proj[K]}/(\langle\beta^a_\omega(\ideala)\rangle_S)\sptilde\in\Hilbf^{p'}_K(\Spec(K)).
\end{align*}

Using the fact that $\langle\beta^1_\omega(\ideala)\rangle_S=\ideala$, we get in particular $p'=p$.

If $a\in K^*$ is a unit, then by Remark \ref{remark about weight orders} it holds $\sheaff^{(\delta(\ideala))}\cong\sheaff^{(\delta(\langle\beta^a_\omega(\ideala)\rangle_S))}$. Since $W$ is closed under isomorphisms and $\delta(\ideala)$, $\delta(\init_\tau\ideala)\in W$, it follows that $\delta(\langle\beta^a_\omega(\ideala)\rangle_S)\in W$ for all $a\in K$. 

Let $\eta:\affline_K\ra\Hilb_K^p$ be the unique morphism of $K$-schemes such that the diagram 
$$\xymatrix@C+2cm{
\Mor_{\sch_K}(\bullet,\,\affline_K)\ar[r]^-{\Mor_{\sch_K}(\bullet,\;\eta)}\ar[d]_{\psi_{\affline}}&\Mor_{\sch_K}(\bullet,\,\Hilb^p_K)\ar[d]^{\psi_{\Hilb}}\\
\afflinef_K\ar[r]^-{\underline{\eta}}&\Hilbf^p_K
}$$
commutes. Arguing as in the proof of part a) we see that the image of the connected set of closed points of $\afflinef_K$ under the map $\eta$ lies in $W$ and the points $x$ and $y$ lie in this image.
\end{proof}

\begin{corollary}\label{connectedness by connecting sequences}
Let $K$ be algebraically closed, $p\in\Q[t]$ a polynomial and $W\subset\closed(\Hilb^p_K)$ a topological subspace, closed under isomorphisms. If for all points $x$, $y\in W$ there is a connecting sequence in $\delta^{-1}(W)$ from $\delta^{-1}(x)$ to $\delta^{-1}(y)$, then $W$ is connected.\ \hfill$\square$
\end{corollary}
\sk
\subsection{Connected subschemes of the Hilbert scheme}\label{subsection main theorem}
We are now ready to define the Hilbert function strata 
$$\Hilb_K^{p,\fbar},\ \Hilb_K^{p,\fbar,\leq f},\ \Hilb_K^{f,\fbar},\ \barHilb_K^{p,\fbar},\ \barHilb_K^{p,\fbar,\leq f},\ \barHilb_K^{f,\fbar}$$ 
by bounding the cohomology functions of the points of $\Hilb^p_K$ (see the definitions in \ref{def subschemes}) and to prove the main Theorem~\ref{main Theorem} which states that 
$$\Hilb^{f,\fbar}_K,\ \barHilb^{p,\,\fbar}_K,\ \barHilb^{p,\fbar,\leq f}_K\text{ and }\barHilb^{f,\fbar}_K$$ 
are connected. It remains an open question whether $\Hilb_K^{p,\fbar}$ and $\Hilb_K^{p,\fbar,\leq f}$ are connected or not (s.~\ref{question}).

\begin{lemma}\label{shifting commutes with preimage}
Let $X$ be a scheme over $K$, $x\in X$ and $\sheaff$ a coherent sheaf of $\sheafo_{\proj[X]}$-modules. Then  $\sheaff(m)_{\kappa(x)}\cong\sheaff_{\kappa(x)}(m)$ for all $m\in\Z$.
\end{lemma}

\begin{proof}
Let $c:\Spec(\kappa(x))\ra X$ denote the canonical inclusion morphism. Let $U=\Spec(A)$ be an affine open neighbourhood of $x$ in $X$. Let $b:U\ra\nolinebreak X$ denote the inclusion. Hence, there is a morphism of $K$-algebras $f:A\ra\kappa(x)$ such that $c=b\circ\Spec(f)$. Let $\gamma:\proj[\kappa(x)]\ra\proj[X]$, $\beta:\proj[U]\ra\proj[X]$ and $\alpha:\proj[\kappa(x)]\ra\proj[A]$ denote the morphisms induced by $c$, $b$ and $\Spec(f)$, respectively. 

We first show that $\beta^*(\sheaff(m))\cong(\beta^*\sheaff)(m)$ for all $m\in\Z$. Since $\beta$ is open, we have $\beta^*(\sheafg\otimes_{\sheafo_{\proj[X]}}\sheafh)\cong\beta^*\sheafg\otimes_{\sheafo_{\proj[U]}}\beta^*\sheafh$ for all sheaves of $\sheafo_{\proj[X]}$-modules $\sheafg$ and $\sheafh$. Therefore, 
$\beta^*(\sheaff(m))=\beta^*(\sheaff\otimes_{\sheafo_{\proj[X]}}\sheafo_{\proj[X]}(m))\cong\beta^*\sheaff\otimes_{\sheafo_{\proj[U]}}\beta^*\sheafo_{\proj[X]}(m)
\cong\beta^*\sheaff\otimes_{\sheafo_{\proj[U]}}\sheafo_{\proj[U]}(m)=(\beta^*\sheaff)(m)$ for all $m\in\Z$.

Since $\beta^*\sheaff$ is a coherent sheaf on $\proj[A]$, there is a graded $A[X_0,\dots,X_n]$-module $M$ such that $\beta^*\sheaff=\widetilde M$. Let $m\in\Z$. By Remark \ref{remark about computing inverse images} it holds 
$\alpha^*((\beta^*\sheaff)(m))=\alpha^*(\widetilde M(m))\cong\widetilde{M\!\otimes_A\!\kappa(x)}(m)\cong(\alpha^*\widetilde M)(m)=(\alpha^*(\beta^*\sheaff))(m).$ 
It follows that $\sheaff(m)_{\kappa(x)}=\gamma^*(\sheaff(m))=\alpha^*(\beta^*(\sheaff(m)))\cong\alpha^*((\beta^*\sheaff)(m))\cong(\alpha^*(\beta^*\sheaff))(m)=(\gamma^*\sheaff)(m)=\sheaff_{\kappa(x)}(m).$
\end{proof}

\begin{notation and remark}\label{computation of cohomological Hilbert functions}
A) Let $p\in\Q[t]$ be a polynomial, $T$ a locally Noetherian scheme over $K$,\,\ $\sheaff\in\Hilbf^p_K(T)$, $y\in T$ and $i\in\N$. Let $\sheafj$ be the kernel of the induced quotient $\sheafo_{\proj[\kappa(y)]}\!\!\ra\sheaff_{\kappa(y)}$. Then $\sheafj$ is coherent. Denote by $I_{\sheaff,y}\subset\kappa(y)[X_0,\dots,X_n]$ the ideal $\bigoplus_{k\in\Z}H^0(\proj[\kappa(y)],\sheafj(k))$. Let $\sheafhilb\in\Hilbf^p_K(\Hilb^p_K)$ be the universal sheaf and $\sheafi\subset\sheafo_{\proj[\Hilb^p_K]}$\negthinspace the kernel of the quotient $\sheafo_{\proj[\Hilb^p_K]}\!\negthickspace\ra\sheafhilb$. Then $\sheafi$ is coherent.

Let $x\in\Hilb^p_K$. Using the notation of \ref{CM-regularity} B), define the following maps from $\Z$ to $\N$:
\begin{align*}
&h^i_{\sheaff,y}:=h^i_{\proj[\kappa(y)],\,\sheaff_{\kappa(y)}}&&h^i_x:=h^i_{\proj[\kappa(x)],\,\sheafhilb_{\kappa(x)}}\\
&\bar h^i_{\sheaff,y}:= h^i_{\proj[\kappa(y)],\,\sheafj}&&\bar h^i_x:=h^i_{\proj[\kappa(x)],\,\sheafi_{\kappa(x)}}\\
&h_{\sheaff,y}:=h_{\kappa(y)[X_0,\dots,X_n]/I_{\sheaff,y}}&&h_x:=h_{\sheafhilb,x}=h_{\kappa(x)[X_0,\dots,X_n]/I_{\sheafhilb,x}}
\end{align*}

\sk
B) Use the same notations as above. By Lemma \ref{shifting commutes with preimage}, $\sheafhilb_{\kappa(x)}(m)\cong\sheafhilb(m)_{\kappa(x)}$ and $\sheafi_{\kappa(x)}(m)\cong\sheafi(m)_{\kappa(x)}$. Therefore, 
\begin{align*}
&h^i_x(m)=\dim_{\kappa(x)}H^i(\proj[\kappa(x)],\sheafhilb(m)_{\kappa(x)}),\\
&\bar h^i_x(m)=\dim_{\kappa(x)}H^i(\proj[\kappa(x)],\sheafi(m)_{\kappa(x)})
\end{align*}
for all $m\in\Z$. Since $\sheafhilb$ is flat over $\Hilb^p_K$ and coherent on $\proj[\Hilb^p_K]$, the sequence
$$0\ra\sheafi_{\kappa(x)}\ra\sheafo_{\proj[\kappa(x)]}\ra\sheafhilb_{\kappa(x)}\ra 0$$
is exact, whence 
\begin{align*}
h_x(m)&=h_{\kappa(x)[X_0,\dots,X_n]/I_{\sheafhilb,x}}(m)=h_{\kappa(x)[X_0,\dots,X_n]}(m)-h_{I_{\sheafhilb,x}}(m)\\
&=\tbinom{m+n}{n}-h^0_{\proj[\kappa(x)],\,\sheafi_{\kappa(x)}}(m)\\
&=\tbinom{m+n}{n}-\dim_{\kappa(x)}H^0(\proj[\kappa(x)],\sheafi(m)_{\kappa(x)})
\end{align*}
for all $m\in\Z$.

\sk
C) Use the same notations as in A). Let $f:T\ra\Hilb^p_K$ be such that $\Hilbf^p_K(f)(\sheafhilb)=\sheaff$ (cf.~\ref{def - Hilbert scheme - def}), and let $f':\proj[T]\ra\proj[\Hilb^p_K]$ be the induced morphism. Then the sequence
$$0\ra((f')^*(\sheafi))_{\kappa(y)}\ra\sheafo_{\proj[\kappa(y)]}\ra\Hilbf^p_K(f)(\sheafhilb)_{\kappa(y)}\ra 0$$
is exact. Therefore, $$\sheafj=((f')^*(\sheafi))_{\kappa(y)}=\sheafi_{\kappa(f(y))},\ \;  I_{\sheaff,y}=I_{\sheafhilb,f(y)},\ \; h_{\sheaff,y}=h_{f(y)},\ \; \bar h^i_{\sheaff,y}=\bar h^i_{f(y)}.$$ 
In particular, if $g:T'\ra T$ is a morphism of locally Noetherian schemes over $K$ and $y'\in T'$, then  $$h_{\Hilbf^p_K(g)(\sheaff),y'}=h_{f(g(y'))}=h_{\sheaff,g(y')},\quad\bar h^i_{\Hilbf^p_K(g)(\sheaff),y'}=\bar h^i_{f(g(y'))}=\bar h^i_{\sheaff,g(y')}.$$
\end{notation and remark}

\begin{definition}
Let $Y$ be a topological space. A map $h:Y\ra\Z$ is called \emph{upper semicontinuous} if $\{y\in Y\mid h(y)\geq i\}$ is a closed subset of $Y$ for all $i\in\Z$.
\end{definition}

\begin{proposition}[{\cite[III. 12.8]{H}, Semicontinuity Theorem}]\label{upper semicontinuity}
Let $X\ra Y$ be a projective morphism of Noetherian schemes, and let $\sheaff$ be a coherent sheaf of $\sheafo_X$-modules which is flat over $Y$. Then for each $i\in\N$, the map 
$$Y\ra\Z,\;y\mapsto\dim_{\kappa(y)}H^i\bigl(X\times_Y\Spec(\kappa(y)),\sheaff_{\kappa(y)}\bigr)$$ 
is upper semicontinuous.\hfill$\qedsymbol$
\end{proposition}

\begin{corollary}\label{closed subsets of Hilb}
Let $p\in\Q[t]$ be a polynomial. Let $f:\Z\ra\N$ be a function and $i\in\N$. Then 
$$\{x\in\Hilb^p_K|h^i_x\geq f\},\ \,\{x\in\Hilb^p_K|\bar h^i_x\geq f\},\ \,\{x\in\Hilb^p_K|h_x\leq f\}$$ 
are closed subsets of\/ $\Hilb^p_K$ and $\{x\in\Hilb^p_K|h_x=f\}\subset\{x\in\Hilb^p_K|h_x\leq f\}$ is open.
\end{corollary}

\begin{proof}
Let $\sheafhilb\in\Hilbf^p_K(\Hilb^p_K)$ be the universal sheaf and $\sheafi\subset\sheafo_{\proj[\Hilb^p_K]}$ the kernel of the epimorphism $\sheafo_{\proj[\Hilb^p_K]}\!\negthickspace\ra\sheafhilb$. Now, $\proj[\Hilb^p_K]\!\ra\Hilb^p_K$ is a projective morphism of Noetherian schemes, and for each $m\in\Z$ the twisted sheaves $\sheafhilb(m)$ and $\sheafi(m)$ are coherent sheaves of $\sheafo_{\proj[\Hilb^p_K]}\!$-modules, flat over $\Hilb^p_K$ (Lemma \ref{lemma about flatness}). For $x\in\Hilb^p_K$ and $m\in\Z$ we have \begin{align*}
h^i_x(m)&=\dim_{\kappa(x)}H^i\bigl(\proj[\Hilb^p_K]\times_{\Hilb^p_K}\Spec(\kappa(x)),\sheafhilb(m)_{\kappa(x)}\bigr),\\
\bar h^i_x(m)&=\dim_{\kappa(x)}H^i\bigl(\proj[\Hilb^p_K]\times_{\Hilb^p_K}\Spec(\kappa(x)),\sheafi(m)_{\kappa(x)}\bigr),\\
h_x(m)&=\tbinom{n+m}{m}-\dim_{\kappa(x)}H^0\bigl(\proj[\Hilb^p_K]\times_{\Hilb^p_K}\Spec(\kappa(x)),\sheafi(m)_{\kappa(x)}\bigr).
\end{align*}
We conclude by Proposition \ref{upper semicontinuity} that the sets $\{x\in\Hilb^p_K|h^i_x(m)\geq f(m)\}$, $\{x\in\Hilb^p_K|\bar h^i_x(m)\geq f(m)\}$ and $\{x\in\Hilb^p_K|h_x(m)\leq f(m)\}$ are closed subsets of $\Hilb^p_K$. Hence 
\begin{align*}
\{x\in\Hilb^p_K|h^i_x\geq f\}&=\bigcap_{m\in\Z}\{x\in\Hilb^p_K|h^i_x(m)\geq f(m)\},\\
\{x\in\Hilb^p_K|\bar h^i_x\geq f\}&=\bigcap_{m\in\Z}\{x\in\Hilb^p_K|\bar h^i_x(m)\geq f(m)\},\\
\{x\in\Hilb^p_K|h_x\leq f\}&=\bigcap_{m\in\Z}\{x\in\Hilb^p_K|h_x(m)\leq f(m)\}
\end{align*}
are closed subsets of $\Hilb^p_K$.

Define the polynomial $q\in\Q[t]$ by $q(t):=\binom{t+n}{n}-p(t)$. Using the fact that the Hilbert polynomial $p_{\sheafi,x}=q$ is independent of the point $x\in\Hilb^p_K$, there is an integer $m_0\in\N$ such that $\sheafi_{\kappa(x)}$ is $m_0$-regular for all $x\in\Hilb^p_K$ by Proposition~\ref{admissible Hilbert polynomials}~b). It follows that $h_x(m)=\binom{m+n}{n}-q(m)$ for all $m\geq m_0$ and all $x\in\Hilb^p_K$ (cf.~\ref{Hilbert polynomial of sheaves}). Assume that the set $\{x\in\Hilb^p_K|h_x=f\}$ is not empty. Then $f(m)=0$ for all $m<0$ and $f(m)=p(m)$ for all $m\geq m_0$. Therefore,
$$\{x\in\Hilb^p_K\!|h_x\!=\!f\}=\{x\in\Hilb^p_K\!|h_x\!\leq\! f\}\setminus\bigcup_{m=0}^{m_0-1}\{x\in\Hilb^p_K\!|h_x(m)\!<\!f(m)\}$$
is an open subset of $\{x\in\Hilb^p_K|h_x\leq f\}$.
\end{proof}

\begin{remark}\label{factorising through subschemes}
A) \cite[4.6.1]{EGA I} Let $X$ be a scheme. For any locally closed subset $Y$ of $X$ there is a unique reduced subscheme $Y_\red$ of $X$ such that the underlying topological space of $Y_\red$ is $Y$. 

\sk
B) Let $X$ be a reduced scheme, $f:X\ra Y$ a morphism of schemes and $Z$ a locally closed subscheme of $Y$ with canonical morphism $j:Z\ra Y$. If $f(x)\in Z$ for all $x\in X$, then there exists a factorising morphism $f':X\ra Z$ such that $j\circ f'=f$. (If $Z$ is closed, this is proved in \cite[4.6.2]{EGA I}, if $Z$ is open, then the statement holds even if $X$ is not reduced (cf.\ \cite[4.2.2 and 4.1.6]{EGA I})).

\sk
C) As a consequence, if $X$ and $Y$ are reduced schemes over a scheme $S$, then $(X\times_S Y)_\red$ is the fibred product of $X$ and $Y$ over $S$ in the category of reduced schemes.
\end{remark}

\begin{definition and remark}\label{def subschemes}
Let $\sch_K^\red$ denote the category of reduced locally Noetherian schemes over $K$. Let $p\in\Q[t]$ be a polynomial, $f:\Z\ra\N$ a function and $\fbar=(f_i)_{i\in\N}$ a sequence of functions $f_i:\Z\ra\N$. 

\sk
A) If $T$ is a locally Noetherian scheme over $K$, set
\begin{align*}
&\Hilbf^{p,\fbar}_K(T):=\{\sheaff\in\Hilbf^p_K(T)\mid h^i_{\sheaff,x}\geq f_i\ \forall i\in\N\ \forall x\in T\},\\
&\Hilbf^{p,\fbar,\leq f}_K(T):=\{\sheaff\in\Hilbf^{p,\fbar}_K(T)\mid h_{\sheaff\!,\,x}\leq f\ \forall x\in T\},\\
&\Hilbf^{f,\fbar}_K(T):=\{\sheaff\in\Hilbf^{p,\fbar}_K(T)\mid h_{\sheaff\!,\,x}=f\ \forall x\in T\},\\
&\barHilbf^{p,\fbar}_K(T):=\{\sheaff\in\Hilbf^p_K(T)\mid\bar h^i_{\sheaff,x}\geq f_i\ \forall i\geq 1\ \forall x\in T\},\\
&\barHilbf^{p,\fbar,\leq f}_K(T):=\{\sheaff\in\barHilbf^{p,\fbar}_K(T)\mid h_{\sheaff\!,\,x}\leq f\ \forall x\in T\},\\
&\barHilbf^{f,\fbar}_K(T):=\{\sheaff\in\barHilbf^{p,\fbar}_K(T)\mid h_{\sheaff\!,\,x}=f\ \forall x\in T\}.
\end{align*}
If $g:T'\ra T$ is a morphism of locally Noetherian schemes over $K$, $\sheaff\in\Hilbf^p_K(T)$ and $x\in T'$, then $\kappa(x)$ is a field extension of $\kappa(g(x))$, and \linebreak[4] $h^i_{\proj[\kappa(x)],\,(\Hilbf^p_K(g)(\sheaff))_{\kappa(x)}}=h^i_{\proj[\kappa(g(x))],\,\sheaff_{\kappa(g(x))}}$ for all $i\in\N$ (cf.\ Remark \ref{consequences of Serre-Grothendieck Correspondence}~C)). Hence, if $\sheaff\in\Hilbf^{p,\fbar}_K(T)$, then  $\Hilbf^p_K(g)(\sheaff)\in\Hilbf^{p,\,\fbar}_K(T')$. Therefore, there is a map 
$$\Hilbf^{p,\fbar}_K(g):\Hilbf^{p,\fbar}_K(T)\ra\Hilbf^{p,\fbar}_K(T'),\;\sheaff\mapsto\Hilbf^p_K(g)(\sheaff).$$
Thus, we have defined a contravariant subfunctor $$\Hilbf_K^{p,\fbar}:\sch_K\ra\sets$$ of the Hilbert functor $\Hilbf^p_K$.

Moreover, if $g:T'\ra T$ is a morphism of locally Noetherian schemes over $K$, $\sheaff\in\Hilbf^{p,\fbar}_K(T)$ and $x\in T'$, then  $h_{\Hilbf^p_K(g)(\sheaff),x}=h_{\sheaff,g(x)}$ (cf.~\ref{computation of cohomological Hilbert functions}~C)). Therefore, we may define the contravariant subfunctors 
$$\Hilbf^{p,\fbar,\leq f}_K:\sch_K\ra\sets,\qquad\Hilbf^{f,\fbar}_K:\sch_K\ra\sets$$
of $\Hilbf^{p,\fbar}_K$ in the same way.
Similarly, there are contravariant subfunctors 
$$\barHilbf_K^{p,\fbar}:\sch_K\ra\sets,\quad\barHilbf^{p,\fbar,\leq f}_K:\sch_K\ra\sets,\quad\barHilbf^{f,\fbar}_K:\sch_K\ra\sets.$$

Restricting these functors to the full subcategory of reduced locally Noe\-therian schemes over $K$ we get functors
\begin{align*}
&\Hilbf_K^{p,\fbar}:\sch_K^\red\ra\sets,&&\barHilbf_K^{p,\fbar}:\sch_K^\red\ra\sets,\\
&\Hilbf^{p,\fbar,\leq f}_K:\sch_K^\red\ra\sets,&&\barHilbf^{p,\fbar,\leq f}_K:\sch_K^\red\ra\sets,\\
&\Hilbf^{f,\fbar}_K:\sch_K^\red\ra\sets,&&\barHilbf^{f,\fbar}_K:\sch_K^\red\ra\sets.
\end{align*}

\sk
B) Define the sets
\begin{align*}
&\Hilb^{p,\fbar}_K:=\{x\in\Hilb^p_K\mid h^i_x\geq f_i\ \forall i\in\N\},\\
&\Hilb^{p,\fbar,\leq f}_K:=\{x\in\Hilb^{p,\fbar}_K\mid h_x\leq f\},\\
&\Hilb^{f,\fbar}_K:=\{x\in\Hilb^{p,\fbar}_K\mid h_x=f\},\\
&\barHilb^{p,\fbar}_K:=\{x\in\Hilb^p_K\mid\bar h^i_x\geq f_i\ \forall i\geq 1\},\\
&\barHilb^{p,\fbar,\leq f}_K:=\{x\in\barHilb^{p,\fbar}_K\mid h_x\leq f\},\\
&\barHilb^{f,\fbar}_K:=\{x\in\barHilb^{p,\fbar}_K\mid h_x=f\}
\end{align*}
and endow them with the induced topology. Then, by Corollary \ref{closed subsets of Hilb}, 
$$\Hilb^{p,\fbar}_K=\bigcap_{i\in\N}\{x\in\Hilb^p_K\mid h^i_x\geq f_i\}$$ 
is a closed subset of $\Hilb^p_K$, $\Hilb^{p,\fbar,\leq f}_K$ is a closed subset of $\Hilb^{p,\fbar}_K$ and $\Hilb^{f,\fbar}_K$ is an open subset of $\Hilb^{p,\fbar,\leq f}_K$. We endow these sets with the reduced induced scheme structure. Similarly, there are closed reduced subschemes $\barHilb^{p,\fbar}_K$, $\barHilb^{p,\fbar,\leq f}_K$ and a locally closed reduced subscheme $\barHilb^{f,\fbar}_K$ of the Hilbert scheme.
\end{definition and remark}

\begin{proposition}\label{the functor Hilb^f is representable}
Let $p\in\Q[t]$ be a polynomial, $f:\Z\ra\N$ a function and $\fbar=(f_i)_{i\in\N}$ a sequence of functions $f_i:\Z\ra\N$. Then, the functors
$$\Hilbf_K^{p,\fbar},\ \Hilbf_K^{p,\fbar,\leq f},\ \Hilbf_K^{f,\fbar},\ \barHilbf_K^{p,\fbar},\ \barHilbf_K^{p,\fbar,\leq f},\ \barHilbf_K^{f,\fbar}:\sch_K^\red\ra\sets$$ 
are representable by the reduced schemes 
$$\Hilb^{p,\fbar}_K,\  \Hilb^{p,\fbar\!,\leq f}_K,\ \Hilb^{f\!,\fbar}_K,\  \barHilb^{p,\fbar}_K,\  \barHilb^{p,\fbar\!,\leq f}_K,\  \barHilb^{f\!,\fbar}_K,$$ 
respectively.
\end{proposition}

\begin{proof}
Let $\sheafhilb\in\Hilbf^p_K(\Hilb^p_K)$ be the universal sheaf and $j:\Hilb^{p,\fbar}_K\ra\Hilb^p_K$ the inclusion morphism. If $g:T\ra\Hilb^p_K$ is a morphism of locally Noetherian schemes over $K$, and if $x\in T$, then $\kappa(x)$ is a field extension of $\kappa(g(x))$, whence 
$h^i_{\proj[\kappa(x)],\,\Hilbf^p_K(g)(\sheafhilb)_{\kappa(x)}}=h^i_{\proj[\kappa(g(x))],\,\sheafhilb_{\kappa(g(x))}}=h^i_{g(x)}.$ 
It follows that $\Hilbf^p_K(g)(\sheafhilb)\in\Hilbf^{p,\fbar}_K(T)$ if and only if $g(x)\in\Hilb^{p,\fbar}_K$ for all $x\in T$. If $T$ is reduced, then by Remark \ref{factorising through subschemes} B),  $\Hilbf^p_K(g)(\sheafhilb)\in\Hilbf^{p,\fbar}_K(T)$ if and only if there exists a factorising morphism of schemes $g':T\ra\Hilb^{p,\fbar}_K$ such that $g=j\circ g'$.

Define the natural transformation 
$$\psi_{\Hilb^\fbar}:\Mor_{\sch_K}(\bullet,\,\Hilb^{p,\fbar}_K)\ra\Hilbf^{p,\fbar}_K$$ 
by 
$$\psi_{\Hilb^\fbar}(T):\Mor_{\sch_K}(T,\,\Hilb^{p,\fbar}_K)\ra\Hilbf^{p,\fbar}_K(T),\;g\mapsto\Hilbf^{p,\fbar}_K(j\circ g)(\sheafhilb)$$ 
for each locally Noetherian $K$-scheme $T$.

Let $T$ be a locally Noetherian scheme over $K$. Then there is a commutative diagram
$$\xymatrix@C+1cm{
\Mor_{\sch_K}(T,\,\Hilb^p_K)\ar[r]^-{\psi_{\Hilb}(T)}_-{\cong}&\Hilbf^p_K(T)\\
\Mor_{\sch_K}(T,\,\Hilb^{p,\fbar}_K)\ar[r]^-{\psi_{\Hilb^\fbar}(T)}\ar@{ >->}[u]_{\Mor_{\sch_K}(T,\,j)}&\Hilbf^{p,\fbar}_K(T)\ar@<3mm>@{ (->}[u]
.}$$
Now assume that $T$ is reduced. Let $\sheaff\in\Hilbf^{p,\fbar}_K(T)$. Let $g:T\ra\Hilb^p_K$ be the unique morphism such that $\psi_{\Hilb}(T)(g)=\sheaff$. Since $\Hilbf^p_K(g)(\sheafhilb)=\psi_{\Hilb}(T)(g)=\sheaff\in\Hilbf^{p,\fbar}_K(T)$, $g$ factorises through $j$. In this way we have proven that $\psi_{\Hilb^\fbar}(T)$ is bijective. It follows that the functor\/ $\Hilbf_K^{p,\fbar}\!:\sch_K^\red\ra\sets$ is represented by the pair $(\Hilb^{p,\fbar}_K,\ \Hilbf^p_K(j)(\sheafhilb))$.

Moreover, we have $h_{\sheaff,x}=h_{g(x)}$ for all $x\in T$ by \ref{computation of cohomological Hilbert functions}~C). Hence, if $T$ is reduced, then $g:T\ra\Hilb^{p,\fbar}_K$ factorises through $\Hilb_K^{p,\fbar,\leq f}$ or $\Hilb_K^{f,\fbar}$ if and only if $\sheaff=\Hilbf^{p,\fbar}_K(g)(\sheafhilb)$ is an element of $\Hilbf^{p,\fbar,\leq f}_K(T)$ or $\Hilbf^{f,\fbar}_K(T)$, respectively. Thus, we may define representing natural equivalences 
\begin{align*}
&\psi_{\Hilb^{\fbar,\leq f}}:\Mor_{\sch_K^\red}(\bullet,\,\Hilb^{p,\fbar,\leq f}_K)\ra\Hilbf^{p,\fbar,\leq f}_K,\\
&\psi_{\Hilb^{\fbar,=f}}:\Mor_{\sch_K^\red}(\bullet,\,\Hilb^{f,\fbar}_K)\ra\Hilbf^{f,\fbar}_K
\end{align*}
as above.

The remaining statements are proved similarly using the last equation of \ref{computation of cohomological Hilbert functions}~C).
\end{proof}

\begin{remark}\label{remark about representable functors}
A) Let $\cat$ be a category with fibred products. For an object $Y$ of $\cat$ let $\cat_Y$ denote the category of objects and morphisms over $Y$. Let $Y'\ra Y$ be a morphism of $\cat$ turning $\cat_{Y'}$ into a subcategory of $\cat_Y$. Let $F:\nolinebreak\cat_Y\ra\sets$ be a contravariant functor, representable by an $Y$-object $X$. Let $F\restr_{\cat_{Y'}}:\cat_{Y'}\ra\sets$ denote the restricted functor, defined by $F\restr_{\cat_{Y'}}\negmedspace(T)=F(T)$ for all objects $T$ of $\cat_{Y'}$. Then $F\restr_{\cat_{Y'}}$ is representable by $X\times_Y Y'$ (\cite[0.1.3.10]{EGA I}).

\sk
B) Let $p\in\Q[t]$ be a polynomial, $f:\Z\ra\N$ a function and $\fbar=(f_i)_{i\in\N}$ a sequence of functions $f_i:\Z\ra\N$. Let $X_K$ be one of the six reduced subschemes of $\Hilb^p_K$ defined in \ref{def subschemes}~B). Let $K\subset k$ be a field extension. Then, by Remark \ref{factorising through subschemes}~C), it holds 
$X_k\cong(X_K\times_{\Spec(K)}\Spec(k))_\red.$
\end{remark}

We now put together all previous results and prove our main Theorem:
\begin{theorem}\label{main Theorem}
Assume that the field $K$ is of characteristic zero. Let $p\in\Q[t]$ be a polynomial, $f:\Z\ra\N$ a function and $\fbar=(f_i)_{i\in\N}$ a sequence of functions $f_i:\Z\ra\N$. Then the topological spaces $\Hilb^{f,\fbar}_K$\!, $\barHilb^{p,\,\fbar}_K$\!, $\barHilb^{p,\fbar,\leq f}_K$\!, and\/ $\barHilb^{f,\fbar}_K$ are connected.
\end{theorem}

\begin{proof}
Observe that $\barHilb^{p,\,\fbar}_K=\barHilb^{p,\fbar,\leq h_S}_K$. So it suffices to prove our statement for each of the reduced schemes $\Hilb^{f,\fbar}_K,\ \barHilb^{p,\fbar,\leq f}_K,\ \barHilb^{f,\fbar}_K$. Let $X_K\in\{\Hilb^{f,\fbar}_K,\ \barHilb^{p,\fbar,\leq f}_K,\ \barHilb^{f,\fbar}_K\}$. Denote the corresponding functor $\Hilbf^{f,\fbar}_K$, $\barHilbf^{p,\fbar,\leq f}_K$ or $\barHilbf^{f,\fbar}_K$ by $\Xf_K$.

If $k$ is the algebraic closure of $K$, then by Remark \ref{remark about representable functors} it holds $X_k\cong(X_K\times_{\Spec(K)}\Spec(k))_\red$. If $X_K$ were not connected, then also $X_k$ were not connected. Hence, we may assume that $K$ is algebraically closed.

In a scheme, any non-empty closed subset contains a closed point. Therefore, it suffices to show that $\closed(X_K)$ is connected.

In what follows we use the notations of \ref{notations for points in Hilb}. Since $K$ is algebraically closed, $\closed(X_K)$ is in bijection with the set $\Mor_{\sch_K}(\Spec(K),\,X_K)$, whence $\closed(X_K)$ is a subset of $\closed(\Hilb^p_K)$. Notice that $\closed(X_K)$ is closed under isomorphisms.

We claim that \emph{for any two points $x$, $y\in W:=\closed(X_K)$ there is a connecting sequence from $\delta^{-1}(\sheaff^{(x)})$ to $\delta^{-1}(\sheaff^{(y)})$ in $\delta^{-1}(W)$.} Then our Theorem follows from Corollary \ref{connectedness by connecting sequences}.

We assume that $0\notin\Ideals_p$, because otherwise $\Hilb^p_K$ has only one closed point. Let $\ideala\in\Ideals_p$. It follows from the Serre-Grothendieck Correspondence (Proposition \ref{Serre-Grothendieck}) that 
\begin{align*}
&h^i_{\proj[K]\!,\,\sheafo_{\proj}/\widetilde\ideala}=h^{i+1}_{S/\ideala}\text{ for all }i\geq 1,\\ &h^0_{\proj[K]\!,\,\sheafo_{\proj}/\widetilde\ideala}=h^1_{S/\ideala}+h_{S/\ideala^\sat},\\ &h^i_{\proj[K]\!,\,\widetilde\ideala}=h^i_{S/\ideala}\text{ for all }1\leq i<n,\\ &h^n_{\proj[K]\!,\,\widetilde\ideala}=h^n_{S/\ideala}+h^{n+1}_S
\end{align*}
(cf.~Remark \ref{consequences of Serre-Grothendieck Correspondence}~B)).
Hence, the following equalities hold:
$$\delta^{-1}(W)=\{\ideala\in\Ideals_p\mid h_{S/\ideala^\sat}=f,\: h^1_{S/\ideala}\geq f_0-f,\: h^i_{S/\ideala}\geq f_{i-1}\ \forall\ i\geq 2\}$$
if $X_K=\Hilb^{f,\fbar}_K$,
$$\delta^{-1}(W)=\{\ideala\in\Ideals_p\mid h_{S/\ideala^\sat}\leq f,\: h^n_{S/\ideala}\geq f_n-h^{n+1}_S,\: h^i_{S/\ideala}\geq f_i\ \forall\ 1\leq i<n\}$$
if $X_K=\barHilb^{p,\fbar,\leq f}_K$,
$$\delta^{-1}(W)=\{\ideala\in\Ideals_p\mid h_{S/\ideala^\sat}=f,\: h^n_{S/\ideala}\geq f_n-h^{n+1}_S,\: h^i_{S/\ideala}\geq f_i\ \forall\ 1\leq i<n\}$$
if $X_K=\barHilb^{f,\fbar}_K$.

The claim now follows from Corollary \ref{corollary of connection sequences (= f)} and Theorem \ref{theorem of connection sequences (<= f)}.
\end{proof}

\begin{open question}\label{question}
It is still an open problem whether the space $\Hilb^{p,\fbar,\leq f}_K$ is connected or not. To formulate the problem, we use the notations of the previous proof: Let $W:=\closed(\Hilb^{p,\fbar,\leq f}_K)$. If $\ideala$, $\idealb\in\Ideals_p$ are two consecutive elements of a connecting sequence in $\delta^{-1}(W)$ coming from the algorithm described in Proposition \ref{existence of Mall bin syst to a lexicographic segment}, then $h^1_{S/\ideala}\leq h^1_{S/\idealb}$ and $h_{S/\ideala^\sat}\geq h_{S/\idealb^\sat}$ (or vice versa). But it is not a priori clear that  $$h^0_{\proj[K]\!,\,\sheafo_{\proj}/\widetilde\ideala}=h^1_{S/\ideala}+h_{S/\ideala^\sat}\leq h^1_{S/\idealb}+h_{S/\idealb^\sat}=h^0_{\proj[K]\!,\,\sheafo_{\proj}/\widetilde\idealb}.$$

The evidence of hundreds of examples up to dimension $n=5$ suggests that $\Hilb^{p,\fbar,\leq f}_K$ is indeed connected. But a proof seems to request complicated combinatorial considerations.
\end{open question}
\clearpage
\addcontentsline{toc}{section}{References}
\bibliographystyle{plain}

\end{document}